\documentclass[12pt,sort&compress]{elsarticle}

\usepackage{afterpage}
\usepackage{enumitem}
\usepackage{mystyle}
\usepackage{algpseudocode}
\usepackage{algorithm}
\usepackage{amsmath}
\usepackage{amssymb}
\usepackage{subfiles}
\graphicspath{{Images/}{../Images/}}
\usepackage{subfigmat} 
\usepackage{relsize}
\usepackage{pst-all}
\usepackage[normalem]{ulem}
\usepackage{cancel}
\usepackage{caption}
\usepackage{cleveref}
\usepackage{comment}

\newcommand\figref[1]{Figure \ref{fig:#1}} 
\newcommand\tabref[1]{Table \ref{tab:#1}} 
\newcommand\eref[1]{Eq. (\ref{eq:#1})} 
\newcommand\p[1]{\partial{#1}}

\begin{document}

\begin{frontmatter}


\title{Data-driven identification of stable differential operators using constrained regression}

\author[CMU]{Aviral Prakash\corref{cor1}}
\ead{aviralp@andrew.cmu.edu}
\author[CMU]{Yongjie Jessica Zhang}
\cortext[cor1]{Corresponding author}

\address[CMU]{Department of Mechanical Engineering, Carnegie Mellon University, Pittsburgh, PA 15213, USA}

\begin{abstract}
    Identifying differential operators from data is essential for the mathematical modeling of complex physical and biological systems where massive datasets are available. These operators must be stable for accurate predictions for dynamics forecasting problems. In this article, we propose a novel methodology for learning sparse differential operators that are theoretically linearly stable by solving a constrained regression problem. These underlying constraints are obtained following linear stability for dynamical systems. We further extend this approach for learning nonlinear differential operators by determining linear stability constraints for linearized equations around an equilibrium point. The applicability of the proposed method is demonstrated for both linear and nonlinear partial differential equations such as 1-D scalar advection-diffusion equation, 1-D Burgers equation and 2-D advection equation. The results indicated that solutions to constrained regression problems with linear stability constraints provide accurate and linearly stable sparse differential operators.
    
\end{abstract}


\end{frontmatter}

\section{Introduction}

Mathematical models for predicting physical systems rely on well-defined partial differential equations (PDEs) or ordinary differential equations (ODEs) that govern the spatio-temporal dynamics. Traditionally, these equations have been defined based on physical insights from domain experts. This approach has prohibited simulations for many complex systems, such as those in biology and finance, where obtaining such equations through insights may not be possible. With advancements in machine learning techniques and the availability of large-scale datasets, there has been a widespread interest in data-driven simulation and modeling \cite{Oden2006, Oden2007, Bajaj2007, Brunton2019, Prakash2022, Prakash2023a}. This interest has led to mathematical techniques that utilize such datasets with partial physical information to infer physical systems and improve existing models. 

For several scenarios, such as physical and biological dynamical systems where PDEs and ODEs are not readily available, there has been growing interest in techniques for system identification \citep{Brunton2016, Raissi2018, Zhang2018}. The first step of modeling complex systems is to determine PDEs from data \cite{Rudy2017, Maddu2022}. The next step is solving these learned PDEs to predict system dynamics. The typical approach to solving these identified equations is using numerical methods such as finite differences \cite{LeVeque2007}, finite volumes \cite{LeVeque2002} or finite elements \cite{Hughes2000}. These numerical methods have been commonly used for solving PDEs as these have provable theoretical estimates of accuracy and stability. However, these methods are often tailored to maintain accuracy and stability for specific problems. Furthermore, the solution of PDEs using these methods also requires the development of a comprehensive and efficient codebase, especially if the targeted application is complex and computationally expensive. When system identification techniques \cite{Rudy2017, Maddu2022} are applied for identifying PDEs, accuracy and stability considerations for the spatial and temporal discretization of learned PDEs to obtain a valid solution are typically not addressed. Instead, appropriate experience and insights from users are expected to address this selection. In such situations, it is desirable to determine discrete differential operators, hereafter referred to as differential operators, directly using data and modern machine learning methods. Furthermore, these differential operators can also enable nonintrusive reduced order modeling \cite{Gkimisis2024} where, despite the knowledge of underlying PDE, the knowledge and access to the discretization scheme are unavailable. 

Identifying appropriate differential operators from data, which is intimately related to learning spatial discretization from data \cite{Bar-Sinai2019} has gained significant interest over the years. Such methods can also be classified as techniques for system identification as they identify the set of discretized equations from data. Several popular approaches use artificial neural networks (ANNs) for discrete representations of systems \cite{Bar-Sinai2019, Long2018, Maddu2023}. Despite their popularity and accuracy for different applications, the typical black-box nature of ANNs often discourages interpretability \cite{Long2019}, which is desired when performing theoretical analysis to identify the accuracy and stability properties of the method. An in-depth study of these properties is essential for gauging the performance of learned differential operators for scenarios not included in the training dataset. Furthermore, these nonlinear learned differential operators could have a high evaluation cost, prohibiting the scalability of such methods for large-scale physical systems. Standard numerical methods often result in sparse linear systems with a lower computational overhead. Decades of research in solving such systems have resulted in sparse linear system solvers that efficiently provide high accuracy with a low memory footprint and enable scalability and portability to different hardware architectures. These ANN-based discrete operators may not allow the efficient use of such sparse linear system solvers, which could require implicit time integration for stiffer dynamics. 

In contrast to previous approaches, which learn noninterpretable differential operators, recent work \cite{Schumann2022} has focused on obtaining interpretable sparse stencils from data for linear PDEs. They posed this as a regression problem where a local solution stencil is learned from the data. In a follow-up work \cite{Schumann2023}, this approach was extended for time-dependent problems and nonlinear PDEs while using a strategy \cite{Maddu2022} to identify stencils and regularization parameters for stable learned differential operators. Recent work on the adjacency-based determination of differential operators \cite{Gkimisis2023, Gkimisis2024} has focused on a similar approach for identifying sparse differential operators and demonstrating their application for nonintrusive reduced order modeling. The accuracy of these methods can be adjusted by selecting the appropriate solution stencil sizes, machine learning techniques and regularization parameters. However, despite high accuracy within the training dataset, the numerical stability of these discrete solution representations is essential to ensure that numerical approximation errors do not grow in time. All these methods, including ANN-based methods, do not theoretically guarantee stability even for linear systems. Instead, feasible stable solutions may only be constrained to scenarios within the validation dataset without any stability guarantees for dynamics forecasting scenarios outside the validation dataset. For the reasons mentioned above, there is an immense need for approaches that determine differential operators from data while providing theoretical stability guarantees. Our experience with learning differential operators through state-of-the-art approaches that use regression to obtain differential operators \citep{Schumann2022, Gkimisis2024} indicated that these methods often learn unstable differential operators for linear PDEs and may not reliably perform even within the training dataset. While this issue is partly addressed in \cite{Schumann2023} using the stability selection procedure \cite{Meinshausen2010}, this approach can be computationally expensive for a system with many degrees of freedom and cannot theoretically guarantee stability. 

In this article, we propose a novel approach for learning sparse differential operators that are provably linearly stable. This approach relies on a set of local conditions for differential operators derived using the stability theory for linear dynamical systems. These conditions are incorporated as inequality constraints in the regression problem to determine the unknown differential operators. The resulting constrained regression problem is solved using a sequential least squares programming optimizer. We further extend this method for learning nonlinear differential operators by formulating constraints based on linearized equations obtained using Taylor-series expansion around an equilibrium point. The applicability of the proposed method for learning stable differential operators is demonstrated by comparing the results against the standard regression-based approach for multiple linear and nonlinear PDEs: 1-D scalar advection-diffusion equation, 1-D Burgers equation and 2-D advection equation. The proposed approach targets identifying suitable differential operators while using the known form of PDE to provide stability constraints. Therefore, this approach differs from other system identification approaches that identify an unknown PDE. 

The outline of this article is given below. Section 2 discusses the relevant mathematical background on differential operators for PDEs and stability theory for linear differential operators. Section 3 first discusses the standard approach for learning sparse differential operators from data. This section also gives the mathematical details for learning stable sparse differential operators from data and an extension of this approach to nonlinear equations. Section 4 includes the results for the three test cases and demonstrates the applicability of the proposed approach for learning differential operators from data. Section 5 concludes this article by highlighting the main contributions and mentioning directions for future research. 

\section{Mathematical background}

In this section, we first introduce the theory of obtaining the semi-discrete form of PDEs and then discuss the linear stability of differential operators in this semi-discrete form. 

\subsection{Differential operators for PDEs}

We restrict the analysis and mathematical formulation to 1-D PDEs, although these concepts can be generalized to higher-dimensional PDEs. Consider a 1-D PDE of the following form:
\begin{equation}
    \frac{\partial u}{\p t} + \mathcal{F} (u) = 0,
    \label{eq:PDEeq1}
\end{equation}
where $\mathcal{F} (\cdot)$ is a continuous differential operator, $u: \Omega \times [0,T] \to \mathbb{R}$ is the PDE solution, $\Omega$ is the simulation domain and $T$ is the final simulation time. Boundary conditions often accompany these PDEs on the boundaries $\Gamma$. Using the method of lines with a suitable spatial discretization, the semi-discrete form of this PDE is obtained as
\begin{equation}
    \frac{d \bm{u}}{d t} + \bm{F} (\bm{u}) = 0,
    \label{eq:PDEeq2}
\end{equation}
where $\bm{u} \in \mathbb{R}^n$ is the discrete solution field, $\bm{F}: \mathbb{R}^n \to \mathbb{R}^n$ is a differential operator and $n$ is the number of degrees of freedom. 
These degrees of freedom often correspond to different locations on the discretized domain, called the simulation grid $\Omega^h$, depending on the spatial discretization approach. The differential operator can be linear or nonlinear, depending on the PDE under consideration. We make this distinction by decomposing the differential operator into a linear component $\bm{L}$ and a nonlinear component $\bm{\hat{N}}$ such that \eref{PDEeq2} is simplified as
\begin{equation}
    \frac{d \bm{u}}{d t} + \bm{L} \bm{u} + \bm{\hat{N}} (\bm{u}) = 0.
    \label{eq:PDEeq3global}
\end{equation}
The form of these linear and nonlinear operators depends on the numerical method chosen to discretize the PDE. Typical discretization methods such as finite difference, finite volume or finite element methods result in a sparse operator, implying that the linear and nonlinear operators are applied to localized degrees of freedom, commonly called solution stencil. The sparsity of these operators allows efficient storage and faster computation of the PDE solution field. The discretized PDE at the $i^{th}$ degree of freedom can be considered an ODE problem
\begin{equation}
    \frac{d u_i}{d t} + (\bm{L}^i)^T \bm{u}_{\Omega^l_i} + \bm{\hat{N}}^i (\bm{u}_{\Omega^n_i}) = 0,
    \label{eq:PDEeq3}
\end{equation}
\noindent where $\bm{L}^i \in \mathbb{R}^{s_l}$ is the local linear operator which contributes to the $i^{th}$ row of $\bm{L}$, superscript $T$ indicates the transpose, $s_l$ is the size of the local linear stencil, $\bm{\hat{N}}^i: \mathbb{R}^{s_n} \to \mathbb{R}$ is the local nonlinear operator and $s_n$ is the size of the local nonlinear stencil. The solution stencils $\bm{u}_{\Omega^l_i}$ and $\bm{u}_{\Omega^n_i}$ are selected based on local degrees of freedom. In the context of this article, we transform the nonlinear term in \eref{PDEeq3} to a matrix-vector product, which is similar to the linear term. The resulting equation is
\begin{equation}
   \frac{d u_i}{d t} + (\bm{L}^i)^T \bm{u}_{\Omega^l_i} + (\bm{N}^i)^T \bm{z} (\bm{u}_{\Omega^n_i}) = 0,
    \label{eq:PDEeq4}
\end{equation}
where $\bm{N}^i \in \mathbb{R}^{s_n}$ is the nonlinear differential operator in a vector form, $s_n$ is the nonlinear stencil size and $\bm{z} (\bm{u}_{\Omega^n_i}) \in \mathbb{R}^{s_n}$ represents nonlinear products of $\bm{u}_{\Omega^n_i}$. To elucidate this notation, we consider an example of the 1-D viscous Burgers equation
\begin{equation}
    \frac{\p u}{\p t} + u \frac{\p u}{\p x} - \nu \frac{\partial^2 u}{\partial x^2} = 0
    \label{eq:BurgerEq}
\end{equation}
with periodic boundary conditions. For this example, we select a $1^{st}$-order backward difference for the nonlinear term and a $2^{nd}$-order centered finite difference for the linear term for discretizing the system on $n$ uniformly spaced grid nodes. This discretization gives the following semi-discrete form
\begin{equation}
    \frac{d u_i}{d t} + u_i \frac{u_{i} - u_{i-1}}{\Delta x} - \nu \frac{u_{i+1} - 2 u_{i} + u_{i-1}}{(\Delta x)^2} = 0.
\end{equation}
This semi-discrete form can be written as \eref{PDEeq4} with
\begin{equation}
    \bm{N}^i = \frac{1}{\Delta x}[-1, 1, 0]^T \quad \text{and} \quad \bm{L}^i = \frac{1}{(\Delta x)^2}[-1, 0, 1]^T \;,
\end{equation}
\noindent where $\bm{u}_{\Omega^l_i} = [u_{i-1}, u_i, u_{i+1}]^T$ and $z (\bm{u}_{\Omega^n_i}) = [u_i u_{i-1}, u^2_i, u_i u_{i+1}]^T$. Note that the nonlinear operator and the associated stencil can also be written in other ways. Typically, linear and nonlinear operators and their corresponding stencils are designed based on accuracy and stability analysis. For example, a $2^{nd}$-order centered finite difference for the nonlinear term will result in $\bm{N}^i = \frac{1}{2\Delta x}[-1, 0, 1]^T$ which is unstable with $1^{st}$-order forward Euler method for time integration, especially for hyperbolic problems such as 1-D advection equation \cite{LeVeque2007}. 

\subsection{Stability of differential operators}

The stability of a numerical scheme is commonly assessed by performing a stability analysis for some canonical linear problems. A well-known strategy for stability analysis of PDEs is to perform Von Neumann analysis \cite{Charney1950, Isaacson1994}, which involves decomposing the solution $u$ as the sum of spectral modes and assessing the growth or decay of these modes. Another similar strategy for stability analysis is to identify a semi-discrete form of the PDE by applying the method of lines with suitable spatial discretization \cite{LeVeque2007} to obtain a set of linear ODEs, which is used to perform stability analysis following common strategies in dynamical system literature \cite{Rugh1996}. We will follow the latter approach, especially as it will allow us to assess the stability property of differential operators learned from data. Consider a system of linear ODEs
\begin{equation}
    \frac{d \bm{u}}{d t} = \bm{A} \bm{u},
    \label{eq:ODEDyn}
\end{equation}
where $\bm{u} = [u_1, u_2, \cdot \cdot \cdot u_n]^T$ and $\bm{A}: \mathbb{R}^n \to \mathbb{R}^n$.  

\bigskip
\noindent \textbf{Definition 2.1}: The ODE system in \eref{ODEDyn} is referred to as stable if and only if the real part of all eigenvalues of $\bm{A}$ are nonpositive. If any eigenvalue of $\bm{A}$ is positive, \eref{ODEDyn} is called an unstable system. 
\bigskip

\begin{figure}
    \centering
    \subfigure[\label{fig:FD_fwd}]{\includegraphics[width=0.49\textwidth, trim={0.0cm 0.0cm 1.5cm 0.5cm},clip]{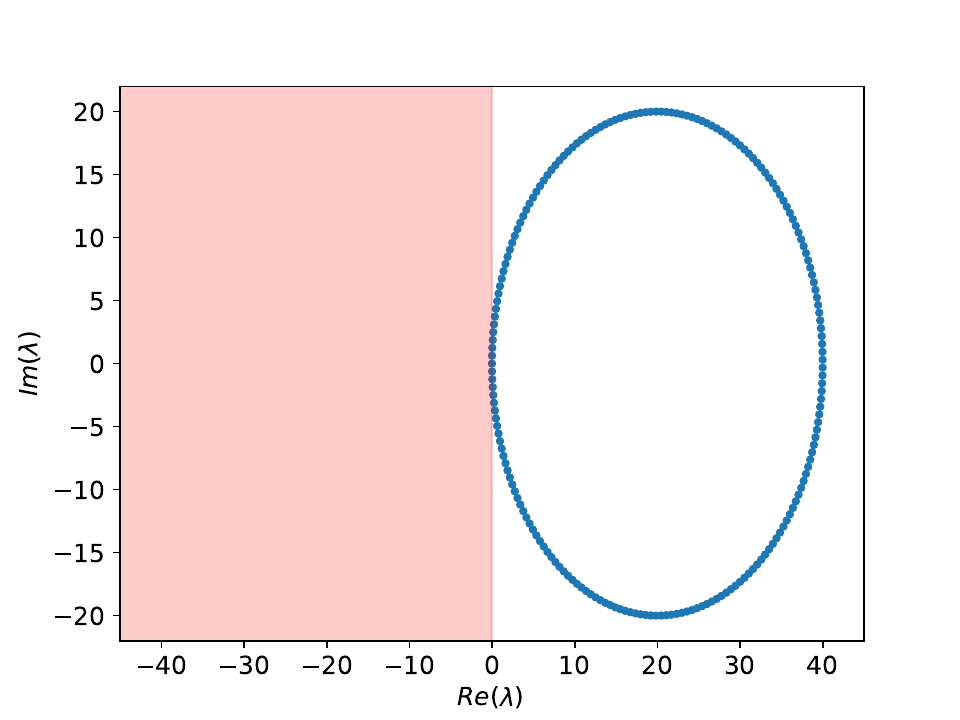}}
    \subfigure[\label{fig:FD_back}]{\includegraphics[width=0.49\textwidth, trim={0.0cm 0.0cm 1.5cm 0.5cm},clip]{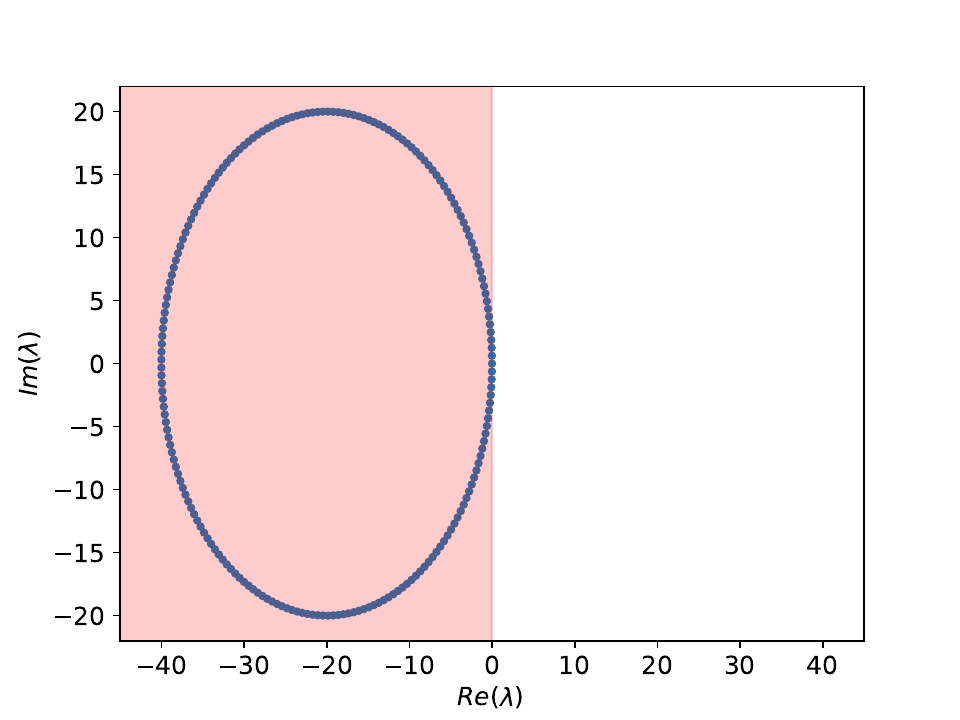}}
    \vspace{-3mm}
    \caption{Eigenvalues of the differential operator $\bm{L}$ obtained using (a) the $1^{st}$-order forward difference and (b) the $1^{st}$-order backward difference. The shaded region in red indicates the stable region.}
    \label{fig:FDStencilCompare}
\end{figure}

Generally speaking, $\bm{A}$ is referred to as asymptotically stable if it has negative eigenvalues. As the above definition holds for linear ODE systems, we refer to this property as linear stability. Such linear systems are also often referred to as globally stable as the stability properties hold irrespective of initial system states. For nonlinear ODE systems, which may arise after discretizing nonlinear PDEs, assessment of global stability may not be possible. Therefore, the local stability can be assessed by linearizing the system around an equilibrium point and performing a linear stability analysis of the resulting linear equations. We illustrate the applicability of stability analysis by considering the 1-D scalar advection equation
\begin{equation}
    \frac{\p u}{\p t} + c \frac{\p u}{\p x} = 0,
\end{equation}
where $c > 0$ is the advection velocity. To illustrate the role of a suitable selection of spatial discretization on stability, we discretize this equation to obtain a set of ODEs of the form \eref{ODEDyn}. For this example problem, we select two spatial discretization schemes, the $1^{st}$-order backward difference and the $1^{st}$-order forward difference, with uniformly spaced $n = 201$ degrees of freedom. We get a linear system with the equation for the $i^{th}$ interior degree of freedom as
\begin{equation}
   \frac{d u_i}{d t} + c (\bm{L}^i)^T \bm{u}_{\Omega^l_i} = 0,
    \label{eq:PDEeq6}
\end{equation}
where $\bm{u}_{\Omega^l_i} = [u_{i-1}, u_i, u_{i+1}]^T$, $\bm{L}^i = \frac{1}{\Delta x}[-1, 1, 0]^T$ for the $1^{st}$-order backward difference and $\bm{L}^i = \frac{1}{\Delta x}[0, -1, 1]^T$ for the $1^{st}$-order forward difference. Note that the boundary conditions will modify the definition of $\bm{L}^i$ for degrees of freedom associated with the boundaries. The eigenvalues of the linear systems $-\bm{L}$ for these two spatial discretizations are shown in \figref{FDStencilCompare}. This figure shows the backward difference yielding a linear system with all eigenvalues in the stable region, whereas the forward difference results in eigenvalues outside the stable region. This stability plot implies that backward differences yield a stable system of equations, whereas forward differences do not. This observation and analysis is common in the literature \cite{LeVeque2007}. The role of this example is solely to demonstrate how we will analyze the eigenvalues to assess the stability properties of differential operators learned from data. Even though a similar theory may not guarantee global stability for nonlinear differential equations, the local linear stability of such equations will be considered around an equilibrium point.

\section{Learning differential operators from data}
\label{sec:MathFrame}

With the increasing availability of simulation data and advancements in modern machine-learning techniques, there is a growing interest in determining differential equations from data. Several proposed methods \cite{Bar-Sinai2019, Long2019, Maddu2023} utilize artificial neural network-based architectures for representing discretized PDEs. On the other hand, interpretable sparse differential operators are determined from data using regression-based methods recently proposed in \cite{Schumann2022, Schumann2023}. The mathematical formulation of these methods is similar to the approach explained in Section \ref{sec:RegressionNormal}. All these methods do not theoretically guarantee stability, even for linear differential equations. This article proposes a novel approach incorporating suitable constraints to provide linearly stable differential operators from data. The mathematical formulation of this proposed approach is explained in Section \ref{sec:RegressionStable}.

\subsection{Regression-based approach for learning differential operator}
\label{sec:RegressionNormal}

In this article, we restrict the discussion to using least squares regression as they allow for a more straightforward and interpretable representation. Using this approach, the differential operators in \eref{PDEeq3} are obtained by solving the following regression problem:

\bigskip
\noindent \textit{Given high-fidelity data, $\bm{u} \in \mathbb{R}^n$ and $\dot{\bm{u}} = \frac{d \bm{u}}{d t} \in \mathbb{R}^n$ at time instances $t = t_j$ for $j = 1, \cdot \cdot \cdot,\; n_t$, find the optimal operators $\bm{L}$ and $\bm{N}$, subject to the objective function}
\begin{equation}
   \underset{\tilde{\bm{L}}, \; \tilde{\bm{N}}}{\text{min}} \; \Big\vert \Big\vert \dot{\bm{u}} (t) + \tilde{\bm{L}} \bm{u} (t) + \tilde{\bm{N}} \bm{z} (\bm{u} (t)) \Big\vert \Big\vert_2^2.
    \label{eq:RegProb0}
\end{equation}
The cost of the regression problem scales as $O(n_t n^2)$, which makes it expensive. Furthermore, solving learned differential equations with dense matrices for $\bm{L}$ and $\bm{N}$ will have significantly higher computational overhead than solving equations obtained using standard discretization schemes, which typically have a sparser stencil. Inspired by determining more practical and efficient differential operators from data, we pose the smaller regression problems for each degree of freedom instead. The resulting regression problem is defined as follows:

\bigskip
\noindent \textit{Given high-fidelity data, $\bm{u} \in \mathbb{R}^n$ and $\dot{\bm{u}} = \frac{d \bm{u}}{d t} \in \mathbb{R}^n$ at time instances $t = t_j$ for $j = 1, \cdot \cdot \cdot,\; n_t$ and solution stencils for linear and nonlinear operators $\bm{u}_{\Omega^l_i} \in \mathbb{R}^{s_l}$ and $\bm{u}_{\Omega^n_i} \in \mathbb{R}^{s_n}$, find the optimal operators $\bm{L}^i$ and $\bm{N}^i$, subject to the objective function}
\begin{equation}
    \underset{\tilde{\bm{L}}^i, \; \tilde{\bm{N}}^i} {\text{min}} \; \Big\vert \Big\vert \dot{u}_i (t) + (\tilde{\bm{L}}^i)^T \bm{u}_{\Omega^l_i} (t) + (\tilde{\bm{N}}^i)^T z (\bm{u}_{\Omega^n_i} (t)) \Big\vert \Big\vert_2^2 \quad \forall \; i = 1, \cdot \cdot \cdot, n.
    \label{eq:RegProb}
\end{equation}
We can determine differential operators by solving this regression problem solely using high-fidelity solution data and its time derivative. The time derivative of the solution can further be extracted from time-resolved high-fidelity data using finite differences. This approach involves solving a regression problem for each point in the domain. Therefore, assuming data at $n_t$ timesteps are available, then the cost of the regression problem at each degree of freedom is $O(n_t (s_l + s_n)^2)$. The complexity of determining the differential operator is $O(n_t n (s_l + s_n)^2)$. For a large stencil size, that is $n \approx s_l + s_n$, the cost can scale as $O(n_t n^3)$, making this problem computationally expensive. However, as typically linear and nonlinear operators are designed to be sparse, we can work with localized sparse stencils such that $n > > s_l + s_n$, resulting in a less expensive regression problem. Furthermore, the boundary conditions can be enforced by appropriately selecting local stencil and setting constraints on the local differential operators.

The training dataset should consist of high-fidelity data available at several time instances or simulation conditions to better pose the regression problem. However, even with large amounts of data, the model could overfit the available data and not generalize to simulation conditions outside the training dataset. Furthermore, the linear system we obtain in this regression problem could also be rank deficient, leading to undefined solutions \cite{Schumann2022}. This issue is overcome by augmenting the least squares regression problem in \eref{RegProb} with a regularization term. The resulting objective function is
\begin{equation}
    \underset{\tilde{\bm{L}}^i, \; \tilde{\bm{N}}^i}{\text{min}} \; \Big\vert \Big\vert  \dot{u}_i (t) + (\tilde{\bm{L}}^i)^T \bm{u}_{\Omega^l_i} (t) + (\tilde{\bm{N}}^i)^T z(\bm{u}_{\Omega^n_i} (t)) \Big\vert \Big\vert_2^2  +
    \beta_1 \vert \vert \tilde{\bm{L}}^i \vert \vert^2_2  \\ + \beta_2 \vert \vert \tilde{\bm{N}}^i \vert \vert^2_2 \quad \forall \; i = 1, \cdot \cdot \cdot, n,
    \label{eq:RegProbReg}
\end{equation}
where $\beta_1$ and $\beta_2$ are regularization constants for the unknown linear and nonlinear operators. Similarly, if the differential equations are linear, then the objective function reduces to
\begin{equation}
    \underset{\tilde{\bm{L}}^i}{\text{min}} \; \Big\vert \Big\vert  \dot{u}_i (t) + (\tilde{\bm{L}}^i)^T \bm{u}_{\Omega^l_i} (t) \Big\vert \Big\vert_2^2 + \beta_1 \vert \vert \tilde{\bm{L}}^i \vert \vert^2_2  \quad \forall \; i = 1, \cdot \cdot \cdot, n.
    \label{eq:RegProbReg_linear}
\end{equation}
By solving the regularized least squares regression problem, we can obtain the sparse operators $\bm{L}^i$ and $\bm{N}^i$. The choice of regularization here corresponds to the Tikhonov regularization \cite{Tikhonov1977} that is commonly used in the literature, resulting in a ridge regression problem. Other regression approaches, such as Lasso \cite{Tibshirani1996} and Elastic-net \cite{Zou2005}, can also be used. In our experience, these appear to give similar results to ridge regression, especially with a sparser solution stencil. The choice of stencil size, regularization approach and regularization parameter significantly influence the stability properties of the learned differential operator. There is some recent but limited work on similar sparse differential operators \cite{Schumann2022, Gkimisis2023} but without a detailed analysis of the impact of stencil size and regularization on stability. The method proposed in \cite{Schumann2023} uses statistical stability technique \cite{Maddu2022} and solves multiple similar regression problems to determine the ideal stencil size and regularization parameter that yields accurate results for a specific validation dataset window. Obtaining stable operators by solving standard unconstrained regression problems may not always be possible. Therefore, a particular combination of these parameters might yield accurate results in the validation dataset but still be theoretically unstable, yielding unphysical results for dynamics forecasting. In Section \ref{sec:Results}, we will demonstrate this behavior through several numerical experiments. Therefore, a strategy is needed to theoretically guarantee the stability of learned differential operators. 

\subsection{Constrained regression-based approach to learning stable differential operator}
\label{sec:RegressionStable}

This section describes a methodology to learn differential operators from data while ensuring stability. As stability is more precisely defined for linear operators, we formulate the method for linear operators and then extend the method for nonlinear operators. Consider a discretized PDE of the form
\begin{equation}
    \frac{d \bm{u}}{d t} + \bm{L} \bm{u} = 0
    \label{eq:PDEeq7}
\end{equation}
with unknown differential operator $\bm{L}$. We can obtain a stable differential operator by formulating a constrained regression problem:

\bigskip
\noindent \textit{Given high-fidelity data, $\bm{u} \in \mathbb{R}^n$ and $\dot{\bm{u}} = \frac{d \bm{u}}{d t} \in \mathbb{R}^n$ at time instances $t = t_j$ for $j = 1, \cdot \cdot \cdot,\; n_t$, find the optimal operator $\bm{L}$, subject to the objective function }
\begin{equation}
    \underset{\tilde{\bm{L}}}{\text{min}} \; \Big\vert \Big\vert \dot{\bm{u}} (t) + \tilde{\bm{L}} \bm{u} (t) \Big\vert \Big\vert_2^2 
    \label{eq:RegProbConstrained}
\end{equation}
\textit{subject to}
\begin{equation}
    \bm{L} \succeq 0.
    \label{eq:Constraint_1}
\end{equation}
The constraint in \eref{Constraint_1} ensures that $\bm{L}$ is a positive semi-definite matrix, which implies that $\bm{L}$ has either positive or zero eigenvalues and is stable following Definition 2.1. Ensuring stability by constraining eigenvalues has been considered in other fields such as linear dynamical systems \cite{Boots2007} and reduced order modeling \cite{Kalashnikova2014b, Sawant2023}  

This regression problem in \eref{RegProbConstrained} and \eref{Constraint_1} can be considered a semi-definite program \cite{Boyd2004}, which is more expensive to compute than \eref{RegProb0}. Furthermore, this constrained regression problem will be much more expensive for large $n$ and would not provide us $\bm{L}$ that follows the sparsity pattern obtained using standard PDE discretizations. Therefore, we pose the regression problem differently to learn a localized differential operator $\bm{L}^i$ that acts on a local solution stencil $\bm{u}_{\Omega^l_i} \in \mathbb{R}^{s_l}$ and satisfies the objective function
\begin{equation}
    \underset{\tilde{\bm{L}}^i}{\text{min}} \; \Big\vert \Big\vert \dot{u}_i (t) + (\tilde{\bm{L}}^i)^T \bm{u}_{\Omega^l_i} (t) \Big\vert \Big\vert_2^2 \quad \forall \; i = 1, \cdot \cdot \cdot, n.
    \label{eq:RegProbConstrained_local}
\end{equation}
These localized differential operators $\bm{L}^i$ are assembled to provide a differential operator $\bm{L}$ that should be positive semi-definite. To describe this regression problem mathematically, we reframe the constraint locally on $\bm{L}^i$ to ensure the assembled $\bm{L}$ is positive semi-definite. To achieve this, we utilize the commonly used Gershgorin circle theorem \cite{Gershgorin1931}.

\bigskip
\noindent \textbf{Theorem 3.1 (Gershgorin circle theorem):} Considering a complex matrix $\bm{A} \in \mathbb{C}^{n \times n}$ with the $ij^{th}$ element as $a_{ij}$. Every eigenvalue, $\lambda$, of matrix $\bm{A}$ satisfies
\begin{equation}
    \vert \lambda - a_{ii} \vert \leq \sum_{j\neq i} \vert a_{ij} \vert \quad \forall \quad i \in {1,2, \cdot \cdot \cdot, n}.
\end{equation}

\bigskip
\noindent Applying the Gershgorin circle theorem to the differential operator $\bm{L}$, every eigenvalue $\lambda$ of $\bm{L}$ must satisfy
\begin{equation}
    \vert \lambda - L^i_{i} \vert \leq \sum_{j\neq i} \vert L^i_{j} \vert \quad \forall \; i \in {1,2, \cdot \cdot \cdot, n},
\end{equation}
where $L^i_{j}$ is the $j^{th}$ element of $\bm{L}^i$ or the $ij^{th}$ element of $\bm{L}$ in common matrix notation. This inequality indicates the bounds for eigenvalues of $\bm{L}$ relating to individual terms of local operators $\bm{L}^i$. Therefore, $\bm{L} $ is guaranteed to be positive semi-definite if
\begin{equation}
     L^i_{i} \geq \sum_{i \neq j} \vert L^i_{j} \vert \quad \forall \; i \in {1,2, \cdot \cdot \cdot, n},
\end{equation}
which is a constraint on values of the local operator $\bm{L}^i$. Using this result, we can pose a constrained regression problem for each degree of freedom:

\bigskip
\noindent \textit{Given high-fidelity data, $\bm{u} \in \mathbb{R}^n$ and $\dot{\bm{u}} = \frac{d \bm{u}}{d t} \in \mathbb{R}^n$ at time instances $t = t_j$ for $j = 1, \cdot \cdot \cdot,\; n_t$ and solution stencil $\bm{u}_{\Omega^l_i} \in \mathbb{R}^{s_l}$, find the differential operator $\bm{L}^i$, subject to the objective function }
\begin{equation}
    \underset{\tilde{\bm{L}}^i}{\text{min}} \; \sum_{j=1}^{n_t} \Big\vert \Big\vert \dot{u}_i (t) + (\tilde{\bm{L}}^i)^T \bm{u}_{\Omega^l_i} (t) \Big\vert \Big\vert_2^2 \quad \forall \; i = 1, \cdot \cdot \cdot, n
    \label{eq:RegProbConstrainedFinal}
\end{equation}
\textit{subject to}
\begin{equation}
    L^i_{i} - \sum_{i \neq j} \vert L^i_{j} \vert \geq 0.
\end{equation}

\noindent This constraint ensures that we can obtain local linear operators from data while ensuring that the assembled global linear operator is stable. By solving this constrained regression problem, we obtain a differential operator that could reflect the sparsity pattern of the spatial discretization method while ensuring the stability of the assembled linear systems. To extend this approach for nonlinear PDEs, we determine the stability constraints from the corresponding linearized equations. These linearized differential equations are obtained by applying Taylor series approximation on \eref{PDEeq3} around an equilibrium point $\bm{u}^0$. This approach provides us with a linearized ODE problem
\begin{equation}
    \frac{d \bm{u}}{d t} + \bm{N}^{L} \bm{u} = 0,
    \label{eq:PDEeq8}
\end{equation}
where $\bm{N}^{L}$ is the linear term arising as a result of linearization and depends on $\bm{N}$, $\bm{L}$ and $\bm{u}^0$. The linearized equations are stable if 
\begin{equation}
\bm{N}^{L} \succ 0,
\label{eq:NonlinConstr}
\end{equation}
implying $\bm{N}^{L}$ should be a positive definite matrix. Note that $\bm{N}^{L}$ having a zero eigenvalue may not imply linear stability as this stability characteristic would depend on the truncated terms during the linearization. With this constraint, the localized stencil regression problem transforms to:

\bigskip
\noindent \textit{Given high-fidelity data, $\bm{u} \in \mathbb{R}^n$ and $\dot{\bm{u}} = \frac{d \bm{u}}{d t} \in \mathbb{R}^n$ at time instances $t = t_j$ for $j = 1, \cdot \cdot \cdot,\; n_t$ and solution stencils for the two operators $\bm{u}_{\Omega^l_i} \in \mathbb{R}^{s_l}$ and $\bm{u}_{\Omega^n_i} \in \mathbb{R}^{s_n}$, find differential operators $\bm{L}^i$ and $\bm{N}^i$, subject to the objective function}
\begin{equation}
    \underset{\tilde{\bm{L}}^i, \; \tilde{\bm{N}}^i}{\text{min}} \; \Big\vert \Big\vert \dot{u}_i (t) + (\tilde{\bm{L}}^i)^T \bm{u}_{\Omega^l_i} (t) + (\tilde{\bm{N}}^i)^T \bm{z} (\bm{u}_{\Omega^n_i} (t)) \Big\vert \Big\vert_2^2 \quad \forall \; i = 1, \cdot \cdot \cdot, n
    \label{eq:RegProbConstrainedFinalNonLin}
\end{equation}
\textit{subject to}
\begin{equation}
    N^{L,i}_{i} - \sum_{i \neq j} \vert N^{L,i}_{j} \vert > 0 \quad \forall \; i = 1, \cdot \cdot \cdot, n,
    \label{eq:ConstrainedFinalNonLin}
\end{equation}
\textit{which is the localized form of \eref{NonlinConstr} where $N^{L,i}_{j}$ is the $j^{th}$ element of $\bm{N}^{L,i}$ and $\bm{N}^{L,i}$ is assembled to form $\bm{N}^{L}$}. 

This constraint ensures that eigenvalues of $\bm{N}^L$ are positive, indicating linear stability for learned nonlinear operators. The above regression problem can theoretically ensure linear stability, implying that the learned \eref{PDEeq3} is stable around an equilibrium point. For both linear and nonlinear equations, this notion of stability classifies the resulting ODEs, or the semi-discrete form of PDEs, as linearly stable or unstable. These stable ODEs will lead to stable solutions if the eigenvalues lie within the absolute stability region of the selected time integration scheme. 

\section{Numerical results}
\label{sec:Results}

In this section, we assess the stability and performance of the differential operator obtained using the two approaches discussed in Section \ref{sec:MathFrame}: 
\begin{itemize}
   \item  Solving the standard regression problem \eref{RegProbReg} or \eref{RegProbReg_linear} which gives us learned differential operators  (LDOs); and
   \item Solving the constrained regression problem \eref{RegProbConstrainedFinal} or \eref{RegProbConstrainedFinalNonLin} which gives us stable learned differential operators (S-LDOs).
\end{itemize}
The constrained regression problem for obtaining S-LDOs is solved using a sequential least squared programming optimizer available in \cite{2020SciPy}. For all the numerical experiments, the tolerance of this optimizer is set to $10^{-6}$ to ensure that the cost of obtaining S-LDOs is comparable to the cost of obtaining LDOs. At such lower tolerances, the performance of S-LDOs is sensitive to the stencil size which is discussed in detail for each test case. We consider three test cases with different PDEs: 1-D scalar advection-diffusion equation, 1-D Burgers equation and 2-D advection equation, to demonstrate the applicability of the proposed approach for learning stable differential operators from data.

\subsection{1-D scalar advection-diffusion equation}

We consider 1-D scalar advection-diffusion PDE, which is one of the most commonly used validation cases for assessing the accuracy and stability characteristics of discretization techniques. This equation is of the form
\begin{equation}
\frac{\p u}{\p t} + c \frac{\p u}{\p x} = \nu \frac{\partial^2 u}{\p x^2},
\end{equation}
where $c$ is the advection velocity and $\nu$ is the diffusion coefficient. The semi-discrete form after the spatial discretization is 
\begin{equation}
\frac{d \bm{u}}{d t} + c \bm{L}_1 \bm{u} - \nu \bm{L}_2 \bm{u} = 0.
\end{equation}
The equation for each degree of freedom is 
\begin{equation}
\frac{d u_i}{d t} + c (\bm{L}^i_{1})^T \bm{u}_{\Omega^1_i} - \nu (\bm{L}^i_{2})^T \bm{u}_{\Omega^2_i} = 0,
\label{eq:IthAdvecDiff}
\end{equation}
where $\bm{u}_{\Omega^1_i}$ and $\bm{u}_{\Omega^2_i}$ are solution stencils for the advection and diffusion terms, whereas $\bm{L}^i_{1}$ and $\bm{L}^i_{2}$ are the local advection and diffusion differential operators that can be assembled to form $\bm{L}_{1}$ and $\bm{L}_{2}$ respectively. For scenarios where $c$ and $\nu$ are unknown, these could be included in the definition of linear differential operators $\bm{L}_1$ and $\bm{L}_2$. To generate the data, we solve this PDE using a $1^{st}$-order backward difference for the advective term and a $2^{nd}$-order centered difference for the diffusion term with $n = 201$ degrees of freedom. This spatial discretization, coupled with the temporal discretization mentioned later, gives accurate and stable results. With this spatial discretization, the differential operators for the chosen discretization are
\begin{equation}
    \bm{L}^i_{1} = \frac{1}{\Delta x}[-1, 1, 0]^T \quad \text{and}\quad \bm{L}^i_{2} = \frac{1}{2\Delta x}[-1, 0, 1]^T.
\end{equation}
In this case, the stencils for advection and diffusion terms are $\bm{u}_{\Omega^1_i} = \bm{u}_{\Omega^2_i} = [u_{i-1}, u_i, u_{i+1}]^T$. We use periodic boundary conditions for this test case. The solution of \eref{IthAdvecDiff} not only provides us with data to determine LDOs and S-LDOs but also serves as reference results for assessing the performance of learned differential operators. This PDE exhibits a mixed hyperbolic-parabolic nature based on the values of $c$ and $\nu$. We first isolate the hyperbolic and parabolic nature to study the performance of learned differential operators for pure diffusive and advective problems. After these tests, we consider learning differential operators when both hyperbolic and parabolic terms are active. 

\subsubsection{Diffusion problem: $c = 0, \nu = 0.02$}

We first consider the diffusion problem, which exhibits a parabolic nature. This problem is obtained by setting the advection velocity to zero, resulting in the following semi-discrete form of PDE for the $i^{th}$ degree of freedom:
\begin{equation}
\frac{d u_i}{d t} - \nu (\bm{L}^i_{2})^T \bm{u}_{\Omega^2_i} = 0.
\label{eq:diff_ODE}
\end{equation}
In this article, we model the semi-discrete form as
\begin{equation}
\frac{d u_i}{d t} - \nu (\bm{L}^{i,m}_{2})^T \bm{u}_{\Omega^2_i} = 0,
\label{eq:diff_ODE_model}
\end{equation}
where the modeled linear differential operator $\bm{L}^{i,m}_{2} \in \mathbb{R}^{s_l}$ is learned from the generated data and $\bm{u}_{\Omega^2_i} \in \mathbb{R}^{s_l}$ is the solution stencil of dimensionality $s_l$. The modeled differential operator of the system $\bm{L}^m_2$ is assembled from $\bm{L}^{i,m}_{2}$. We obtain $500$ snapshots of data by solving \eref{diff_ODE} using the forward Euler method with a timestep size of $0.04$. This data is used to determine LDOs and S-LDOs by solving the regression problems mentioned earlier. We normalize these local regression problems by the Euclidean norm of local solution over time to ensure similar regression problems are solved for different locations in the domain. 

\begin{figure}
    \centering
    \subfigure[\label{fig:Stab_diff_2}]{\includegraphics[width=0.32\textwidth, trim={0.0cm 0cm 1.5cm 0.5cm},clip]{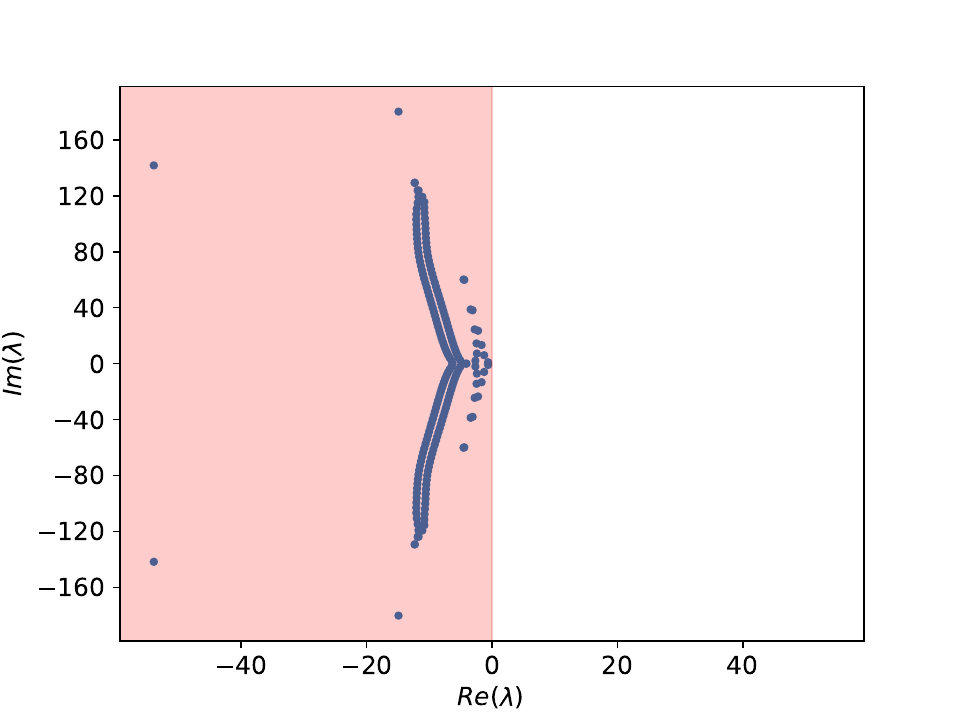}}
    \subfigure[\label{fig:Stab_diff_3}]{\includegraphics[width=0.32\textwidth, trim={0.0cm 0cm 1.5cm 0.5cm},clip]{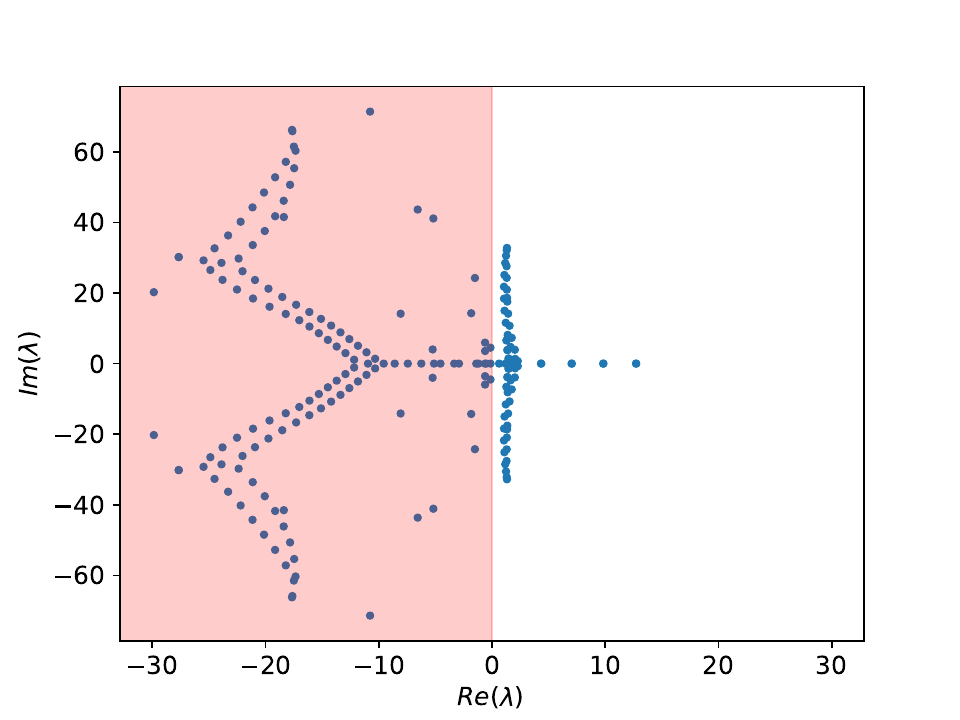}}
    \subfigure[\label{fig:Stab_diff_5}]{\includegraphics[width=0.32\textwidth, trim={0.0cm 0cm 1.5cm 0.5cm},clip]{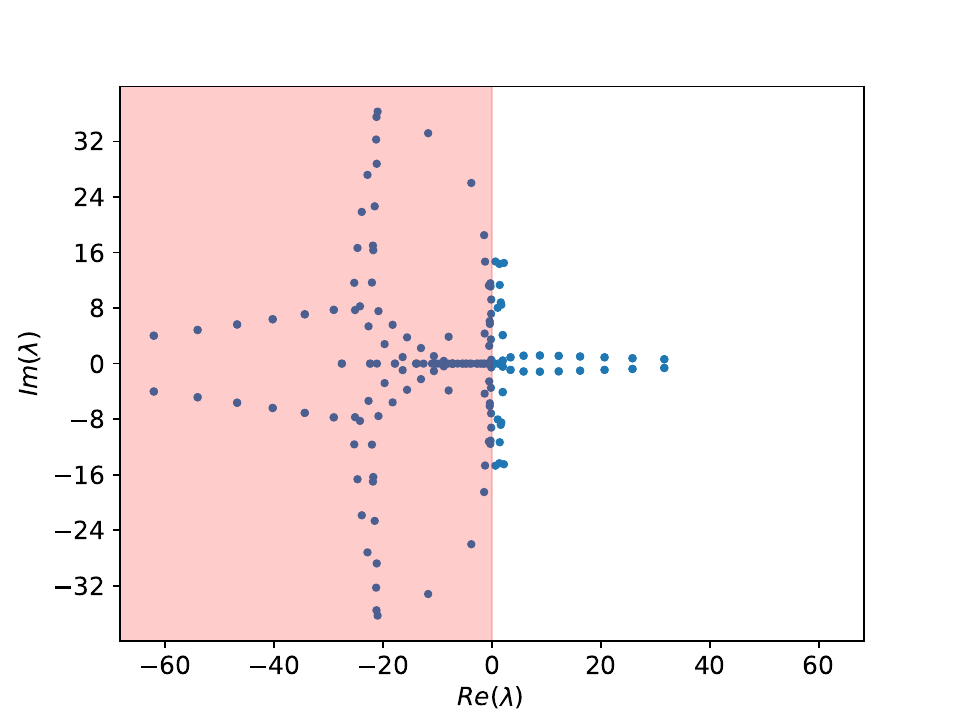}}
    \subfigure[\label{fig:Stab_diff_10}]{\includegraphics[width=0.32\textwidth, trim={0.0cm 0cm 1.5cm 0.5cm},clip]{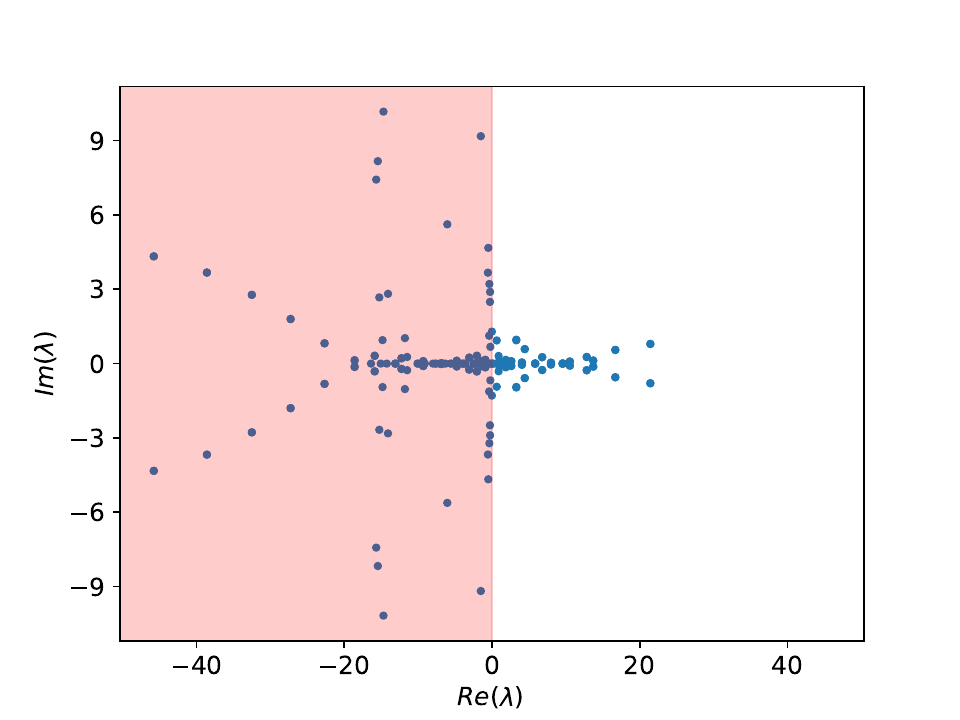}}
    \subfigure[\label{fig:Stab_diff_20}]{\includegraphics[width=0.32\textwidth, trim={0.0cm 0cm 1.5cm 0.5cm},clip]{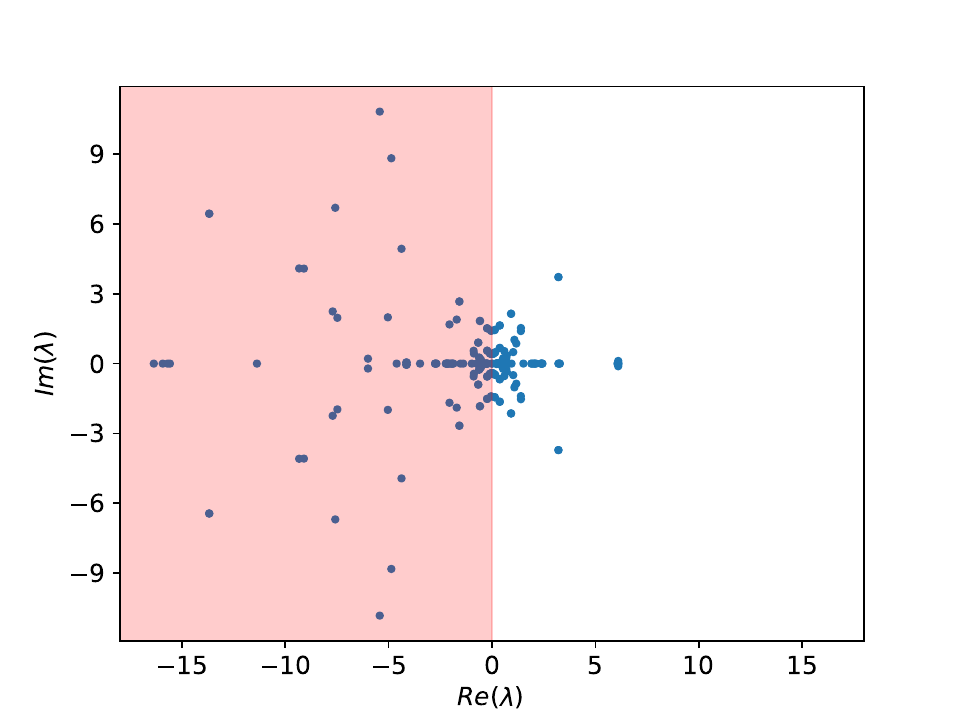}}    
    \subfigure[\label{fig:Stab_diff_40}]{\includegraphics[width=0.32\textwidth, trim={0.0cm 0cm 1.5cm 0.5cm},clip]{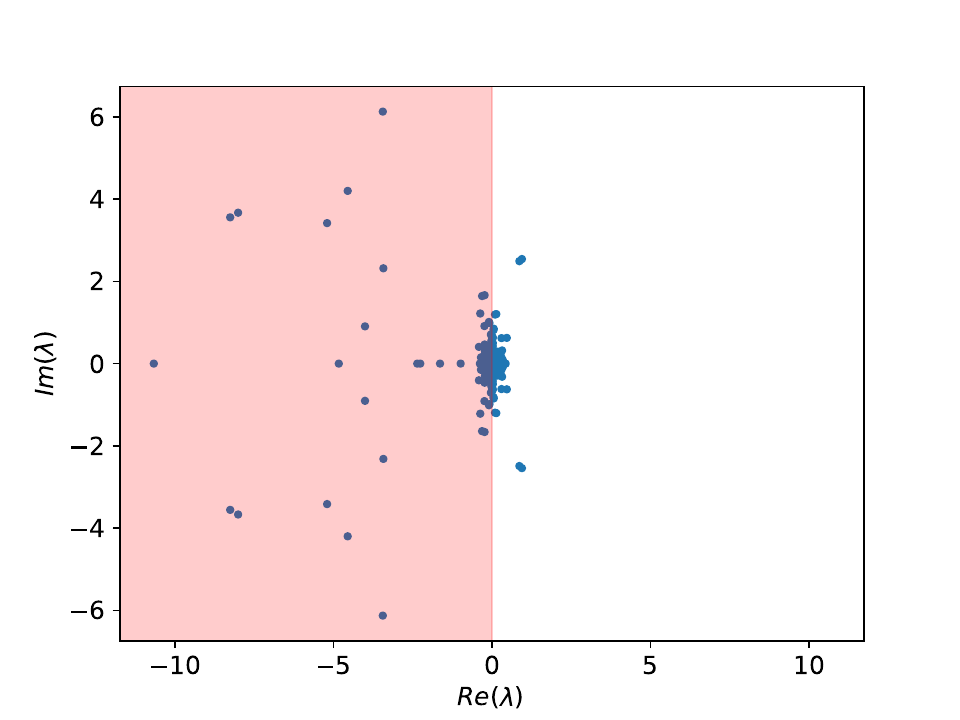}}
    \vspace{-3mm}
    \caption{Diffusion problem ($c = 0$, $\nu = 0.02$): Eigenvalues of LDOs ($\bm{L}^{m}_2$) with stencil size of (a) 3, (b) 5, (c) 7, (d) 11, (e) 21 and (f) 41 and regularization parameter $\beta_1 = 10^{-3}$. The shaded region in red indicates the stable region.}
    \label{fig:Stability_diff}
\end{figure}

\begin{figure}
    \centering
    \subfigure[\label{fig:Stab_diff_reg_1}]{\includegraphics[width=0.49\textwidth, trim={0.2cm 0cm 1.5cm 0.5cm},clip]{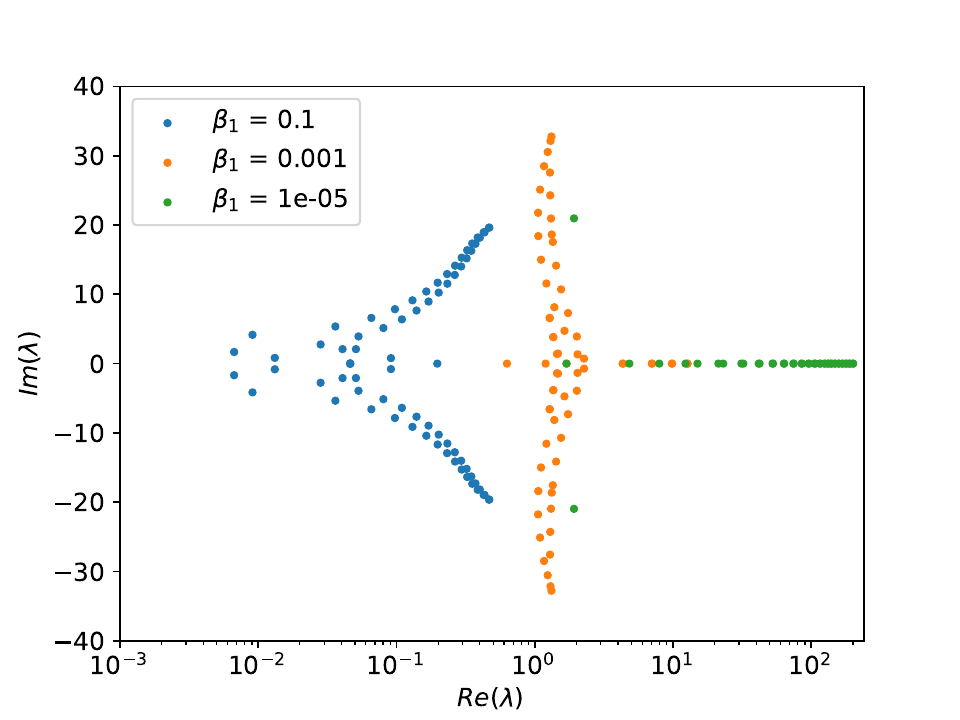}}
    \subfigure[\label{fig:Stab_diff_reg_3}]{\includegraphics[width=0.49\textwidth, trim={0.2cm 0cm 1.5cm 0.5cm},clip]{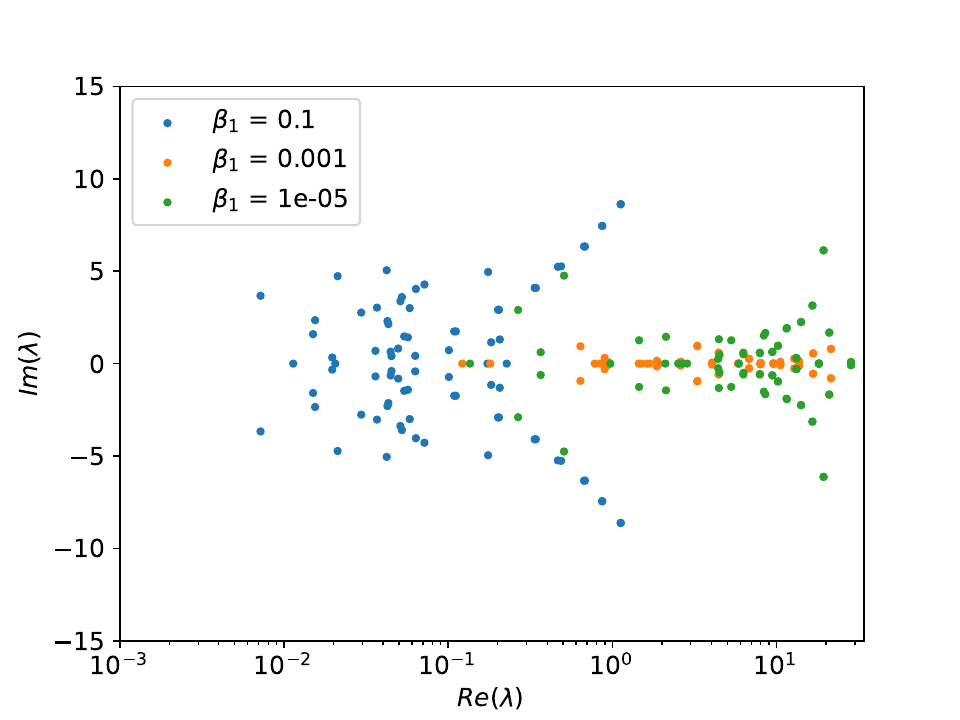}}
    \vspace{-3mm}
    \caption{Diffusion problem ($c = 0$, $\nu = 0.02$): Eigenvalues for LDO ($-\bm{L}^{m}_1$) for stencil size of (a) 5 and (b) 11 with several regularization parameters. The stable region is not shown as the $x$-axis is in the log scale.}
    \label{fig:Stability_diff_Reg}
\end{figure}
\begin{figure}
    \centering
    \includegraphics[width=0.5\textwidth]{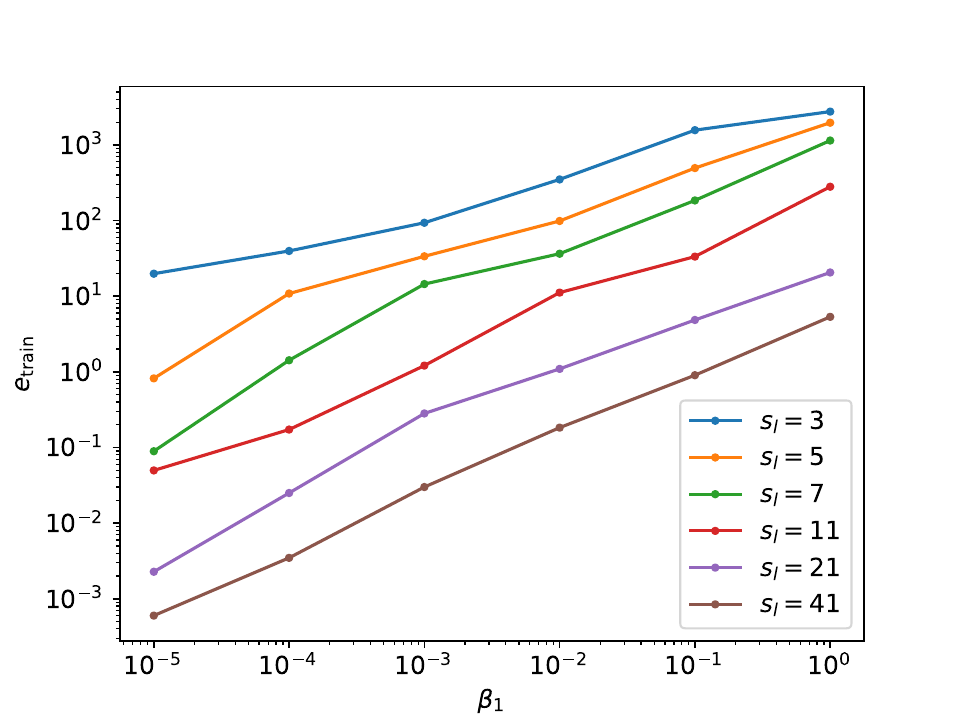}
    \vspace{-3mm}
    \caption{Diffusion problem ($c = 0$, $\nu = 0.02$): Error in the regression objective function $e_{\text{train}}$ for LDOs of different stencil sizes.}
    \label{fig:Error_diff}
\end{figure}

We first assess the stability properties of the learned differential operator by analyzing the eigenvalues of $\bm{L}^{m}_2$. These eigenvalues for LDOs of different stencil sizes are shown in \figref{Stability_diff}. All eigenvalues have a negative real part for the smallest stencil size ($s_l = 3$). On the other hand, for all other stencil sizes ($s_l > 3$), there are several eigenvalues with a positive real part. We also observe a reduction in the magnitude of eigenvalues with increased stencil size. This behavior indicates that the increase in stencil size improves the stability. Despite this improvement in stability, as several stencil sizes ($s_l > 3$) have positive eigenvalues, these learned operators are unstable. Even though the smallest stencil size ($s_l = 3$) does not follow this behavior and is found to be stable for this pure diffusion problem, we will show using the following test case that this behavior does not hold for the advection-dominated scenario.

The regression problem in \eref{RegProbReg_linear} also involves a regularization term often tuned to improve the well-posedness of the regression problem and avoid overfitting of the model. We assess the impact of the regularization parameter $\beta_1$ on the stability characteristics of the LDO. The eigenvalues of LDOs $\bm{L}^{m}_2$ for different values of $\beta_1$ are shown in \figref{Stability_diff_Reg}. The results indicate that the increase in $\beta_1$ reduces the positive real part of the eigenvalues. Therefore, an increase in $\beta_1$ improves the stability characteristics of the learned differential operator. However, these learned differential operators are unstable as positive real part of eigenvalues exist even at higher values of $\beta_1$. The error in the regression objective function 
\begin{equation}
    e_{\text{train}} = \sum_{l=0}^{n_t} \sum_{i=0}^{n} \Big( \frac{d u_i}{d t} (t_l) - \nu (\bm{L}^{i,m}_{2})^T \bm{u}_{i,\Omega^2_i} (t_l)\Big)^2
\end{equation}
is shown in \figref{Error_diff}. The error increases with the increase in $\beta_1$ for all the stencil sizes. This observation indicates that improved stability characteristics with increased $\beta_1$ come at the cost of reduced accuracy. We also observe that the increase in the stencil size not only improves stability but also reduces the error while learning LDOs. Therefore, a wider stencil would imply a more accurate and improved stability of LDOs. However, such LDOs are also accompanied by a larger matrix bandwidth, which makes them more expensive due to the higher cost of matrix-vector products. 

\begin{figure}
    \centering
    \subfigure[\label{fig:Stab_diff_1_stableC}]{\includegraphics[width=0.32\textwidth, trim={0.2cm 0cm 1.5cm 0.5cm},clip]{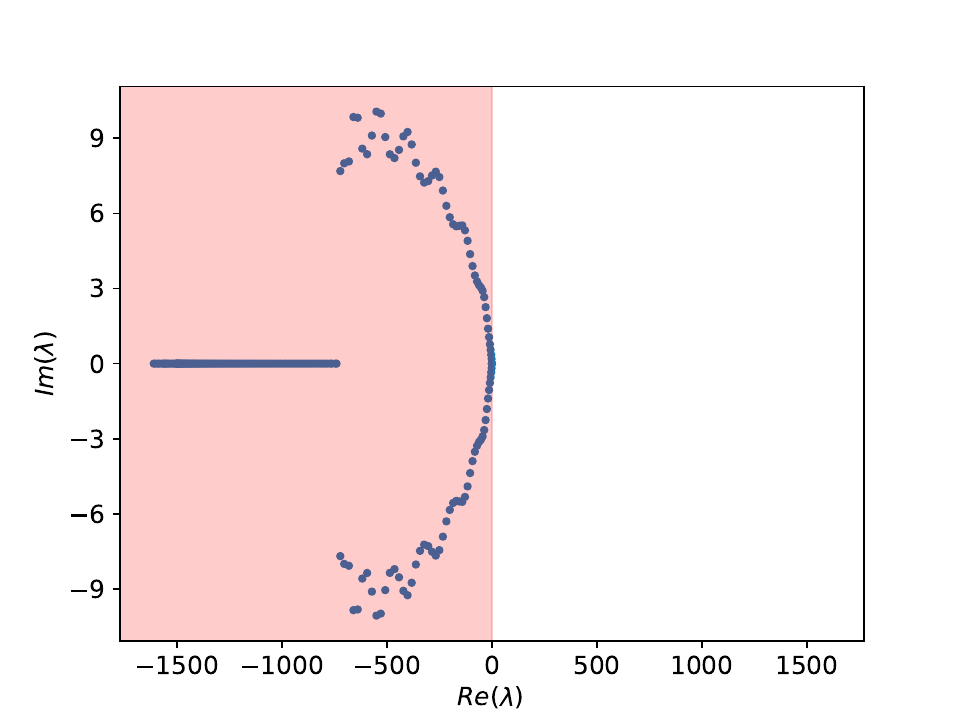}}
    \subfigure[\label{fig:Stab_diff_2_stableC}]{\includegraphics[width=0.32\textwidth, trim={0.2cm 0cm 1.5cm 0.5cm},clip]{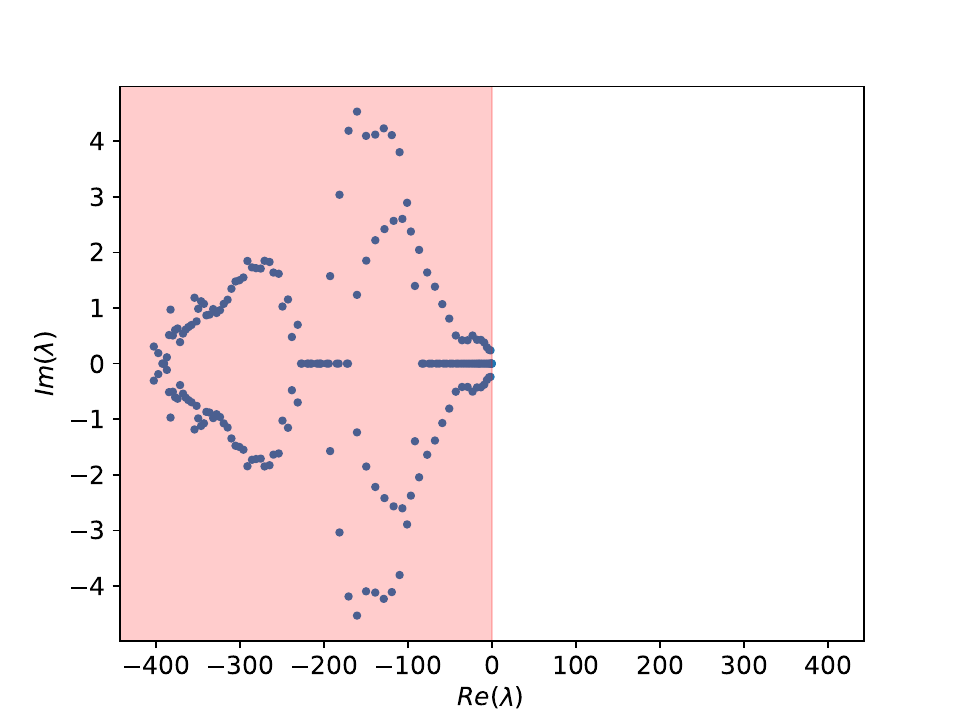}}
    \subfigure[\label{fig:Stab_diff_3_stableC}]{\includegraphics[width=0.32\textwidth, trim={0.2cm 0cm 1.5cm 0.5cm},clip]{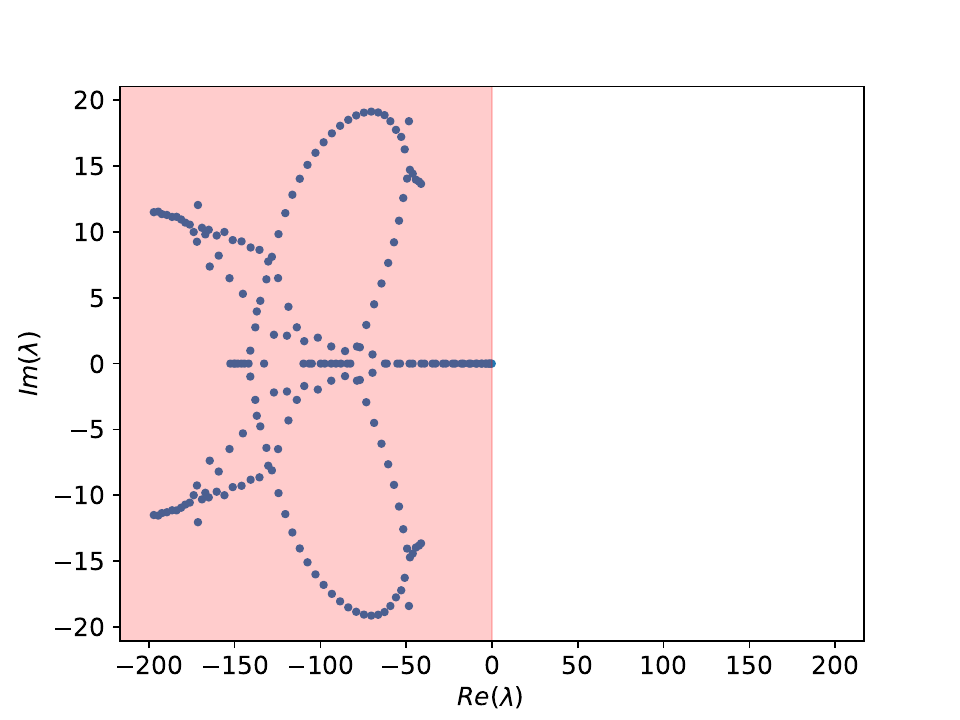}}
    \subfigure[\label{fig:Stab_diff_5_stableC}]{\includegraphics[width=0.32\textwidth, trim={0.2cm 0cm 1.5cm 0.5cm},clip]{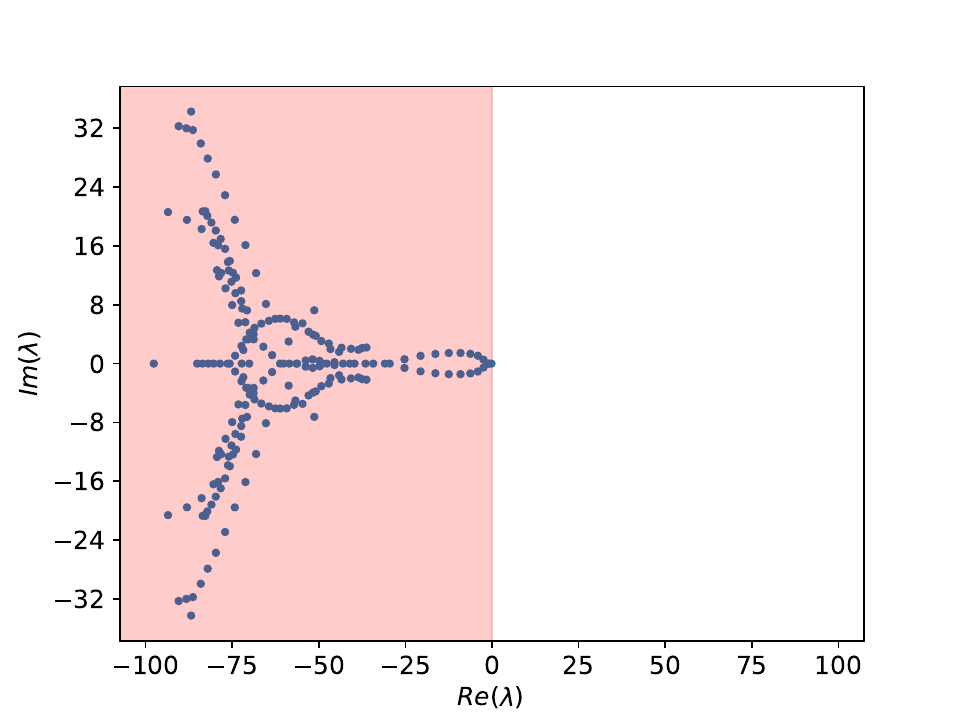}}
    \subfigure[\label{fig:Stab_diff_10_stableC}]{\includegraphics[width=0.32\textwidth, trim={0.2cm 0cm 1.5cm 0.5cm},clip]{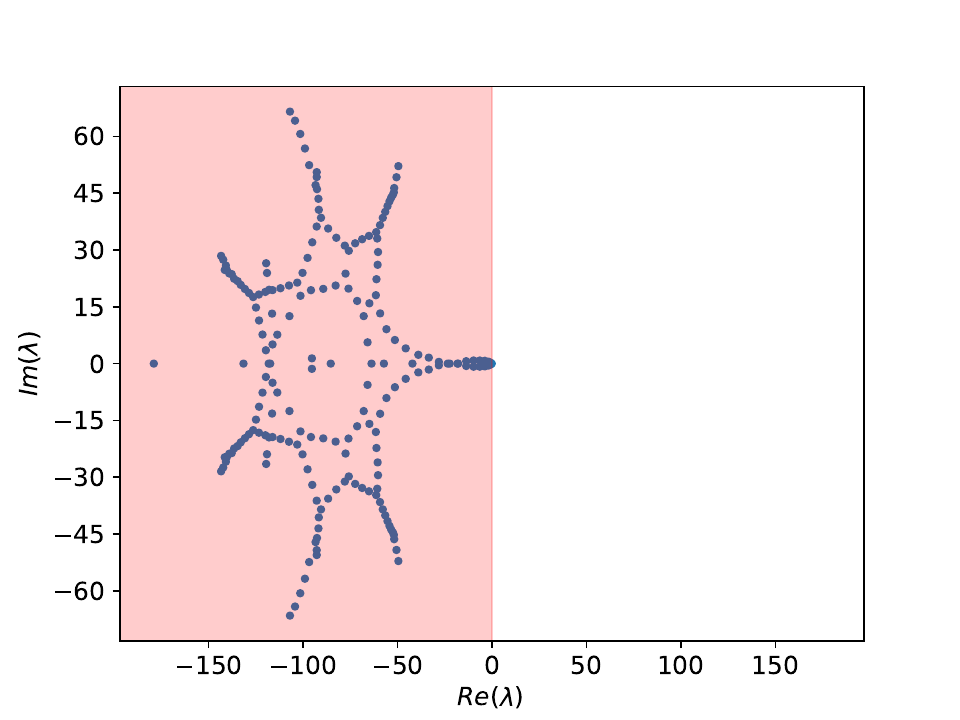}}    
    \subfigure[\label{fig:Stab_diff_20_stableC}]{\includegraphics[width=0.32\textwidth, trim={0.2cm 0cm 1.5cm 0.5cm},clip]{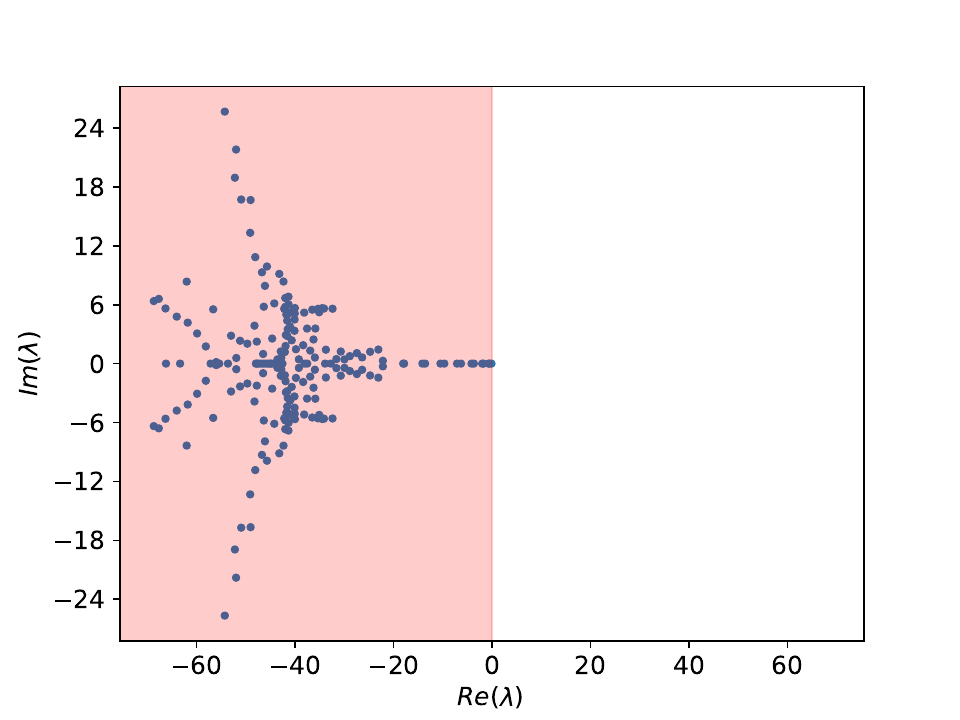}}
    \vspace{-3mm}
    \caption{Diffusion problem ($c = 0$, $\nu = 0.02$): Eigenvalues of S-LDOs ($\bm{L}^{m}_2$)  for stencil size of (a) 3, (b) 5, (c) 7, (d) 11, (e) 21 and (f) 41 with a regularization parameter of $\beta_1 = 10^{-3}$. The shaded region in red indicates the stable region.}
    \label{fig:Stability_diff_StableC}
\end{figure}

The results have indicated that LDOs for most stencil sizes are unstable even for a diffusion problem. We expect S-LDOs to resolve this issue and provide stable differential operators. We assess the stability characteristics by analyzing the eigenvalues of S-LDOs $\bm{L}^m_2$ for different stencil sizes as shown in \figref{Stability_diff_StableC}. We observe that all the eigenvalues have a negative real part, implying that S-LDOs are stable for all tested stencil sizes. With the increase in stencil size, we observed a reduction in the magnitude of the eigenvalues. A larger ratio of eigenvalue magnitudes implies stiffer system behavior and possibly decreased accuracy. The results indicate that the increase in stencil size reduces the stiffness of the system. Despite the variability in stiffness behavior, the stability of the system is still guaranteed in such situations.

\begin{figure}[t]
    \centering
    \subfigure[\label{fig:u_diff_LDO_t20}]{\includegraphics[width=0.33\textwidth, trim={0.2cm 0cm 1.5cm 0.5cm},clip]{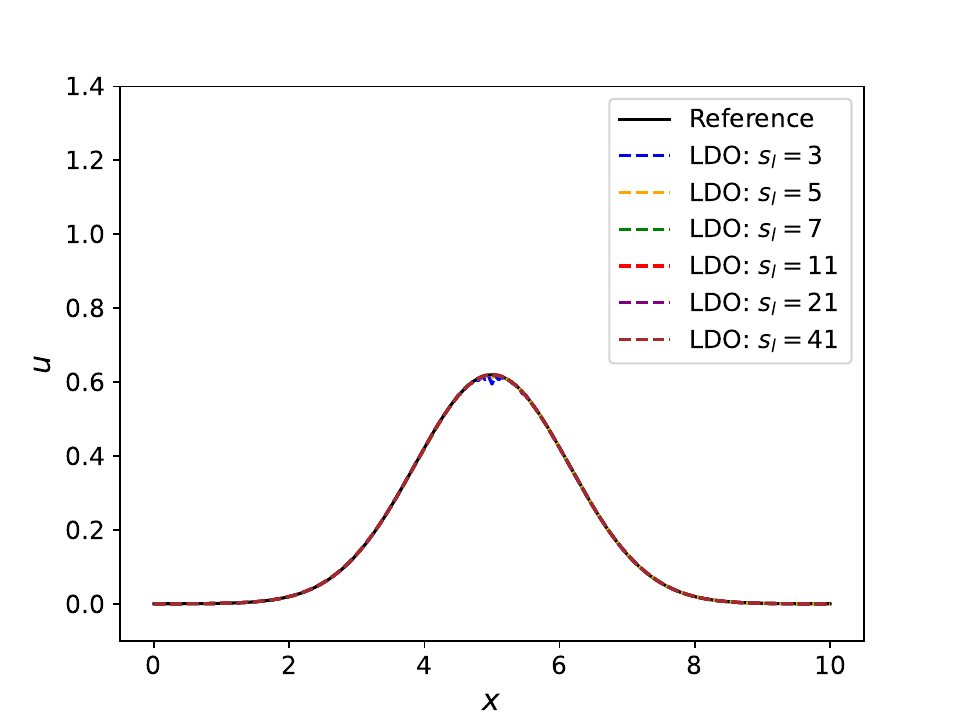}}
\subfigure[\label{fig:u_diff_LDO_t40}]{\includegraphics[width=0.33\textwidth, trim={0.2cm 0cm 1.5cm 0.5cm},clip]{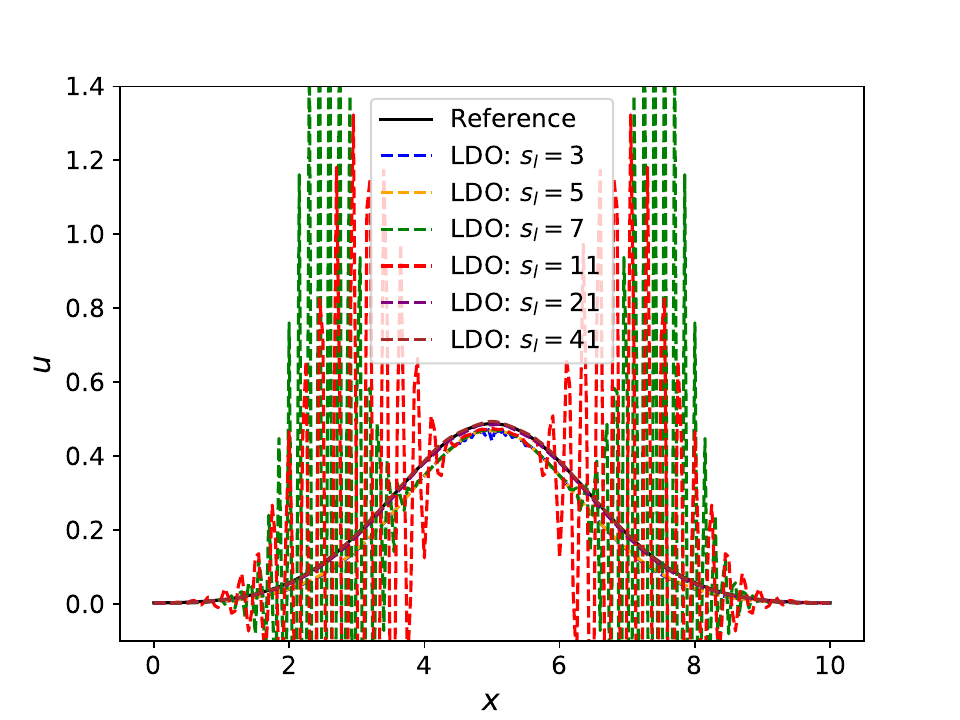}}\subfigure[\label{fig:error_diff_LDO}]{\includegraphics[width=0.33\textwidth, trim={0.2cm 0cm 1.5cm 0.5cm},clip]{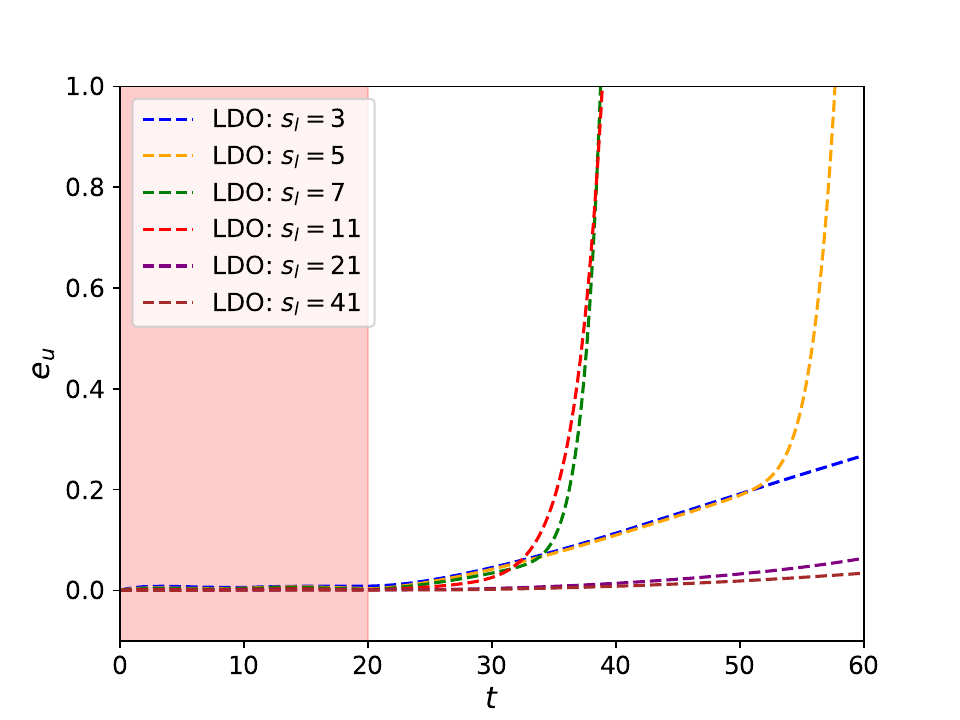}}
    \subfigure[\label{fig:u_diff_SLDO_t40}]{\includegraphics[width=0.33\textwidth, trim={0.2cm 0cm 1.5cm 0.5cm},clip]{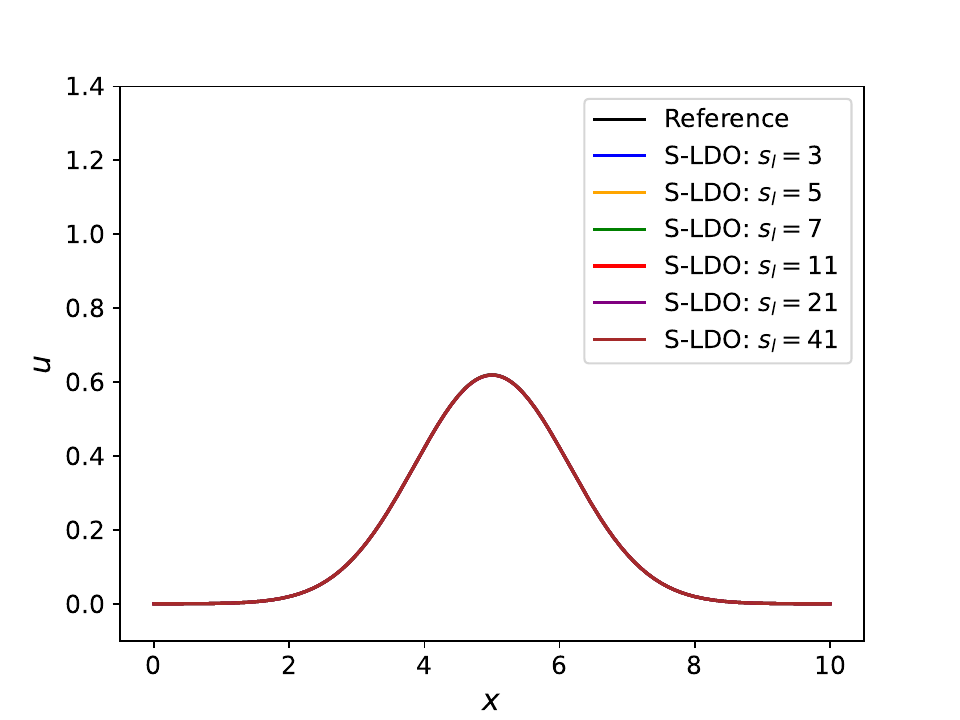}}
\subfigure[\label{fig:u_diff_SLDO_t20}]{\includegraphics[width=0.33\textwidth, trim={0.2cm 0cm 1.5cm 0.5cm},clip]{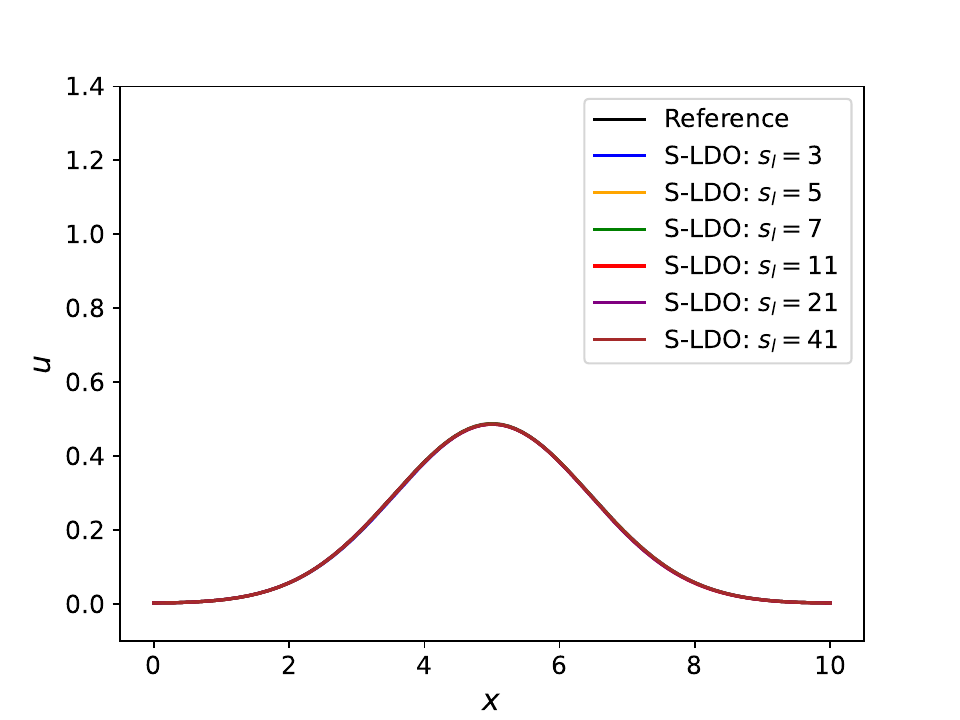}}\subfigure[\label{fig:error_diff_SLDO}]{\includegraphics[width=0.33\textwidth, trim={0.2cm 0cm 1.5cm 0.5cm},clip]{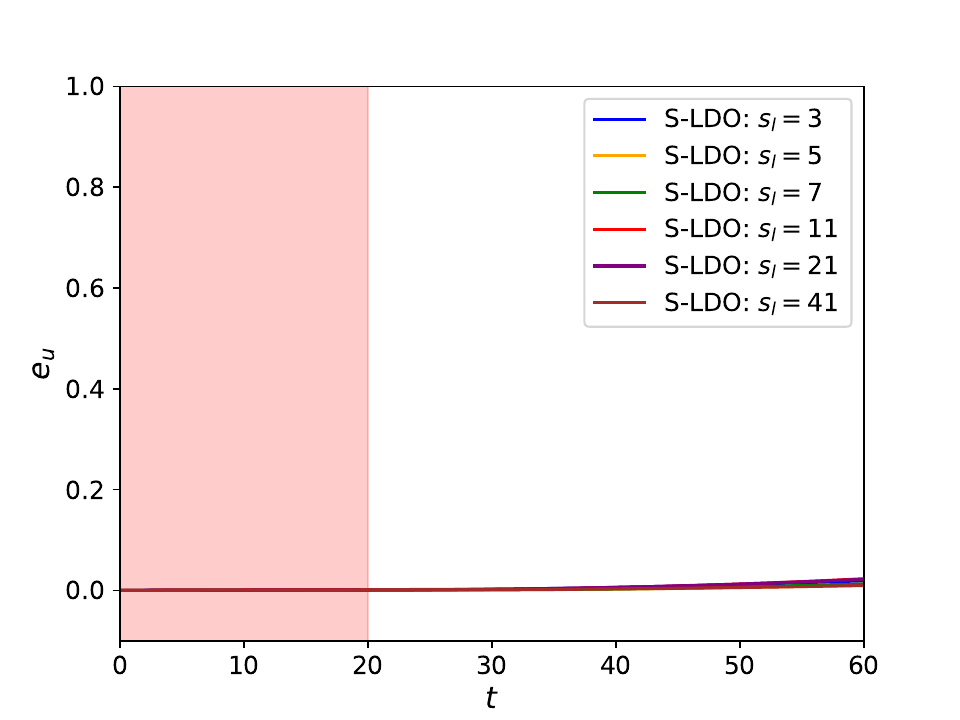}}
    
    \vspace{-3mm}
    \caption{Diffusion problem ($c = 0$, $\nu = 0.02$): Predicted solution at (a) $t = 20s$, (b) $t = 40s$ and (c) error in time using LDOs (with $\beta_1 = 10^{-3}$). Predicted solution at (d) $t = 20s$, (e) $t = 40s$ and (f) error in time using S-LDOs. The results for different stencil sizes for S-LDOs overlap. In the error plots, the unshaded region is the region of extrapolation.}
    \label{fig:u_diff_t500}
\end{figure}

Analysis of eigenvalues has allowed the identification of stability characteristics of LDOs and S-LDOs. We now assess the performance of these learned operators on replicating the dynamics of true solutions. These solution dynamics are obtained by solving \eref{diff_ODE_model} using the same time integration scheme for generating the data. We also quantify the error in solution prediction defined as
\begin{equation}
    e_u (t) = \frac{\vert \vert \bm{u} (t) - \bm{u}^m (t) \vert \vert_2^2 }{\vert \vert \bm{u} (t) \vert \vert_2^2},
\end{equation}
where $\bm{u} (t)$ is the reference solution and $\bm{u}^{m} (t)$ is the solution predicted using LDOs or S-LDOs. The predicted solutions using LDOs and S-LDOs and the corresponding errors are shown in \figref{u_diff_t500}. Both LDOs and S-LDOs provide a good estimation of the solution at peak and tails at $t = 20s$ for all stencil sizes. The smallest stencil size ($s_l = 3$) for LDOs exhibits some oscillations at the peak for $t = 20s$, whereas S-LDOs exhibit very close prediction to the reference results for all stencil sizes. At $t = 40s$, we observe oscillations in solutions predicted by LDOs, which are a consequence of the unstable nature of LDOs. On the other hand, solutions predicted by S-LDOs are very accurate even at $t = 40s$ which is outside the training dataset. These observations are also reiterated by assessing the variation of error in time. Despite the low error for LDOs in the initial time interval, the data from which is used for learning differential operators, high error and unstable behavior is observed for multiple stencil sizes outside this time interval. The error is maintained at a low value for the largest stencil sizes ($s_l \geq 21$). However, this error is expected to increase in future time instances due to the unstable nature of LDOs. On the other hand, we observe that S-LDOs exhibit much smaller errors, which do not grow substantially outside the initial time interval. Furthermore, this error is low for all stencil sizes, indicating that even sparser S-LDOs are more accurate and stable than LDOs with much larger stencil sizes. These results suggest that S-LDOs are better suited for diffusion problems than LDOs. 

\subsubsection{Advection problem: $c = 1.25, \nu = 0$}

For the diffusion problem, although LDOs do not yield as good performance as S-LDOs, errors do not rapidly blow up due to the unstable system behavior for large stencil sizes. This behavior is peculiar to this test case as diffusive problems have fewer stability issues than other advection-dominated problems. Therefore, to assess the performance of LDOs and S-LDOs for hyperbolic PDEs, we consider the advection problem with the following semi-discrete form of PDE for the $i^{th}$ degree of freedom:
\begin{equation}
\frac{d u_i}{d t} + c (\bm{L}^i_{1})^T \bm{u}_{\Omega^1_i} = 0.
\label{eq:advec_ODE}
\end{equation}
In this article, this semi-discrete form is modelled as
\begin{equation}
\frac{d u_i}{d t} + c (\bm{L}^{i,m}_{1})^T \bm{u}_{\Omega^1_i} = 0,
\label{eq:advec_ODE_model}
\end{equation}
where the modeled local differential operator $\bm{L}^{i,m}_{1} \in \mathbb{R}^{s_l}$ is learned from generated data and $\bm{u}_{\Omega^1_i} \in \mathbb{R}^{s_l}$ is the solution stencil of dimensionality $s_l$. The modeled differential operator of the system $\bm{L}^m_1$ is assembled from $\bm{L}^{i,m}_{1}$. The details of the time integration for data generation are the same as those of the diffusion problem. 

\begin{figure}
    \centering
    \subfigure[\label{fig:Stab_advec_1}]{\includegraphics[width=0.32\textwidth, trim={0.2cm 0cm 1.5cm 0.5cm},clip]{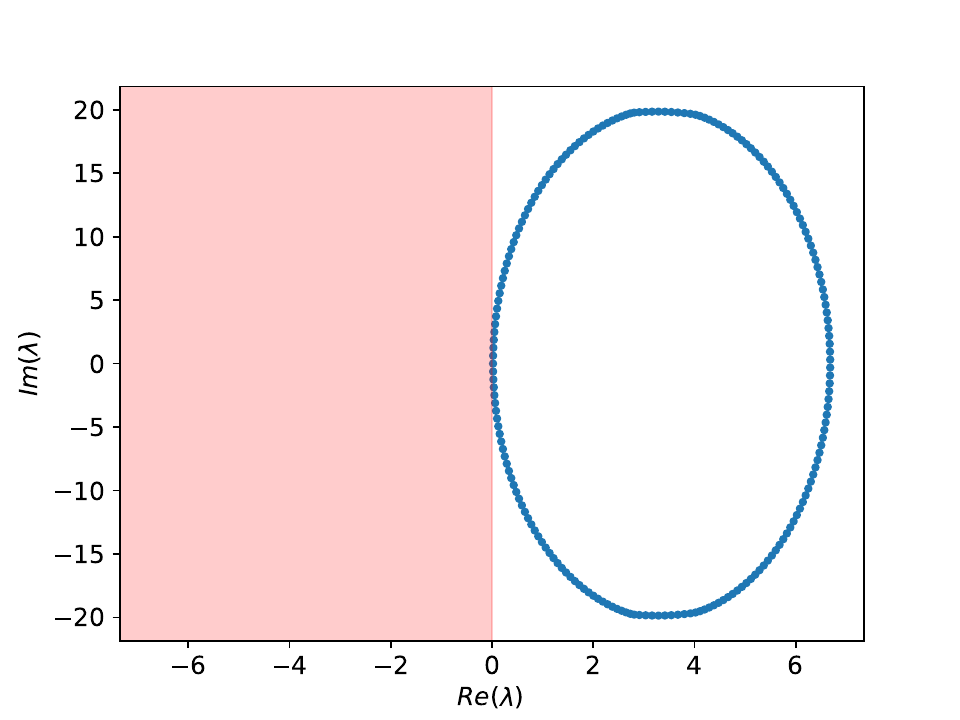}}
    \subfigure[\label{fig:Stab_advec_2}]{\includegraphics[width=0.32\textwidth, trim={0.2cm 0cm 1.5cm 0.5cm},clip]{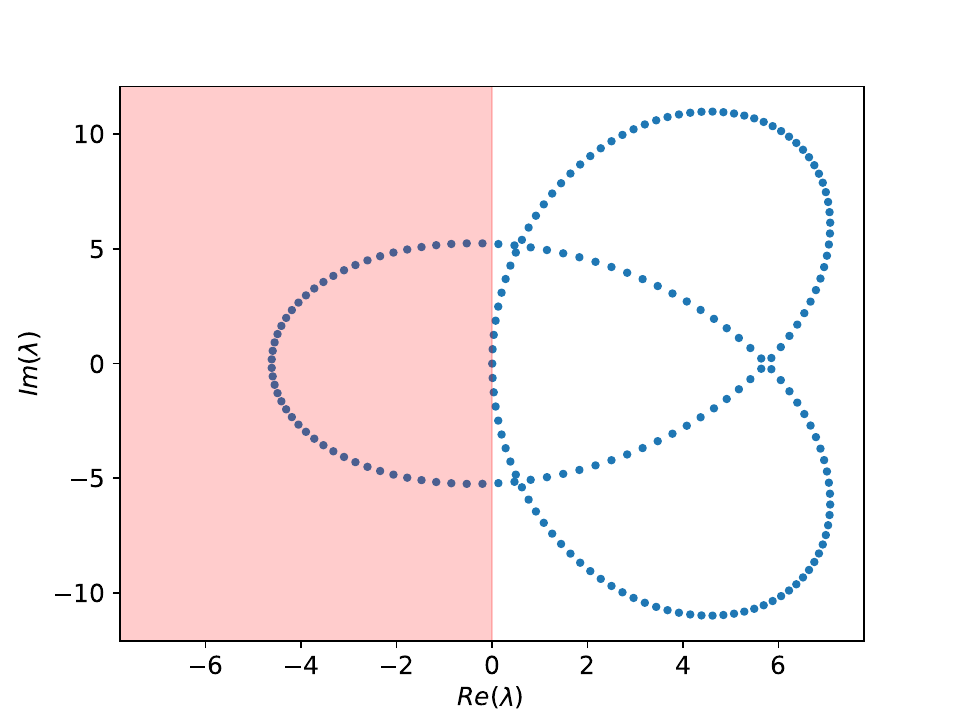}}
    \subfigure[\label{fig:Stab_advec_3}]{\includegraphics[width=0.32\textwidth, trim={0.2cm 0cm 1.5cm 0.5cm},clip]{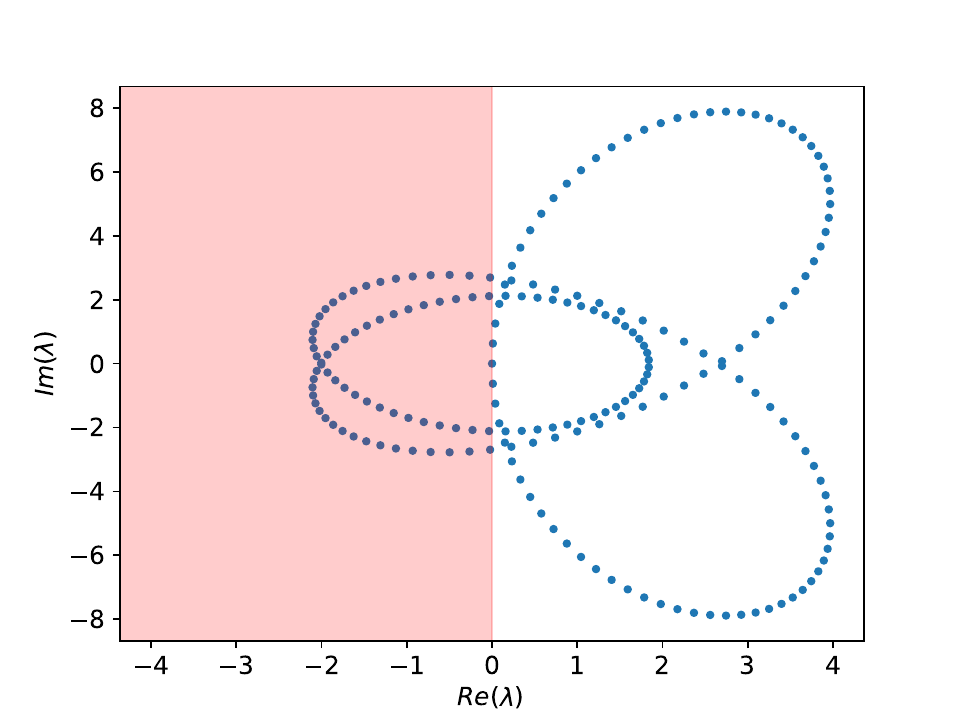}}
    \subfigure[\label{fig:Stab_advec_5}]{\includegraphics[width=0.32\textwidth, trim={0.2cm 0cm 1.5cm 0.5cm},clip]{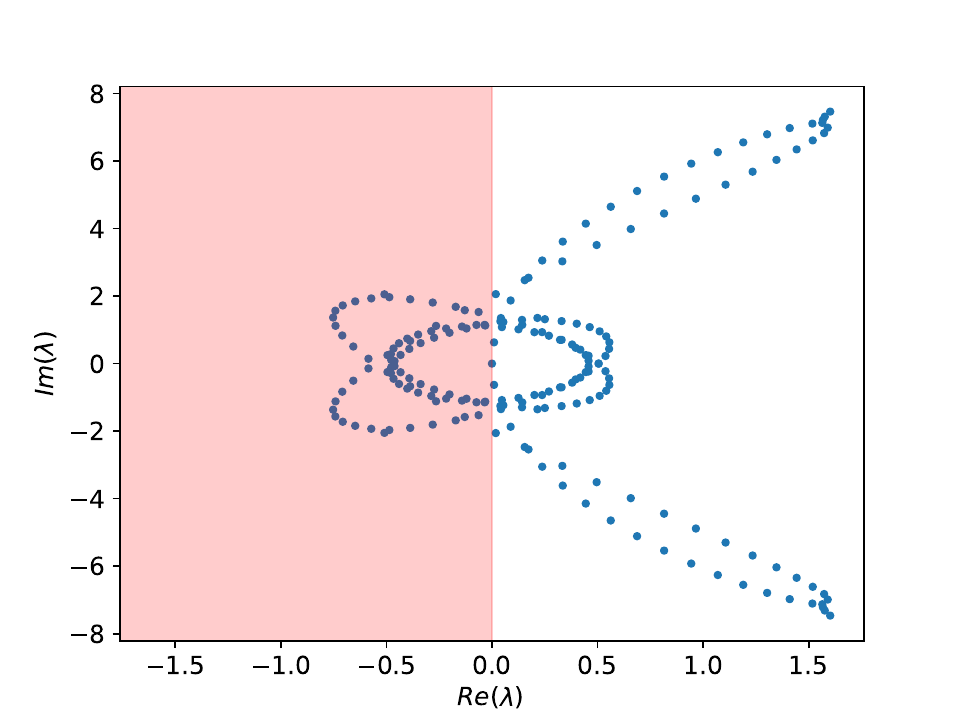}}
    \subfigure[\label{fig:Stab_advec_10}]{\includegraphics[width=0.32\textwidth, trim={0.2cm 0cm 1.5cm 0.5cm},clip]{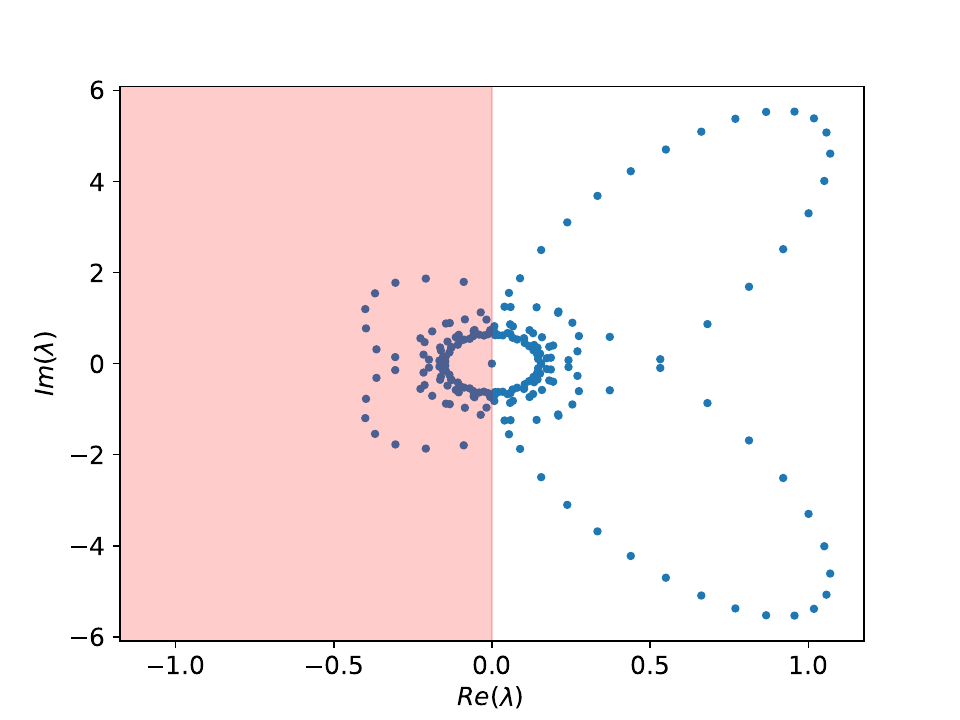}}    
    \subfigure[\label{fig:Stab_advec_20}]{\includegraphics[width=0.32\textwidth, trim={0.2cm 0cm 1.5cm 0.5cm},clip]{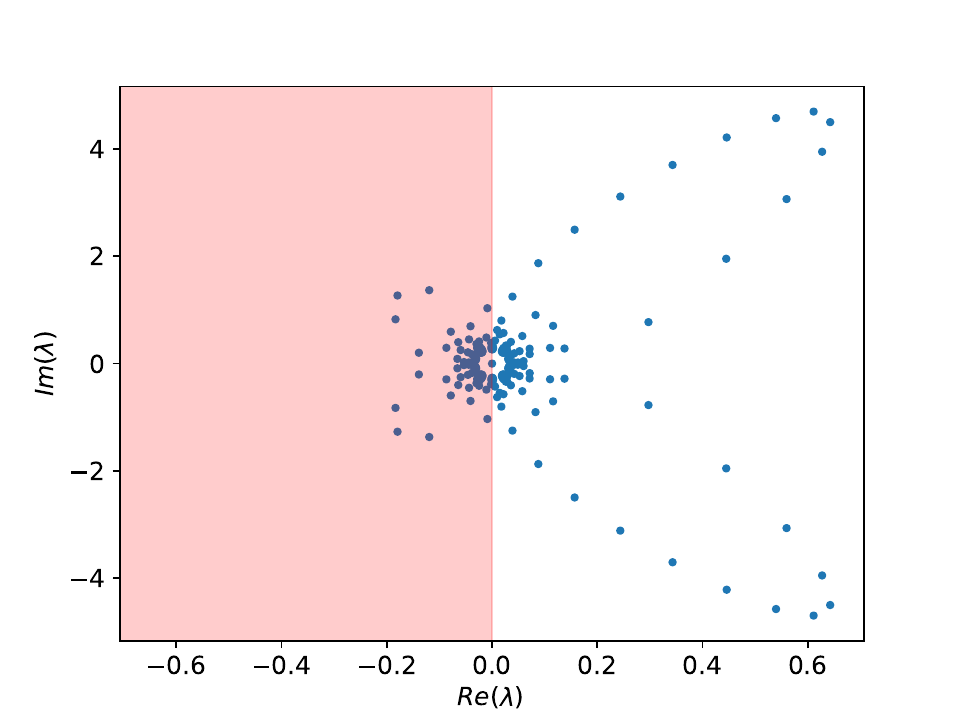}}
    \vspace{-3mm}
    \caption{Advection problem ($c = 1.25$, $\nu = 0$): Eigenvalues of LDOs ($-\bm{L}^{m}_1$) for stencil size of (a) 3, (b) 5, (c) 7, (d) 11, (e) 21 and (f) 41 with a regularization parameter of $\beta_1 = 10^{-3}$. The shaded region in red indicates the stable region.}
    \label{fig:Stability_advec}
\end{figure}

We assess the stability of the LDOs by analyzing the eigenvalues of $-\bm{L}^{m}_1$ for different stencil sizes as shown in \figref{Stability_advec}. The results indicate many eigenvalues with a positive real part for all stencil sizes. At the smallest stencil, all the eigenvalues have a positive real part. With increased stencil size, we observe more eigenvalues with negative real parts, while the rest have a smaller positive real part. This behavior matches with the diffusion problem where a larger stencil was observed to yield better stability properties in LDOs. Despite improved stability characteristics with an increased stencil size, LDOs are unstable for all these stencil sizes. 

\begin{figure}[t]
    \centering
    \subfigure[\label{fig:Stab_advec_reg_1}]{\includegraphics[width=0.32\textwidth, trim={0.2cm 0cm 1.5cm 0.5cm},clip]{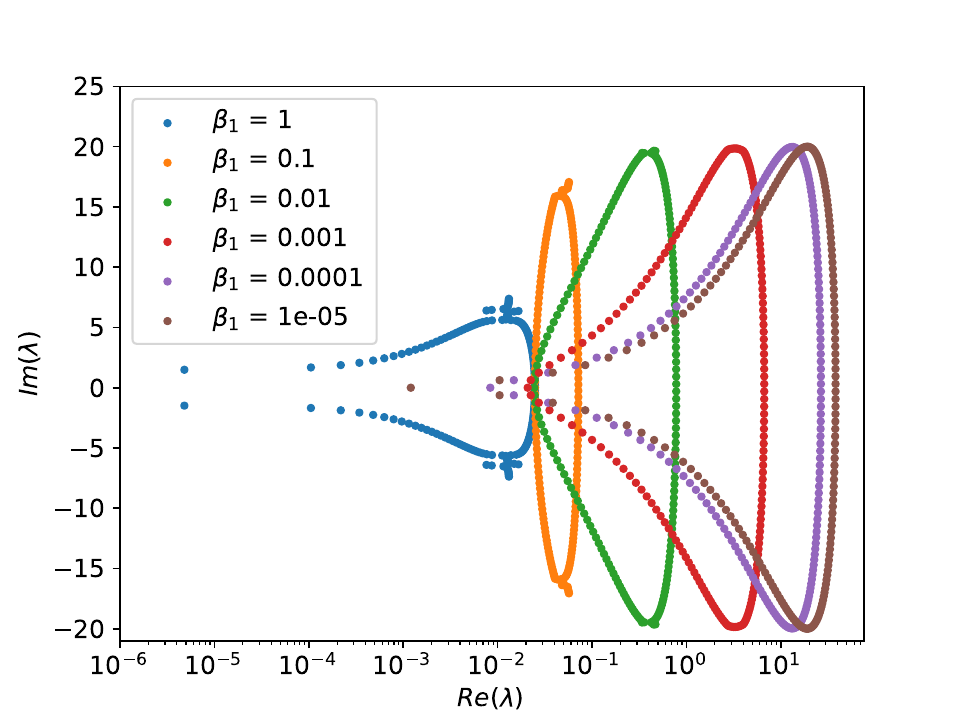}}
    \subfigure[\label{fig:Stab_advec_reg_2}]{\includegraphics[width=0.32\textwidth, trim={0.2cm 0cm 1.5cm 0.5cm},clip]{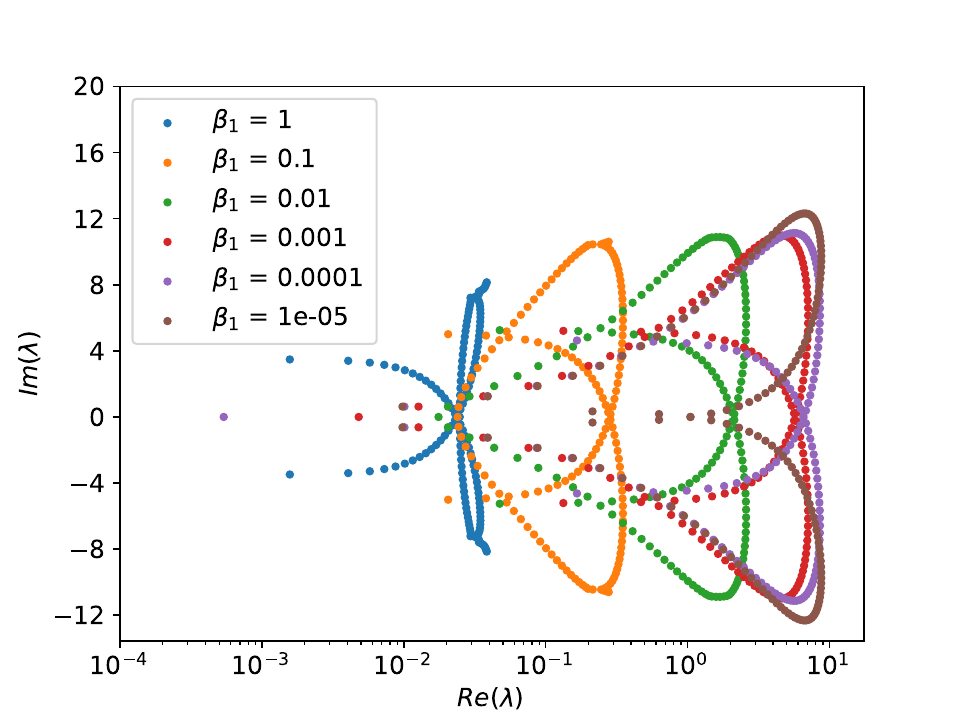}}
    \subfigure[\label{fig:Stab_advec_reg_5}]{\includegraphics[width=0.32\textwidth, trim={0.2cm 0cm 1.5cm 0.5cm},clip]{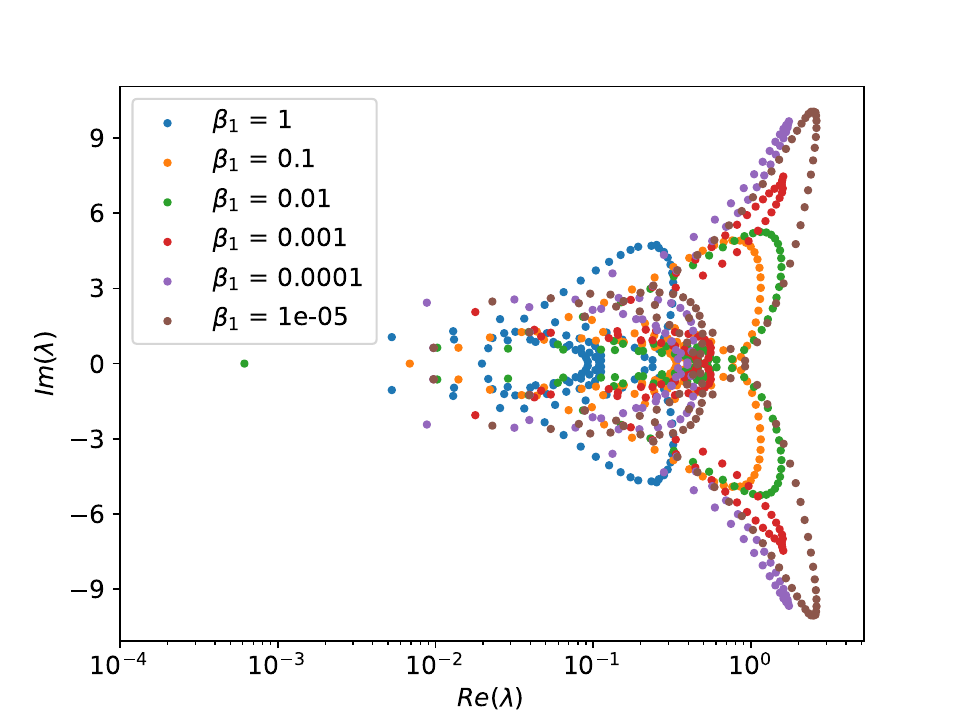}}
    \vspace{-3mm}
    \caption{Advection problem ($c = 1.25$, $\nu = 0$): Eigenvalues for LDO ($-\bm{L}^{m}_1$) for stencil size of (a) 3, (b) 5 and (c) 11 with several regularization parameters. The stable region is not shown as the $x$-axis is in the log scale.}
    \label{fig:Stability_advec_reg}
\end{figure}
\begin{figure}
    \centering
    \includegraphics[width=0.5\textwidth]{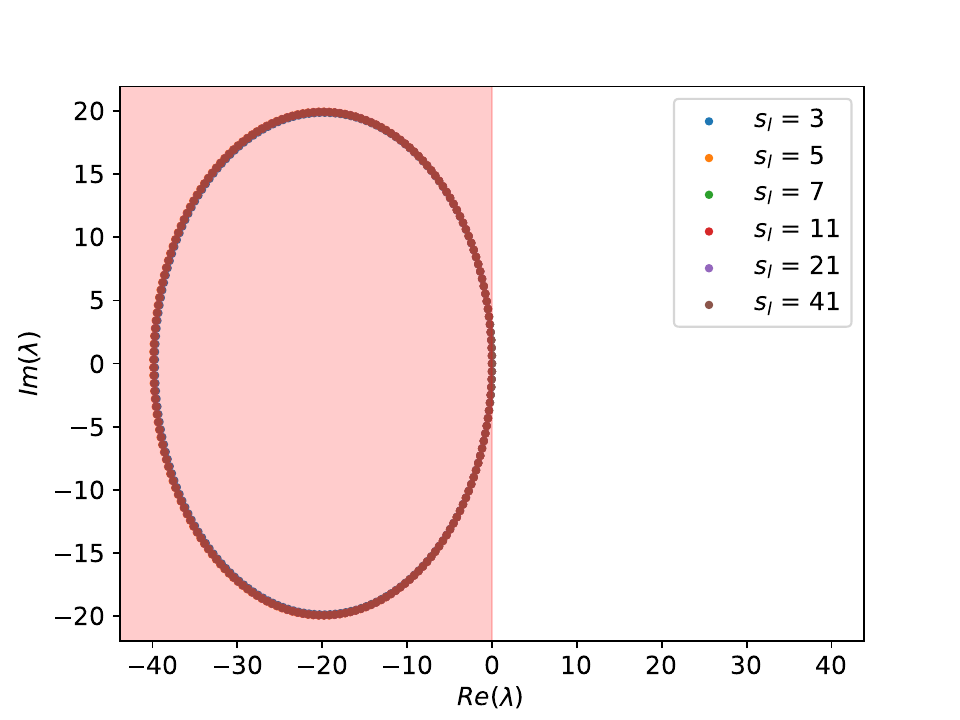}
    \vspace{-3mm}
    \caption{Advection problem ($c = 1.25$, $\nu = 0$): Eigenvalues of S-LDOs ($-\bm{L}^{m}_1$) for various stencil sizes. The shaded region in red indicates the stable region. All the tested stencil sizes lead to overlapping eigenvalues.}
    \label{fig:Stability_advec_SLDO}
\end{figure}

The effect of the regularization parameter $\beta_1$ on the eigenvalues of LDOs is shown in \figref{Stability_advec_reg}. The results indicate that an increase in $\beta_1$ brings the real part of eigenvalues closer to the origin, thereby improving the stability characteristics of the LDOs. This behavior holds for different stencil sizes as shown in the \figref{Stability_advec_reg}. As we still observe the positive real part of the eigenvalues irrespective of $\beta_1$, these LDOs are also unstable. This improvement in stability characteristics also comes at the cost of error in approximation as discussed for the pure diffusion problem. For brevity, we do not include a detailed discussion on this behavior for the pure advection case.  

The results indicated that LDOs with different stencil sizes and regularization parameters are unstable for the advection problem. We now assess the stability characteristics of S-LDOs by analyzing the eigenvalues of $-\bm{L}^{m}_1$ as shown in \figref{Stability_advec_SLDO}. The eigenvalues for different stencil sizes are very similar and overlapping. We observe that all the eigenvalues have a negative real part, implying that S-LDOs are stable. This behavior indicates that adding the constraint in the regression formulation exhibits the desired effect and enables learning a stable differential operator from data. 

\begin{figure}
    \centering
    \subfigure[\label{fig:u_advec_LDO}]{\includegraphics[width=0.33\textwidth, trim={0.2cm 0cm 1.5cm 0.5cm},clip]{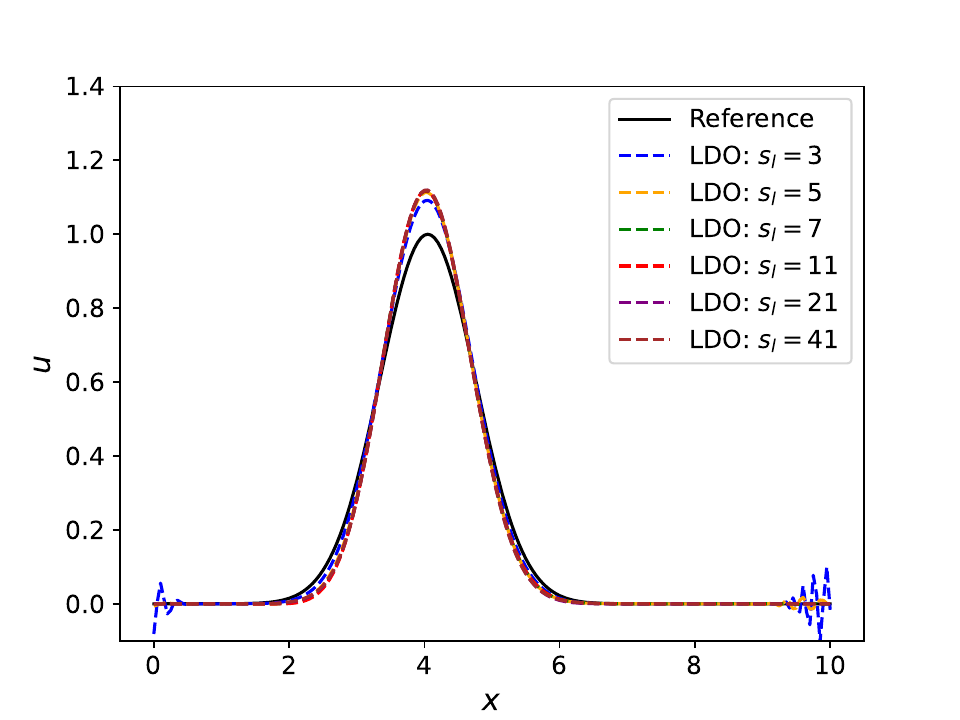}}\subfigure[\label{fig:erroru_advec_LDO}]{\includegraphics[width=0.33\textwidth, trim={0.2cm 0cm 1.5cm 0.5cm},clip]{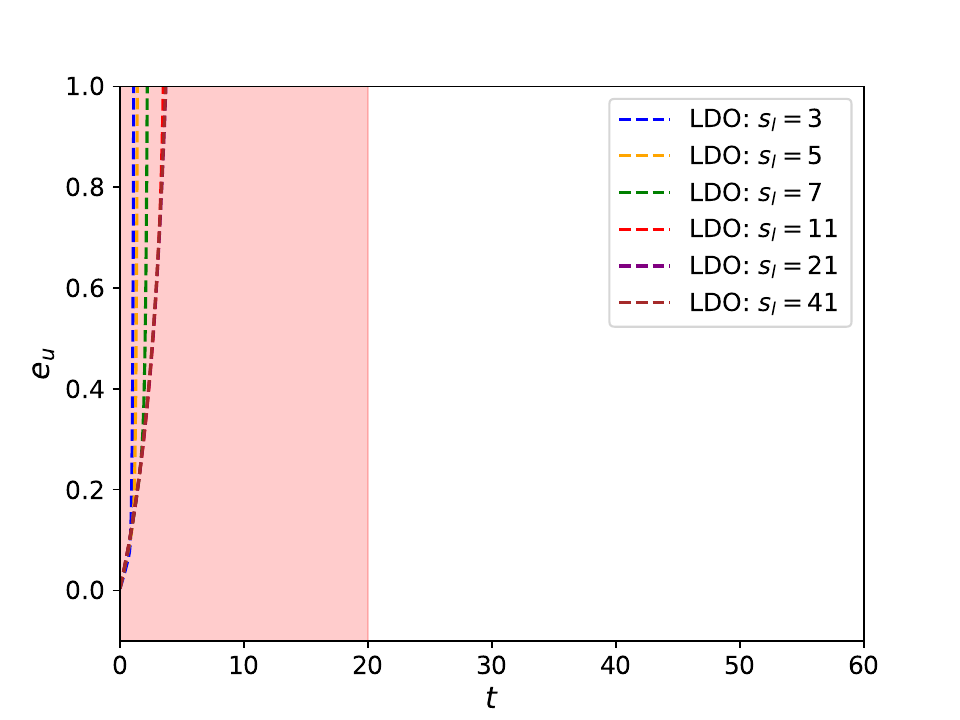}}
    
    \subfigure[\label{fig:u_advec_SLDO}]    
    {\includegraphics[width=0.33\textwidth, trim={0.2cm 0cm 1.5cm 0.5cm},clip]{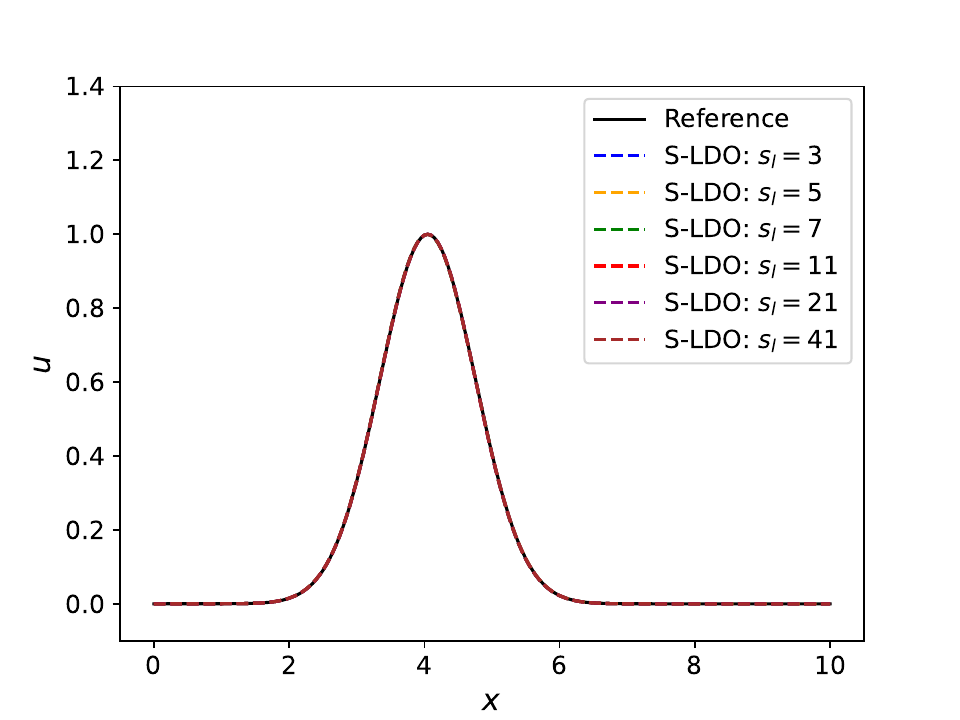}}\subfigure[\label{fig:u_advec_SLDO_t500}]{\includegraphics[width=0.33\textwidth, trim={0.2cm 0cm 1.5cm 0.5cm},clip]{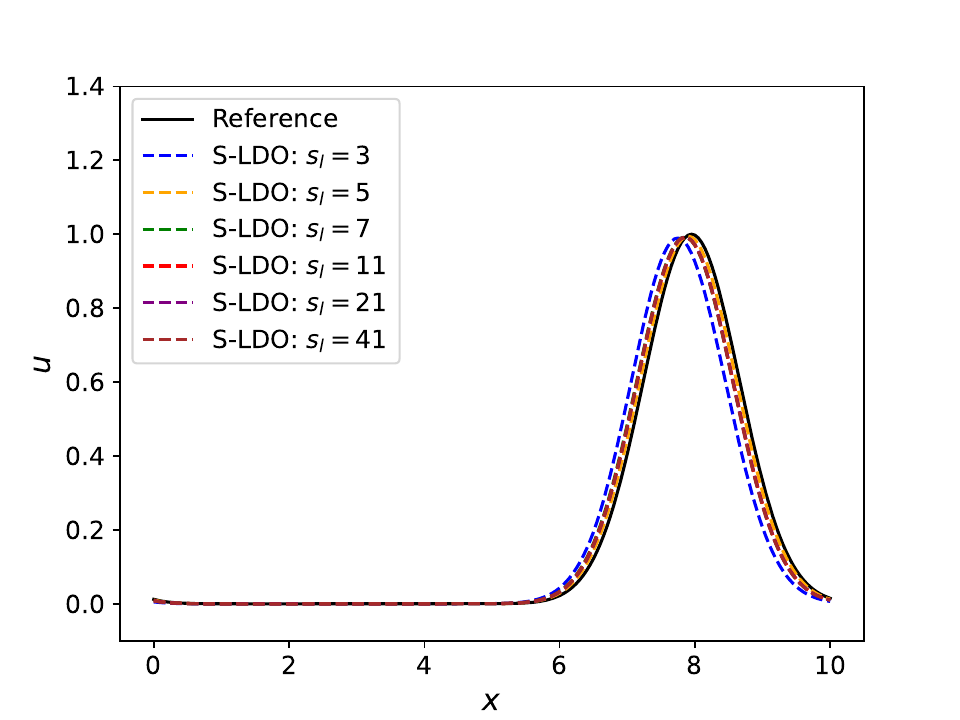}}\subfigure[\label{fig:erroru_advec_SLDO}]{\includegraphics[width=0.33\textwidth, trim={0.2cm 0cm 1.5cm 0.5cm},clip]{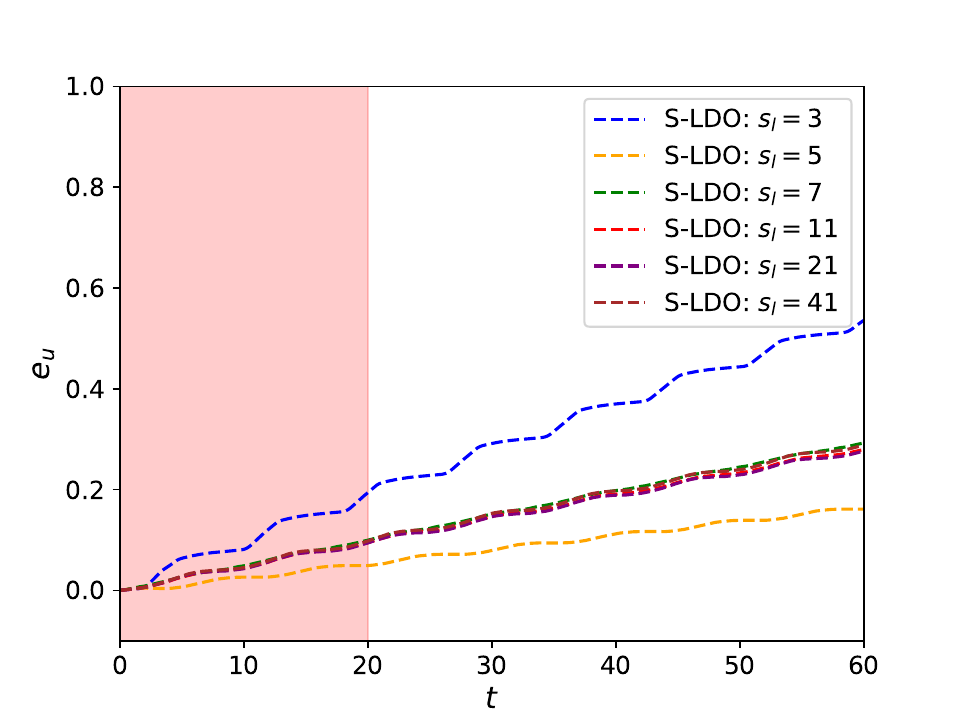}}
    \vspace{-3mm}
    \caption{Advection problem ($c = 1.25$, $\nu = 0$): Predicted velocity at (a) $t = 0.4s$ and (b) corresponding errors using LDO (with $\beta_1 = 10^{-3}$). Predicted velocity at (c) $t = 0.4s$, (d) $t = 20s$ and (e) corresponding errors using S-LDOs for different stencil sizes. In the error plots, the unshaded region is the region of extrapolation.}
    \label{fig:u_advec_comp}
\end{figure}

We now assess the performance of these learned operators on replicating the dynamics of true solutions. These solution dynamics are obtained by solving \eref{advec_ODE_model} using the same time integration scheme used for generating the data. The solutions obtained using LDOs and S-LDOs for different stencil sizes and corresponding errors are shown in \figref{u_advec_comp}. The results indicate that LDOs exhibit oscillations near the tails, especially for the coarsest stencil size, even at an early time $t = 0.4s$. The solutions predicted by LDOs at later instances are not shown because they exhibit large oscillations for all stencil sizes. The variation of prediction errors in time solidifies this point as we observe a rapid increase in errors for all stencil sizes at an early time, which is within the dataset used for determining LDOs. This behavior is also exhibited at other regularization parameters not discussed for brevity. Conversely, S-LDOs exhibit very high accuracy without any instability in the solution. Even though LDOs become unstable and yield wildly inaccurate solutions at later time instances, S-LDOs provide accurate predictions even at $t = 20s$. The error plots indicate that S-LDOs exhibit a much smaller error in time, which grows linearly in time even outside the initial time region used to determine S-LDOs. This error is maximum when the smallest stencil is used and becomes lower when larger solution stencils are selected. An exception to this behavior is $s_l = 5$, which gives the lowest errors that could be attributed to initialization and tolerance of the optimizer. 

These results indicate the importance of adding stability constraints in the regression problem for learning differential operators. Without such constraints, the learned operator may not be stable and fail to perform accurately even within the training data period. These observations for stability and accuracy can also be explained by looking at learned differential operators (scaled by $\Delta x$ and spatially averaged) using both approaches as shown in \tabref{table_advec_stencil}. S-LDOs are observed to be similar to differential operators corresponding to the $1^{st}$-order backward difference, which explains the stable behavior of these operators. With increased stencil size, this behavior is still retained, and the solution stencil does not correspond to a higher-order backward difference stencil. Therefore, despite providing better accuracy, a larger stencil does not guarantee a higher order of accuracy for the learned operator. Instead, it recovers the true differential operator with a better accuracy. On the other hand, LDOs do not directly resemble any standard finite difference stencil, especially for larger stencils.

\begin{table}[t]
    \centering
        \caption{Advection problem ($c = 1.25$, $\nu = 0$): Spatially averaged coefficients for different solution stencils for LDOs (with $\beta_1 = 10^{-3}$) and SLDOs compared to the $1^{st}$-order backward difference (BD).}
    \label{tab:table_advec_stencil}
    \begin{tabular}{|c|c|c|c|c|c|c|c|c|}
         \hline
         \multirow{2}{*}{\textbf{Approach}} &  \multirow{2}{*}{\textbf{Stencil size}} & \multicolumn{6}{c}{\textbf{Coefficients}} & \\
         \cline{3-9}
         & & $u_{i-3}$ & $u_{i-2}$ & $u_{i-1}$ & $u_{i}$ & $u_{i+1}$ & $u_{i+2}$ & $u_{i+3}$ \\
         \hline
          & $s_{l} = 3$ & - & - & $-0.415$ & $-0.167$ & $0.582$ & - & -\\
         LDO & $s_{l} = 5$ & -& $-0.084$ & $-0.160$ & $-0.116$ & $0.044$ & $0.315$ & - \\
          & $s_{l} = 7$ & $-0.041$ & $-0.080$ & $-0.079$ & $-0.046$ & $0.010$ & $0.081$ & $0.156$ \\
         \hline
          & $s_{l} = 3$ & - & - & $-0.993$ & $0.993$ & $\approx$ 0 & - & - \\
         SLDO & $s_{l} = 5$ & - & $\approx 0$ & $-0.998$ & $0.998$ & $\approx 0$ & - & - \\
          & $s_{l} = 7$ & $\approx 0$ & $\approx 0$ & $-0.996$ & $0.996$ & $\approx 0$ & $\approx 0$ & $\approx 0$\\
         \hline
         BD &  & - & - & -1 & 1 & - & - & -\\ 
         \hline         
    \end{tabular}
\end{table}

\subsubsection{Case with both advection and diffusion: $c = 0.2$, $\nu = 0.02$}

In this third scenario, we explore a mixed hyperbolic-parabolic regime with active advection and diffusion terms. Common scalar transport equations often have this form. The resulting differential equations for the $i^{th}$ degree of freedom is given in \eref{IthAdvecDiff}. In this article, we model the semi-discrete form as
\begin{equation}
\frac{d u_i}{d t} + c (\bm{L}^{i,m}_{1})^T \bm{u}_{\Omega^1_i} - \nu (\bm{L}^{i,m}_{2})^T \bm{u}_{\Omega^2_i} = 0,
\label{eq:AdvecDiff_ODE_model}
\end{equation}
where the modeled linear differential operators $\bm{L}^{i,m}_{1} \in \mathbb{R}^{s_{l_1}}$ and $\bm{L}^{i,m}_{2} \in \mathbb{R}^{s_{l_2}}$ are learned from generated data. The local velocity stencils are denoted by $\bm{u}_{\Omega^1_i} \in \mathbb{R}^{s_{l_1}}$ and $\bm{u}_{\Omega^2_i} \in \mathbb{R}^{s_{l_2}}$, where $s_{l_1}$ and $s_{l_2}$ denote the dimensionality of these stencils. We obtain $1,000$ snapshots of data by solving \eref{diff_ODE} using the forward Euler method with a timestep size of $0.04$. Assessment of stability properties by identifying eigenvalues yielded similar results and discussion to those observed for pure advection and diffusion cases. Even though we are not discussing them for brevity, the results indicated that LDOs are unstable. In contrast, S-LDOs result in a stable system as we ensure that eigenvalues for $-c \bm{L}^{m}_1 + \nu \bm{L}^{m}_2$ are negative by solving a constrained regression problem. Instead, we assess the performance of these learned operators on replicating the dynamics of the true solution, which is referred to as the reference solution in the figures. The predicted dynamics are obtained by solving \eref{AdvecDiff_ODE_model} using the same time integration scheme used for generating the data.

\begin{figure}
    \centering
    \subfigure[\label{fig:u_both_LDO_500}]{\includegraphics[width=0.33\textwidth, trim={0cm 0cm 1.5cm 0.5cm},clip]{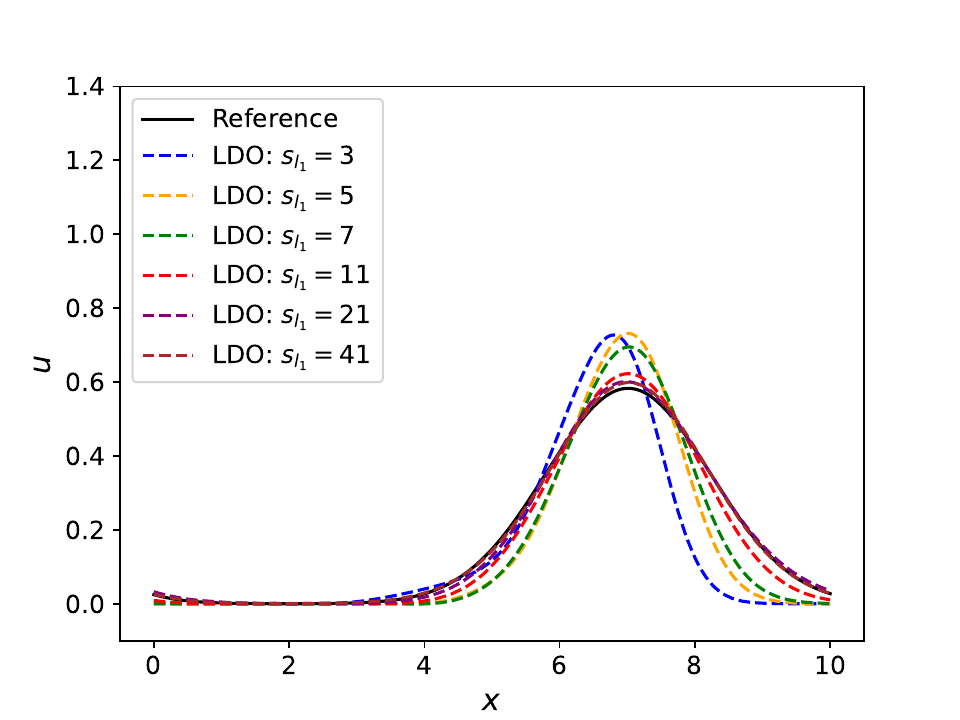}}\subfigure[\label{fig:u_both_LDO_1000}]{\includegraphics[width=0.33\textwidth, trim={0cm 0cm 1.5cm 0.5cm},clip]{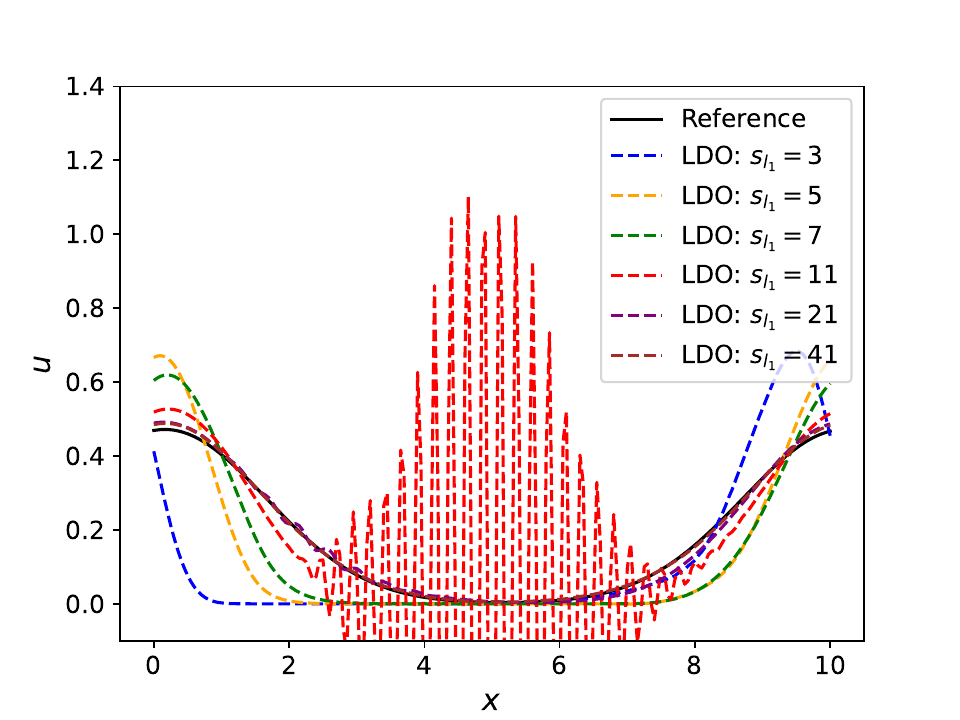}}\subfigure[\label{fig:error_both_LDO}]{\includegraphics[width=0.33\textwidth, trim={0cm 0cm 1.5cm 0.5cm},clip]{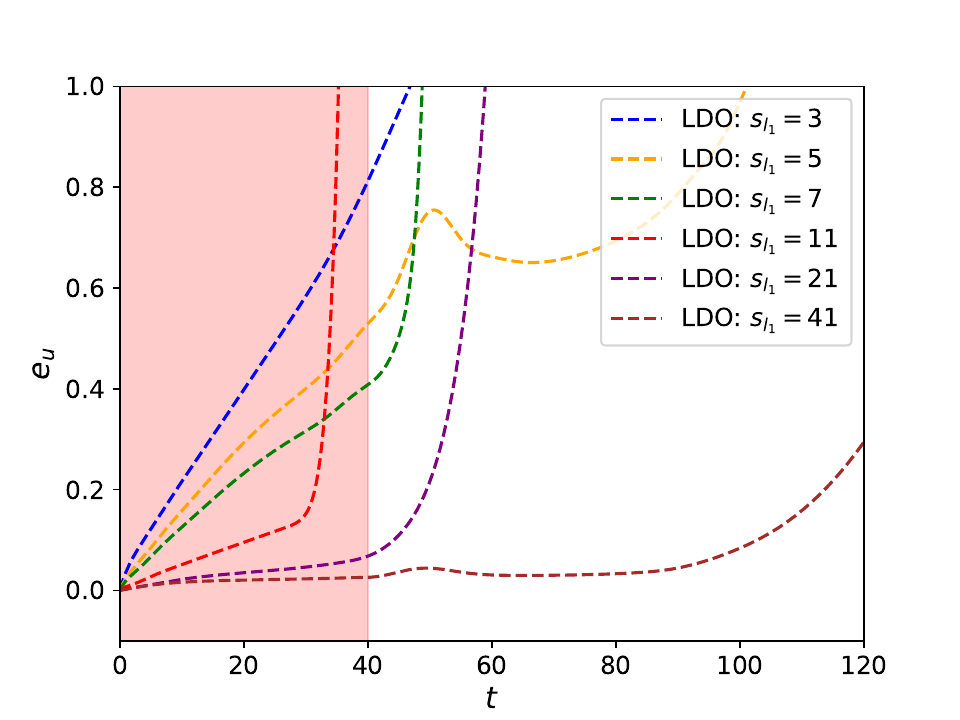}}
    \subfigure[\label{fig:u_both_SLDO_500}]{\includegraphics[width=0.33\textwidth, trim={0cm 0cm 1.5cm 0.5cm},clip]{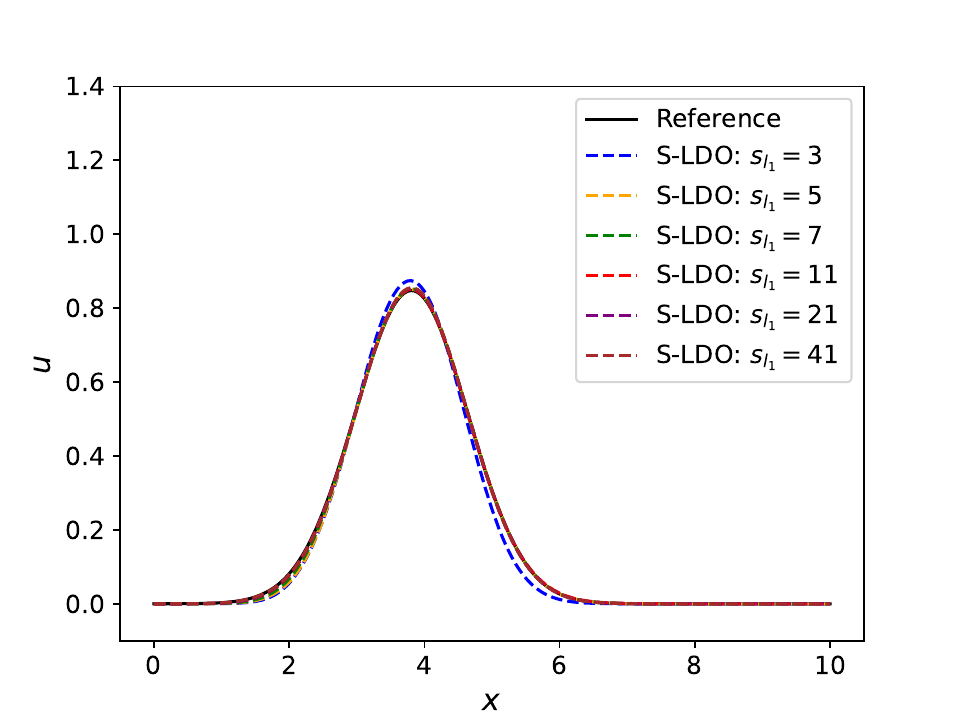}}\subfigure[\label{fig:u_both_SLDO_1000}]{\includegraphics[width=0.33\textwidth, trim={0cm 0cm 1.5cm 0.5cm},clip]{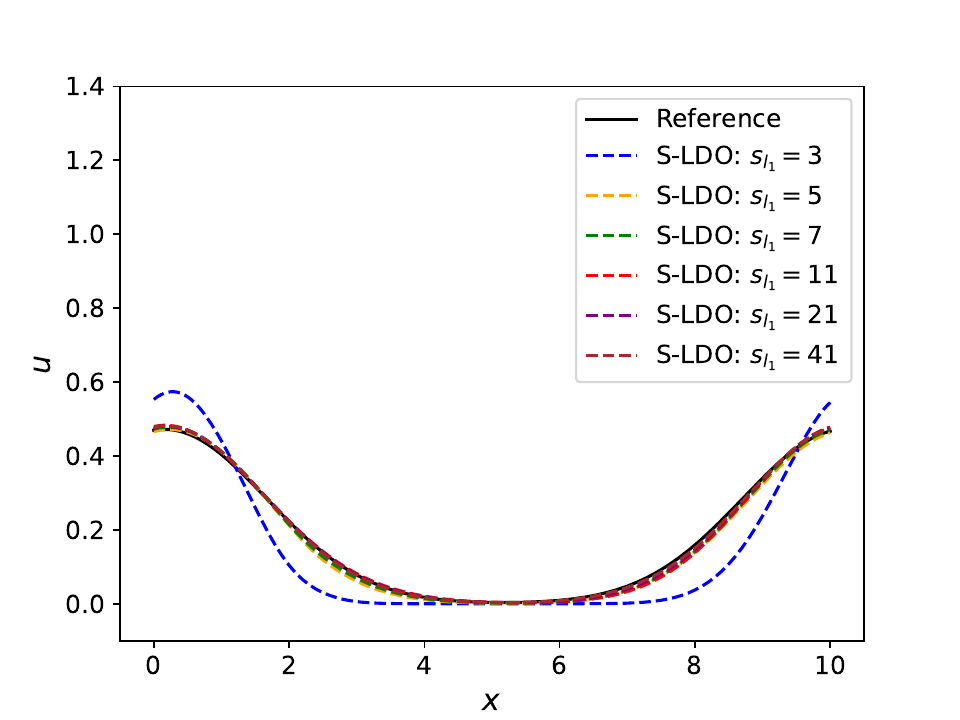}}\subfigure[\label{fig:error_both_SLDO}]{\includegraphics[width=0.33\textwidth, trim={0cm 0cm 1.5cm 0.5cm},clip]{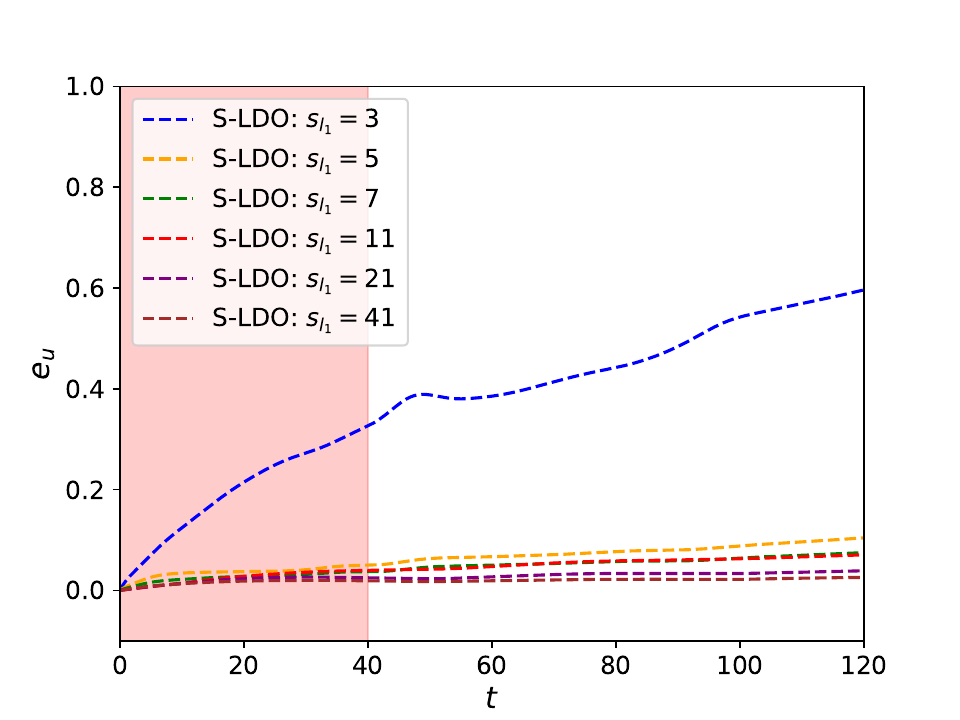}}
    \vspace{-3mm}
    \caption{Advection-diffusion problem ($c = 0.2$, $\nu = 0.02$): Predicted solutions at (a) $t = 20s$, (b) $t = 36s$ and (c) errors in time using LDOs (with $\beta_1 = 10^{-3}$). Predicted solutions at (d) $t = 20s$, (e) $t = 36s$ and (f) errors in time using SLDOs for several stencil sizes $s_{l_1}$ and with $s_{l_2} = 3$. In the error plots, the unshaded region is the region of extrapolation. }
    \label{fig:u_both_comp}
\end{figure}

The predicted solutions obtained using LDOs and S-LDOs, along with the corresponding errors, are compared to the reference data in \figref{u_both_comp} for different stencil sizes of the advection term. LDOs exhibit incorrect solution predictions for most stencil sizes at $t = 20s$. Smaller stencil sizes lead to a higher error in prediction and exhibit wildly unstable solution time instances even at $t = 36s$, which is within the dataset used for learning these operators. This error is reduced with increased stencil sizes and LDOs using the largest solution stencil, which exhibits good accuracy for $t = 20s$ and $36s$. This behavior is also reflected by a large growth of errors in solutions predicted by LDOs for $s_{l_1} \leq 21$. In contrast, S-LDOs exhibit a much higher accuracy where all but the smallest stencil size $s_{l_1} = 3$ yield close predictions to the reference result. This behavior holds for both time instances, indicating that S-LDOs are robust and less prone to stability issues. The temporal variation of error confirms this behavior as S-LDOs exhibit a very low error for all the stencil sizes except the smallest one. As errors remain low even for time instances not used to determine the operators, these results demonstrate the applicability of S-LDOs for long-time dynamics forecasting. 

\begin{figure}
    \centering
    \subfigure[\label{fig:Materror_both_LDO_beta1}]{\includegraphics[width=0.28\textwidth, trim={1.0cm 0.0cm 4.0cm 0cm},clip]{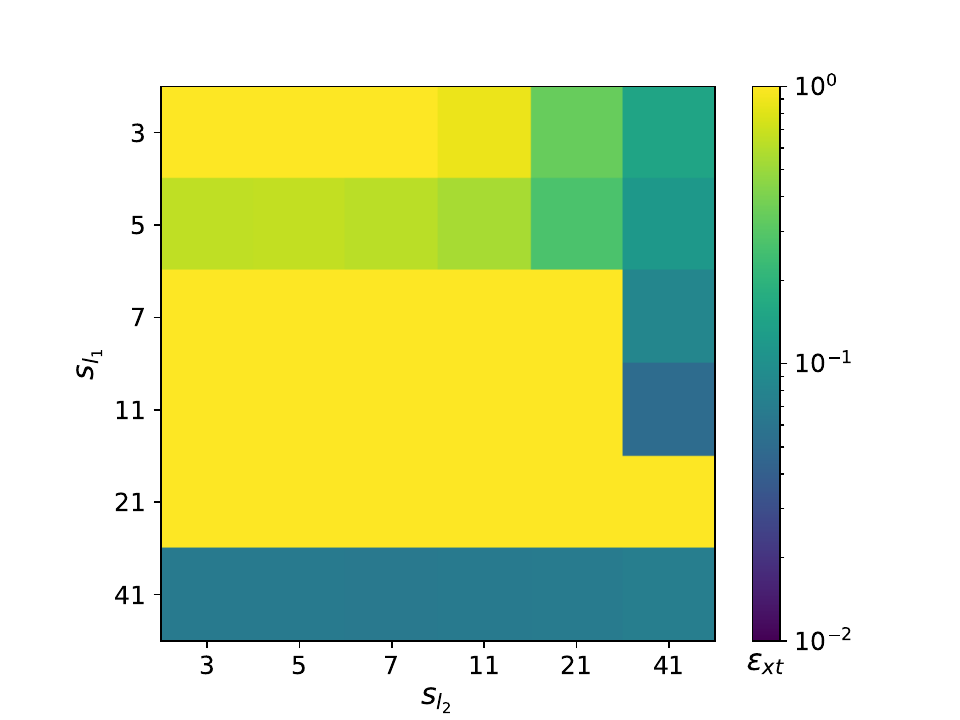}}
    \subfigure[\label{fig:Materror_both_LDO_beta2}]{\includegraphics[width=0.28\textwidth, trim={1.0cm 0.0cm 4.0cm 0cm},clip]{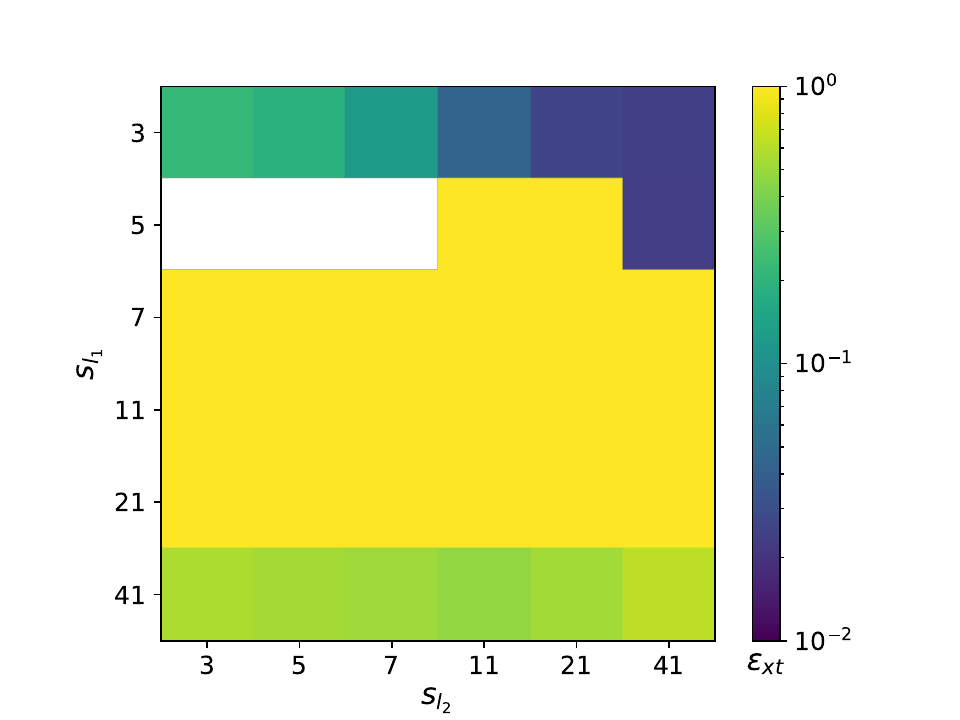}}    
    \subfigure[\label{fig:Materror_both_SLDO}]{\includegraphics[width=0.28\textwidth, trim={1.0cm 0.0cm 4.0cm 0cm},clip]{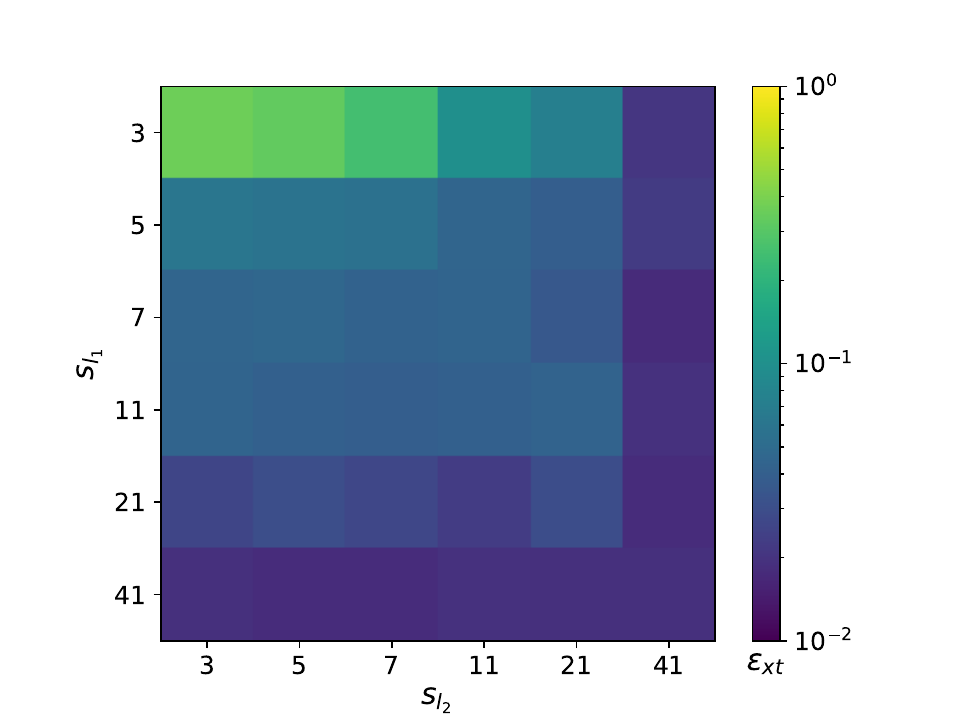}}    \includegraphics[width=0.075\textwidth, trim={12.5cm 0.5cm 1cm 0cm},clip]{errorMat_SLDO_both.pdf}  
    \vspace{-3mm}
    \caption{Advection-diffusion problem ($c = 0.2$, $\nu = 0.02$): Total error in the solution obtained using (a) LDOs (with $\beta_1 = 10^{-3}$), (b) LDOs (with $\beta_1 = 10^{-6}$) and (c) S-LDOs. The white regions indicate a very high error that is an undefined number.}
    \label{fig:Materror_both_comp}
\end{figure}

In the results until now, the stencil size of the diffusion operator was fixed at $s_{l_2} = 3$. To obtain a holistic picture of the operator performance for several combinations of advection and diffusion stencils, we compare the total error in both space and time
\begin{equation}
    \varepsilon_{xt} = \frac{\vert \vert \bm{u} (\cdot,\cdot) - \bm{u}^m (\cdot,\cdot) \vert \vert_F }{\vert \vert \bm{u} (\cdot,\cdot) \vert \vert_F},
    \label{eq:def_totalerror}
\end{equation}
where $\bm{u}$ is the reference solution and $\bm{u}^m$ is the solution predicted using LDOs or S-LDOs. This error metric is obtained by taking the $L_2$ norm in space and time, which corresponds to the Frobenius norm $\vert\vert \cdot \vert \vert_F$ if the solution is stored as a 2-D matrix. The total error for different advection and diffusion term stencil sizes is shown in \figref{Materror_both_comp}. We observe a high error for LDOs for most combinations of stencil sizes except when $s_{l_1} = 41$ or $s_{l_2} = 41$. These results indicate that LDOs require larger stencils for low error. Even in such scenarios, the system can become unstable and yield high errors for future time instances. This behavior is also observed at different regularization parameters as shown in the figure. Alternately, S-LDOs exhibit a very low error for most combinations of stencil sizes except when $s_{l_1} = 3$ and $s_{l_2} \leq 7$. The accuracy appears to increase with both advection and diffusion stencil sizes. The accuracy of S-LDOs for small stencil sizes can be further improved by adjusting the initial guess or selecting a lower tolerance for the optimizer. In this study, we kept this strategy the same for all stencils for a fairer comparison. These results indicate that S-LDOs perform much better than LDOs and highlight the importance of ensuring stability while learning differential operators to maintain stable and accurate results for dynamics forecasting. 

\subsection{1-D Burgers equation}

Having demonstrated the applicability of learned differential operators for linear PDEs, we now assess the validity of this approach for nonlinear PDEs. Therefore, we consider the 1-D Burgers equation in \eref{BurgerEq}. Unlike the scalar advection-diffusion problem, this equation has a nonlinear transport term. Using appropriate spatial discretization schemes, we obtain the semi-discrete form of the equations
\begin{equation}
\frac{d \bm{u}}{d t} + \bm{N} \bm{z} (\bm{u}) - \nu \bm{L} \bm{u} = 0,
\end{equation}
where $\bm{N}$ is the nonlinear operator, $\bm{L}$ is the linear diffusion operator and $\bm{z} (\bm{u})$ is composed of quadratic products of $\bm{u}$. We generate the data by using the $2^{nd}$-order centered difference for both nonlinear and linear diffusion terms with $129$ degrees of freedom. The differential equation for the $i^{th}$ degree of freedom is
\begin{equation}
\frac{d u_i}{d t} + (\bm{N}^i)^T \bm{z}_{\Omega^n_i} - \nu (\bm{L}^i)^T \bm{u}_{\Omega^l_i} = 0,
\label{eq:Burgers_eqFull}
\end{equation}
where
\begin{equation}
    \bm{N}^i = \frac{1}{2 \Delta x}[-1, 0, 1]^T \quad \text{and} \quad \bm{L}^i = \frac{1}{(\Delta x)^2}[1, -2, 1]^T. 
\end{equation}
The stencil for the nonlinear term is $\bm{z}_{\Omega^n_i} = u_i [u_{i-1}, u_i, u_{i+1}]^T$, whereas it is $\bm{u}_{\Omega^l_i} = [u_{i-1}, u_i, u_{i+1}]^T$ for the linear diffusion term. In this article, we model the semi-discrete form as
\begin{equation}
\frac{d u_i}{d t} + (\bm{N}^{i,m})^T \bm{z}_{\Omega^n_i} - \nu (\bm{L}^{i,m})^T \bm{u}_{\Omega^l_i} = 0,
\label{eq:BurgerDForm_NL}
\end{equation}
where $\bm{N}^{i,m} \in \mathbb{R}^{s_{l_1}}$ and $\bm{L}^{i,m} \in \mathbb{R}^{s_{l_2}}$ are nonlinear and linear operators determined from the data, while $\bm{z}_{\Omega^n_i} \in \mathbb{R}^{s_{l_1}}$ and $\bm{u}_{\Omega^l_i} \in \mathbb{R}^{s_{l_2}}$ are the solution stencils of dimensionality $s_{l_1}$ and $s_{l_2}$. For nonlinear problems such as the 1-D Burgers equation, a single stable equilibrium point does not exist as any constant solution is a valid equilibrium point. Due to the viscous term, the solution decays to a stable equilibrium close to the initial condition. As the modeled semi-discrete form in \eref{BurgerDForm_NL} is nonlinear, we cannot guarantee the global stability of these equations. Instead, we focus on local linear stability around the equilibrium point $\bm{u}^0$ by considering the $1^{st}$-order Taylor series approximation to the equation and obtaining the linearized equation
\begin{equation}
\frac{d u_i}{d t} + (\bm{N}^{i,m,L})^T \bm{u}_{\Omega^n_i} = 0,
\label{eq:BurgerDForm_NL_TS}
\end{equation}
with
\begin{equation}
    \bm{N}^{i,m,L} =  \Bigg\vert\frac{\partial (\bm{N}^{i,m})^T \bm{z}_{\Omega^n_i}}  {\partial \bm{u}_{\Omega^n_i}} - \nu \frac{\partial (\bm{L}^{i,m})^T \bm{u}_{\Omega^l_i}}  {\partial \bm{u}_{\Omega^n_i}}  \Bigg\vert_{\bm{u}^0},
\end{equation}
where $\vert \cdot \vert_{\bm{u}^0}$ implies the term is evaluated at equilibrium point $\bm{u}^0$ and $\bm{N}^{i,m,L}$ is the local linearized operator that is assembled to give the global linearized operator $\bm{N}^{m,L}$. This linearized equation is used in \eref{ConstrainedFinalNonLin} to determine linear stability constraints for the regression problem in \eref{RegProbConstrainedFinalNonLin} for obtaining S-LDOs. 

The reference data is used to learn differential operators and assess the performance of LDOs and S-LDOs. This data is obtained by integrating \eref{Burgers_eqFull} in time using a $1^{st}$-order explicit time integration scheme and a time step size of $0.002$.
This test case is initialized with Gaussian random perturbations with a mean of $0.3$ and a standard deviation of $0.2$. For this test case, we can assess the eigenvalues of the operators $\bm{N}^{m,L}$, which can provide us with the notion of linear stability of the system around an equilibrium point. This analysis aligns with one for the advection-diffusion equation, which showed that S-LDOs are linearly stable, whereas LDOs are not. A detailed discussion on this aspect is excluded from this article for brevity. Instead, we discuss the consequences of this behavior on the evolution of the solution over time.

\begin{figure}
    \centering
    \includegraphics[width=0.5\textwidth]{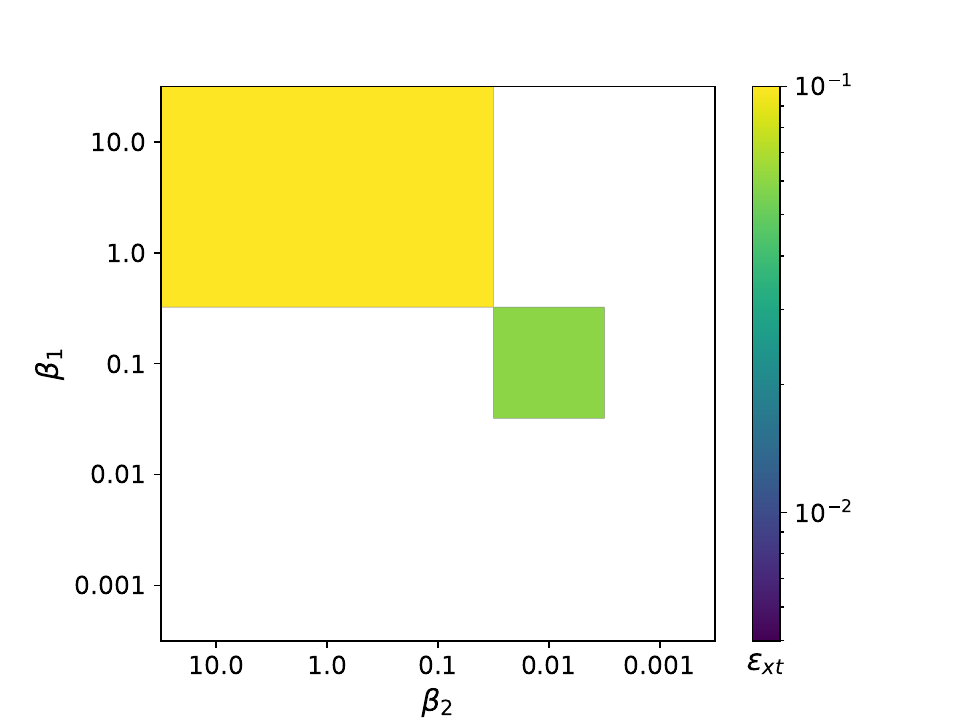}
    \caption{1-D Burgers problem: Total error in the solution obtained for LDOs with $s_{l_1} = 5$ and $s_{l_2} = 5$. The white regions indicate a very high error that is an undefined number.}
    \label{fig:Materror_LDO_Reg}
\end{figure}

We first assess the impact of regularization parameters $\beta_1$ and $\beta_2$ on the predicted solution for LDOs with a fixed stencil size of $s_{l_1} = s_{l_2} = 5$. The total error in the predicted solution for LDOs is shown in \figref{Materror_LDO_Reg}. The results indicate the sensitivity of LDOs to selected regularization parameters and typically give predictions with very high errors. As $\beta_1 = 0.1$ and $\beta_2 = 0.01$ lead to reasonable errors where solutions are well undefined, we select these regularization parameters for LDOs for further comparison.

\begin{figure}
    \centering
    \subfigure[\label{fig:uBurg_LDO_t100}]{\includegraphics[width=0.33\textwidth, trim={0.0cm 0cm 1.5cm 0.5cm},clip]{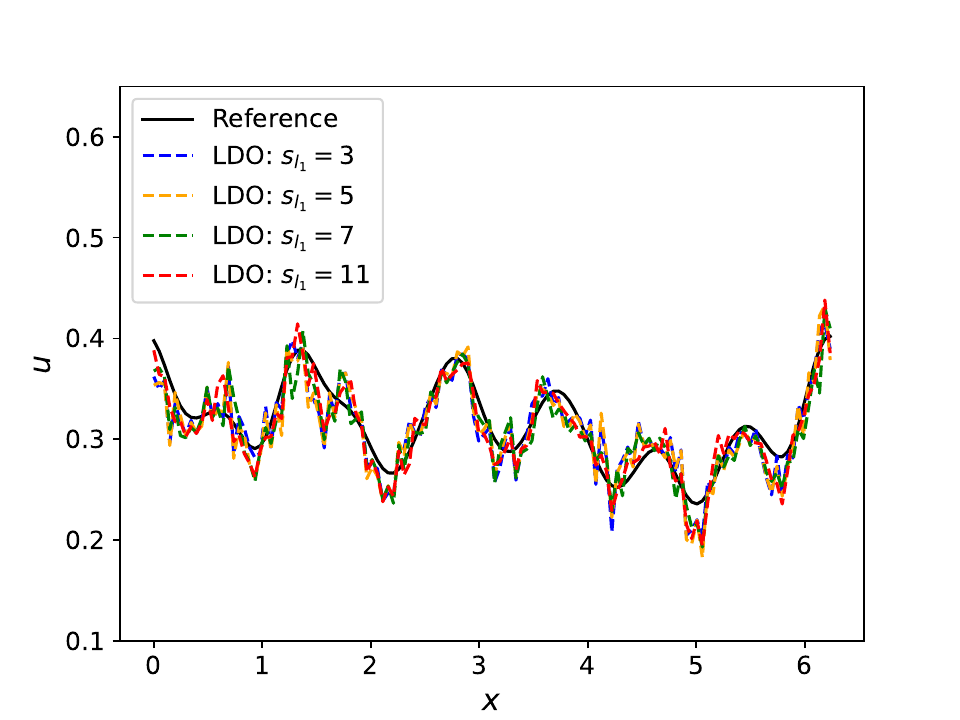}}\subfigure[\label{fig:uBurg_LDO_t1000}]{\includegraphics[width=0.33\textwidth, trim={0.0cm 0cm 1.5cm 0.5cm},clip]{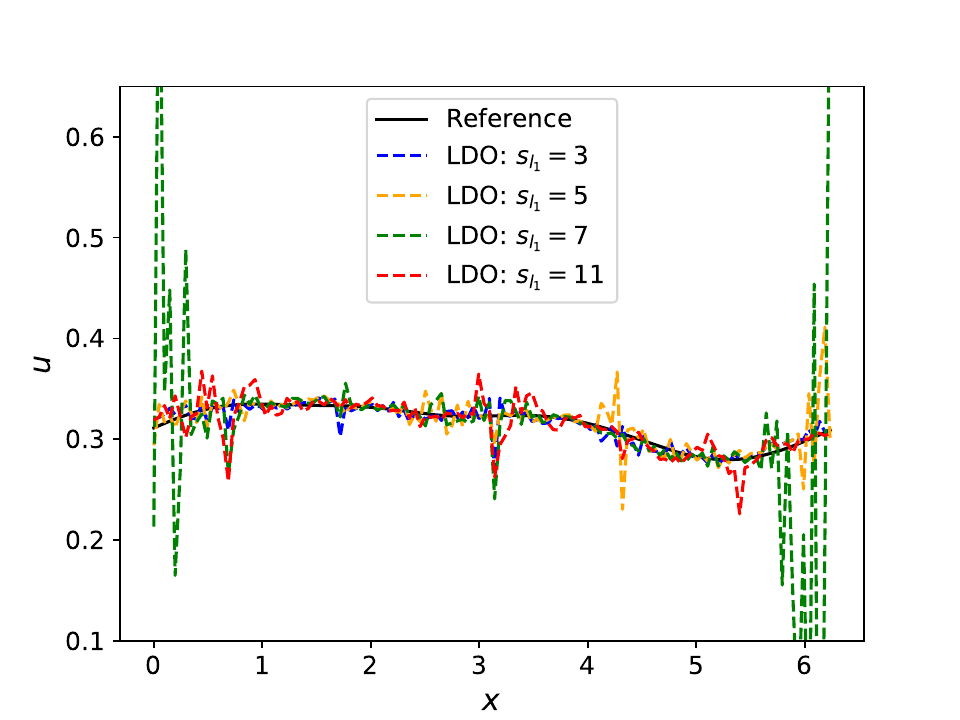}}\subfigure[\label{fig:errorBurg_LDO}]{\includegraphics[width=0.33\textwidth, trim={0.0cm 0cm 1.2cm 0.5cm},clip]{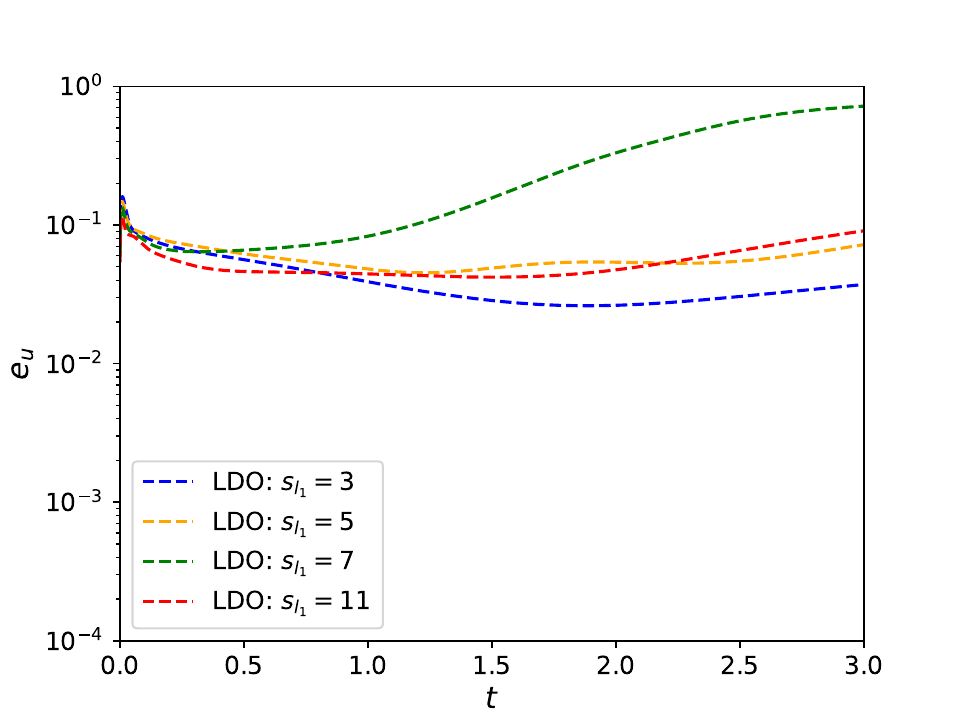}}
    
    \subfigure[\label{fig:uBurg_SLDO_t100}]{\includegraphics[width=0.33\textwidth, trim={0.0cm 0cm 1.5cm 0.5cm},clip]{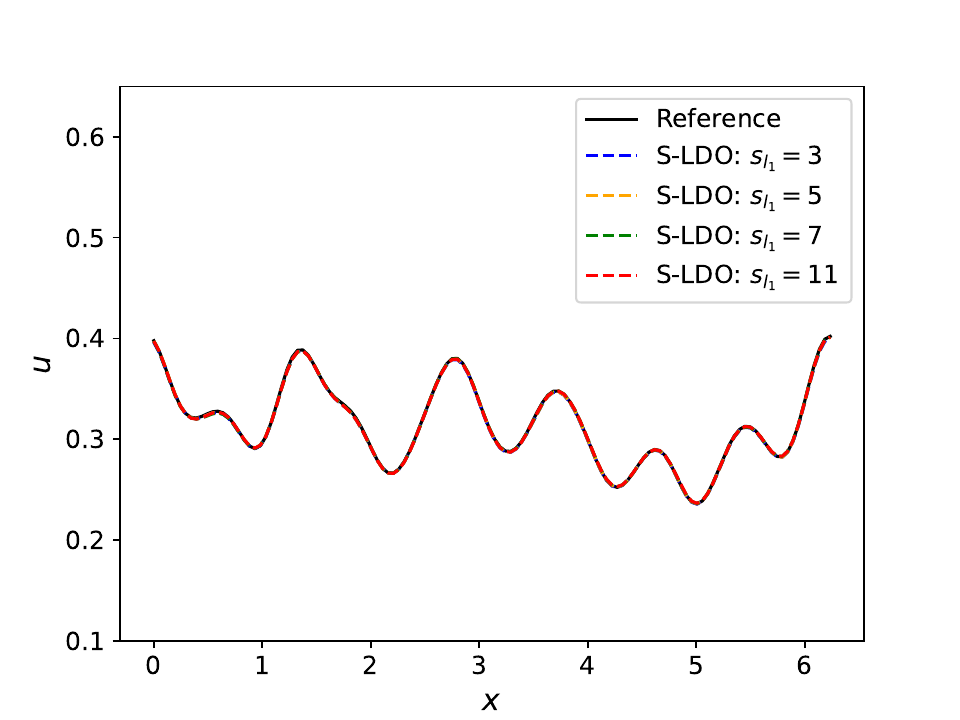}} \subfigure[\label{fig:uBurg_SLDO_t1000}]{\includegraphics[width=0.33\textwidth, trim={0.0cm 0cm 1.5cm 0.5cm},clip]{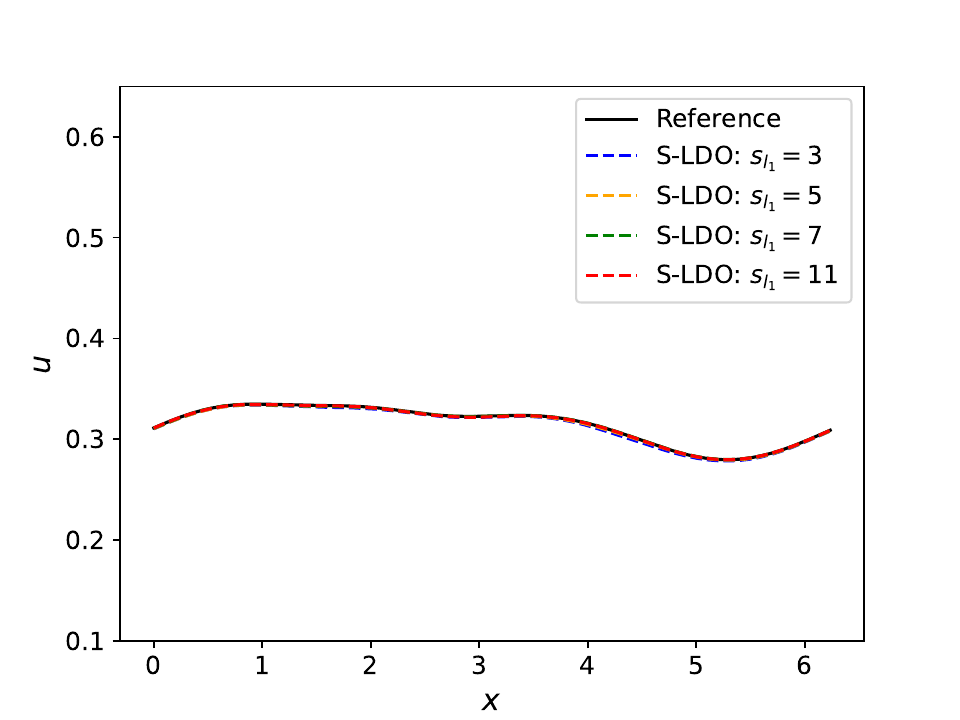}}\subfigure[\label{fig:errorBurg_SLDO}]{\includegraphics[width=0.33\textwidth, trim={0.0cm 0cm 1.2cm 0.5cm},clip]{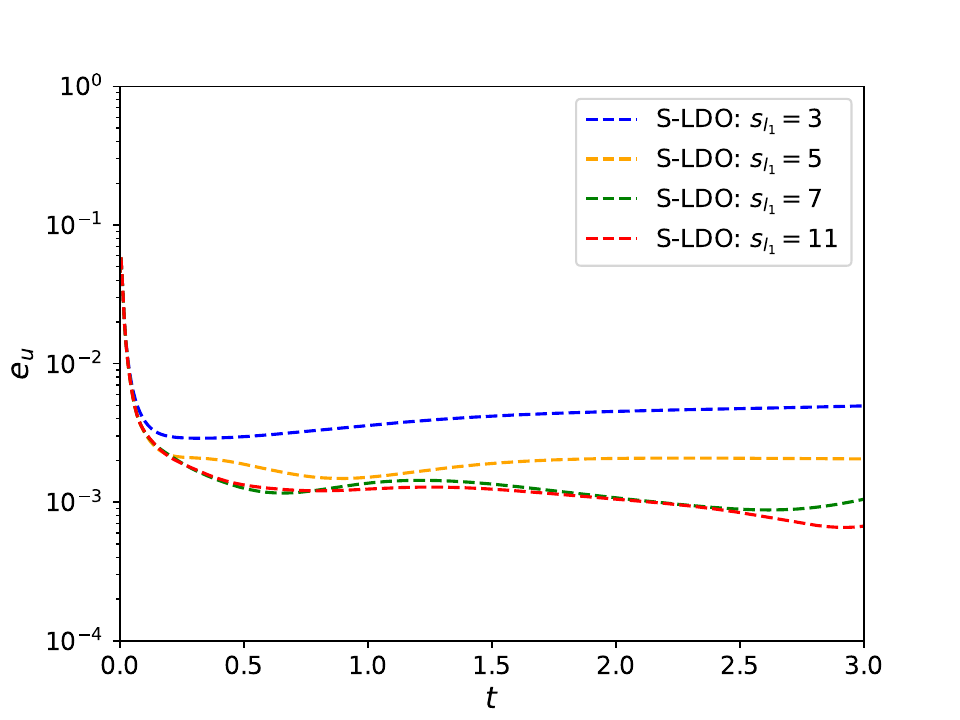}}
    \vspace{-3mm}
    \caption{1-D Burgers problem: Predicted solutions at (a) $t = 0.2 s$, (b) $t = 2 s$ and (c) corresponding errors for LDOs (with $\beta_1 = 0.1$ and $\beta_2 = 0.01$). Predicted solutions at (d) $t = 0.2 s$, (e) $t = 2 s$ and (f) corresponding errors for S-LDOs. These are evaluated for different $s_{l_1}$ while keeping $s_{l_2} = 5$.}
    \label{fig:uBurg_comp}
\end{figure}

\begin{figure}
    \centering
    \subfigure[\label{fig:Materror_both_LDO_Burger}]{\includegraphics[width=0.4\textwidth, trim={1.0cm 0.0cm 4.0cm 0cm},clip]{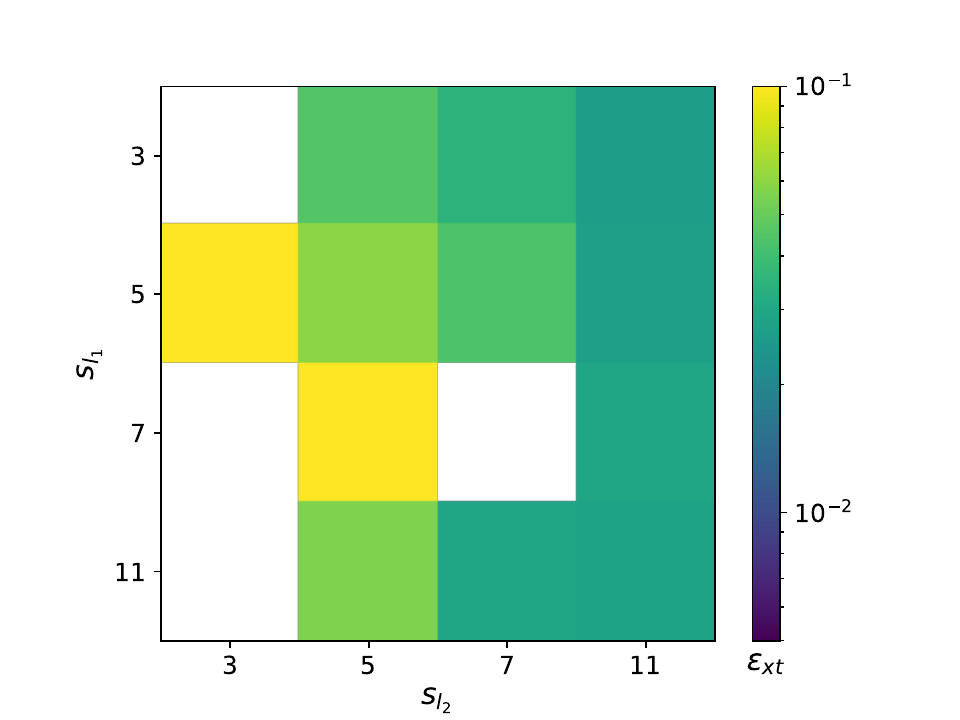}}    
    \subfigure[\label{fig:Materror_both_SLDO_Burger}]{\includegraphics[width=0.4\textwidth, trim={1.0cm 0.0cm 4.0cm 0cm},clip]{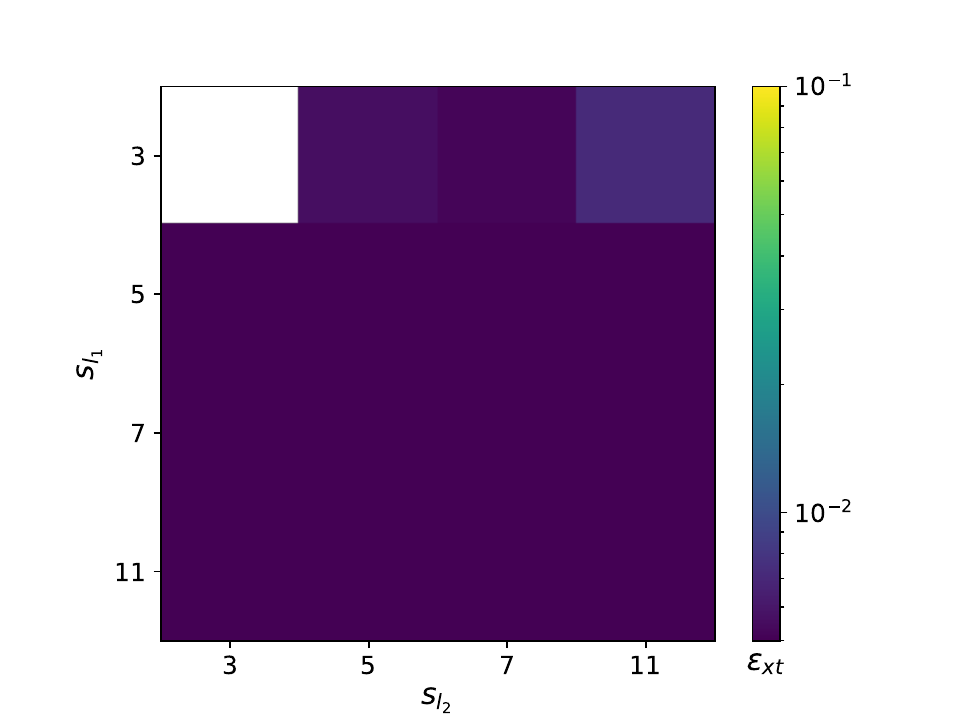}}
    \includegraphics[width=0.105\textwidth, trim={12.5cm 0.5cm 1cm 1.0cm},clip]{MaterrorUBurger_SLDO.pdf}  
    \vspace{-3mm}
    \caption{1-D Burgers problem: Total error in the solution obtained using (a) LDOs (with $\beta_1 = 0.1$ and $\beta_2 = 0.01$) and (b) S-LDOs. The white regions indicate a very high error that is an undefined number.}
    \label{fig:Materror_Burger_comp}
\end{figure}

The predicted solutions obtained using LDOs and S-LDOs and corresponding errors are compared to the reference data in \figref{uBurg_comp}. We observe that even at an early time $t = 0.2s$, LDOs for multiple stencil sizes exhibit oscillations. These oscillations persist and grow larger as evident from the results at $t = 2s$. These results correspond to the best selection of the stencil sizes for a given regularization parameter, and still, LDOs do not provide stable and accurate results. This behavior is also highlighted in the error plots, which indicate high errors for LDOs for different stencil sizes. On the contrary, S-LDOs for all stencil sizes exhibit high accuracy at both time instances as the results are close to the reference data. Furthermore, the errors in solution predicted by S-LDOs are consistently low for all stencil sizes, even at longer times. These results indicate that S-LDOs exhibit superior stability properties and provide feasible, stable and physically accurate solutions. 

The matrix plot with total errors for different combinations of stencil sizes for LDOs and S-LDOs is shown in \figref{Materror_Burger_comp}. We observe that LDOs exhibit high errors for all combinations of stencil sizes. The total error is lower for the largest stencil sizes, similar to the behavior observed for the advection-diffusion problem. However, this error is still much larger compared to errors exhibited by S-LDOs. For the smallest stencil size for S-LDOs, some instability in the results is observed, which gives high errors. As S-LDOs can only ensure linear stability, whereas the system is nonlinear, there is no guarantee that all stencil sizes will always have stable solutions. Unstable behavior is also observed at several stencil sizes larger than those considered in this study. The inability of differential operators to perform well in such scenarios is not a big issue as larger stencil sizes may not even be considered due to higher computational costs. This behavior highly depends on the linearization strategy and selection of the equilibrium point. The errors for most combinations of stencil sizes are very low, demonstrating the applicability of S-LDOs in giving accurate solutions. Although S-LDOs are designed to provide linearly stable operators, they perform remarkably well for nonlinear problems, as demonstrated using this test case. These results highlight that adding linear stability constraints while learning nonlinear operators can be a viable solution even for nonlinear PDEs. 

\begin{figure}
    \centering    
    \subfigure[\label{fig:uCont_t250_true}]{\includegraphics[width=0.29\textwidth,trim={0.5cm 0.5cm 1.0cm 0cm},clip]{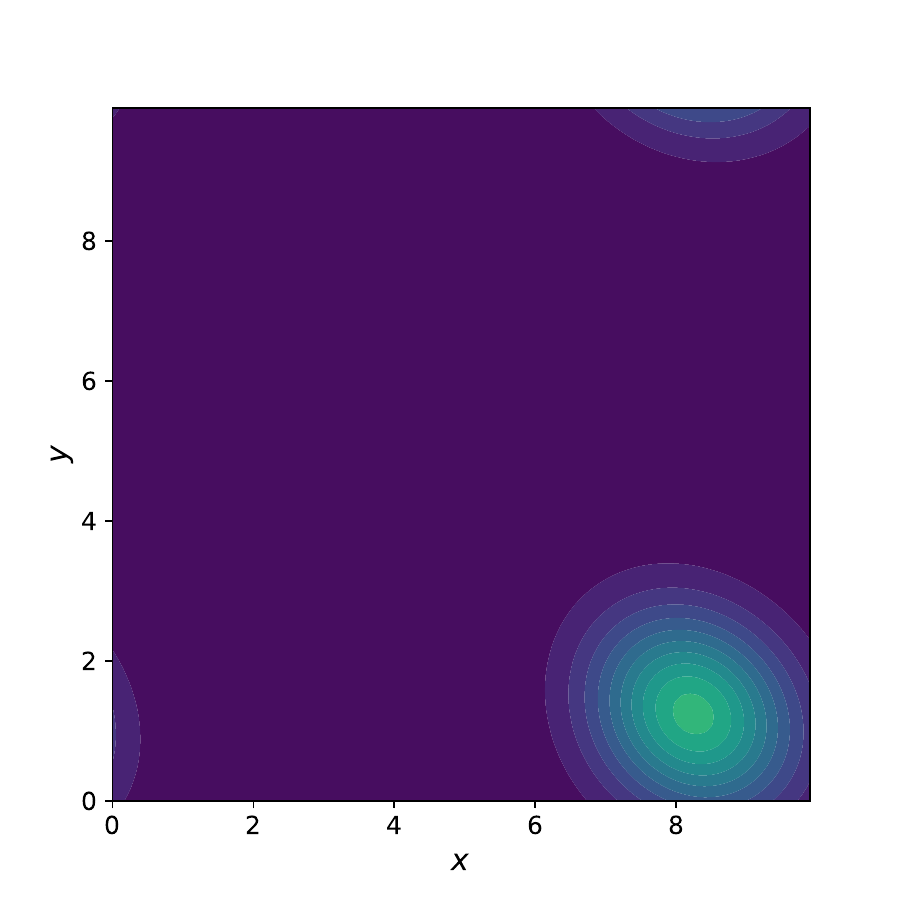}}    
    \subfigure[\label{fig:uCont_t250_LDO}]{\includegraphics[width=0.29\textwidth,trim={0.5cm 0.5cm 1.0cm 0cm},clip]{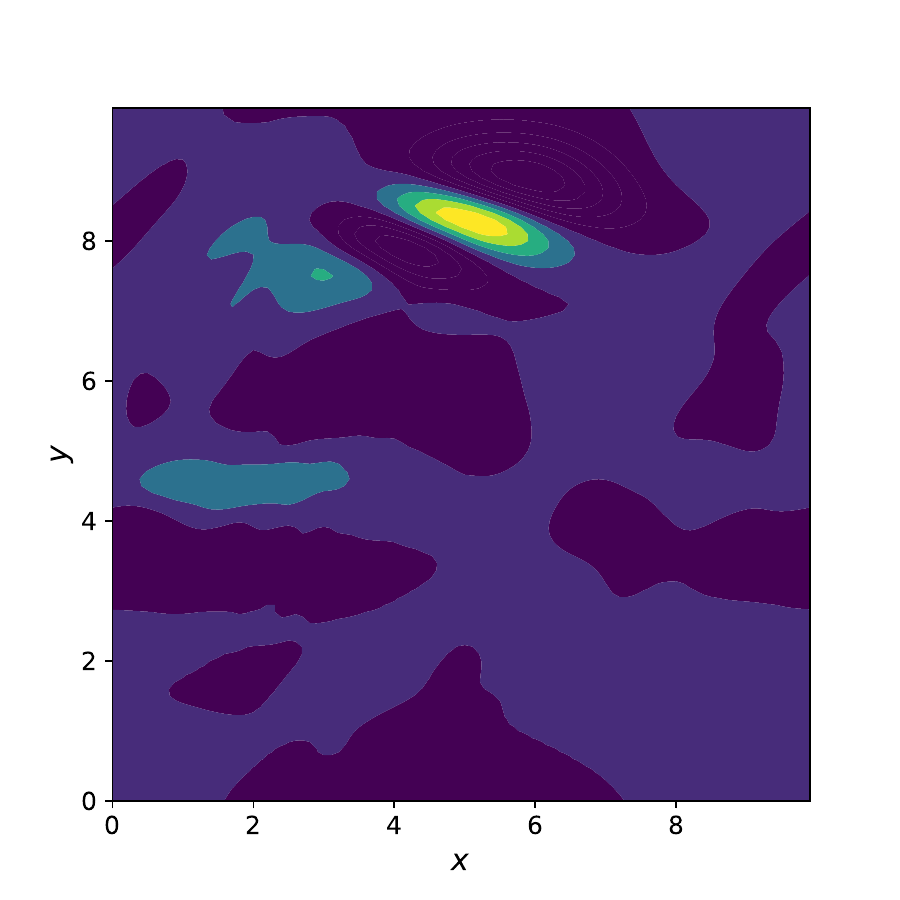}}    
    \subfigure[\label{fig:uCont_t250_SLDO}]{\includegraphics[width=0.29\textwidth,trim={0.5cm 0.5cm 1.0cm 0cm},clip]{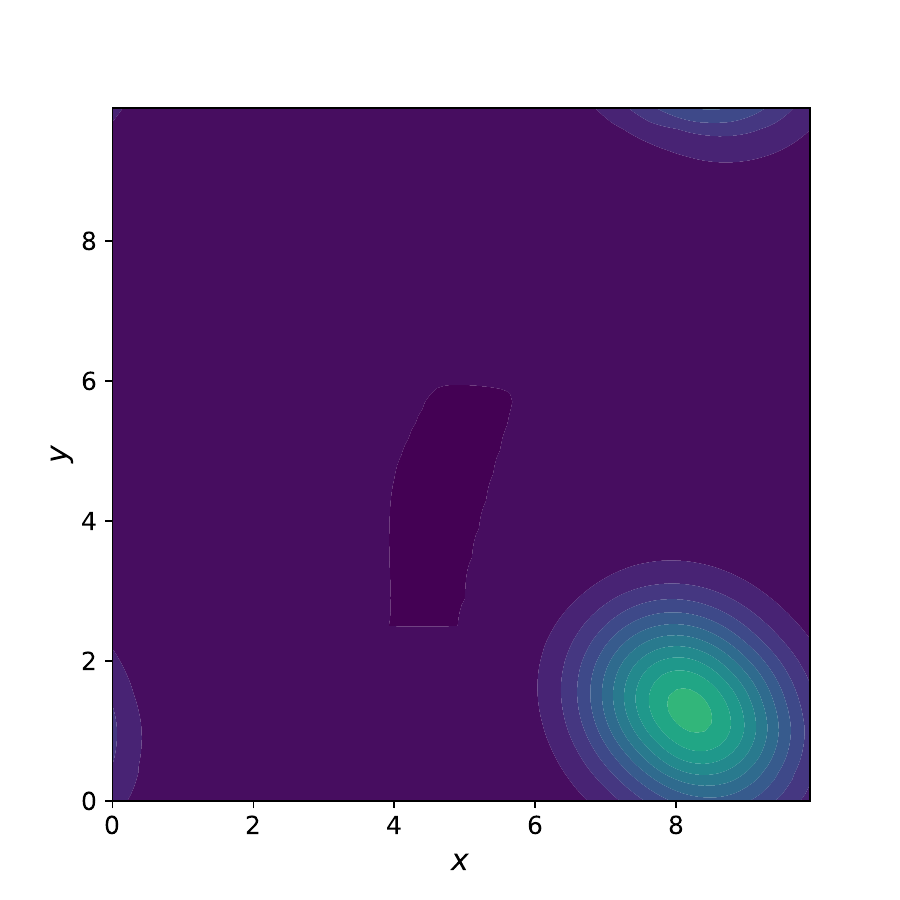}} \includegraphics[width=0.083\textwidth, trim={11.5cm 0.5cm 0 0cm},clip]{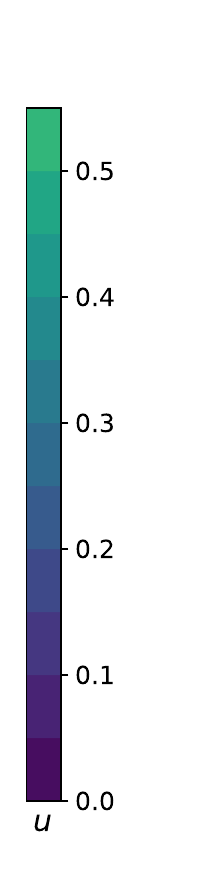}  
    \vspace{-3mm}
    \caption{2-D advection problem: Solutions at $t = 12.5s$ for (a) reference data, (b) LDO with $\beta = 10^{-2}$ and (c) S-LDO for stencil sizes $s_{l_x} = 7$ and $s_{l_y} = 7$.}
    \label{fig:uCont_t250}
\end{figure}
\subsection{2-D scalar advection equation}

We select a 2-D scalar advection equation for the third test case to demonstrate the applicability of learned differential operators to operators for 2-D PDEs. The equation is given as follows:
\begin{equation}
\frac{\partial u}{\p t} + \bm{c} \cdot \nabla u = 0,
\end{equation}
where $\bm{c}$ is the advection velocity. Using an appropriate spatial discretization scheme, the PDE can be converted to a set of ODEs
\begin{equation}
\frac{d \bm{u}}{d t} + c_x \bm{L}^x \bm{u} + c_y \bm{L}^y \bm{u} = 0,
\label{eq:ODE_2DAdvec}
\end{equation}
where $\bm{L}^x$ and $\bm{L}^y$ are the linear operators on the solution components. For data generation, we use a $1^{st}$-order backward difference for approximating these operators with $101$ degrees of freedom in each direction. This selection implies that the differential equations for the $i^{th}$ degree of freedom is
\begin{equation}
\frac{d u_i}{d t} + c_x (\bm{L}^{i,x})^T \bm{u}_{\Omega^x_{i}} + c_y (\bm{L}^{i,y})^T \bm{u}_{\Omega^y_{i}} = 0,
\label{eq:ODE_2DAdvec_local}
\end{equation}
where 
\begin{equation}
    \bm{L}^{i,x} = \frac{1}{\Delta x} [-1, 1, 0]^T \quad \text{and} \quad \bm{L}^{i,y} = \frac{1}{\Delta y} [-1, 1, 0]^T.
\end{equation}

The solution stencil for the $\bm{L}^{i,x}$ operator is denoted by $\bm{u}_{\Omega^x_{i}}$, whereas the stencil for the $\bm{L}^{i,y}$ operator is denoted by $\bm{u}_{\Omega^{y}_{i}}$. This test case is initialized with the initial condition
\begin{equation}
    u(x,y,0) = \exp{\Big(-(x-2)^2 + (y-5)^2\Big)}
\end{equation}
and integrated in time using a $1^{st}$-order forward Euler method with a timestep of $0.05$. The generated data is used to learn the differential operators and as a reference result for assessing the performance of LDOs and S-LDOs. The advection velocity is chosen to be $c = 0.5 \hat{e}_1 + 0.5 \hat{e}_2$ where $\hat{e}_1$ and $\hat{e}_2$ are unit vectors aligned to Cartesian $x$ and $y$ directions respectively. In this article, we determine the unknown operators $\bm{L}^{i,x,m}$ and $\bm{L}^{i,y,m}$ from data. Therefore, \eref{ODE_2DAdvec_local} is modeled as
\begin{equation}
\frac{d u_i}{d t} + c_x (\bm{L}^{i,x,m})^T \bm{u}_{\Omega^x_{i}} + c_y (\bm{L}^{i,y,m})^T \bm{u}_{\Omega^y_{i}} = 0,
\label{eq:ODE_2DAdvec_model}
\end{equation}
where $\bm{L}^{i,x,m} \in \mathbb{R}^{s_{l_x}}$ and $\bm{L}^{i,y,m} \in \mathbb{R}^{s_{l_y}}$ are the modeled linear operators. The stencils for these operators $\bm{u}_{\Omega^x_{i}} \in \mathbb{R}^{s_{l_x}}$ and $\bm{u}_{\Omega^y_{i}} \in \mathbb{R}^{s_{l_y}}$ have a dimensionality of $s_{l_x}$ and $s_{l_y}$ respectively. The constraints in S-LDO formulation ensure a positive real part for the eigenvalues of $c_x \bm{L}^m_x + c_y \bm{L}^m_y$. A detailed analysis showing this behavior is excluded for brevity, although the behavior is similar to the one observed for the 1-D advection equation. 

\begin{figure}
    \centering
    \subfigure[\label{fig:error2DAdvec_LDO_slx}]{\includegraphics[width=0.49\textwidth, trim={0.0cm 0cm 1.5cm 0.5cm},clip]{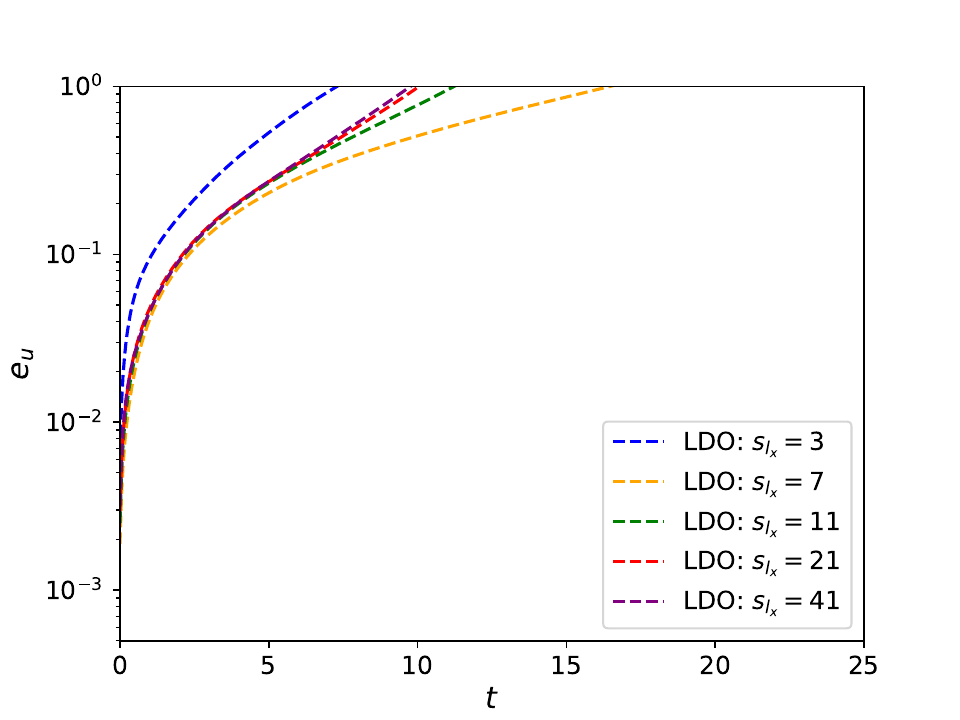}}
    \subfigure[\label{fig:error2DAdvec_SLDO_slx}]{\includegraphics[width=0.49\textwidth, trim={0.0cm 0cm 1.5cm 0.5cm},clip]{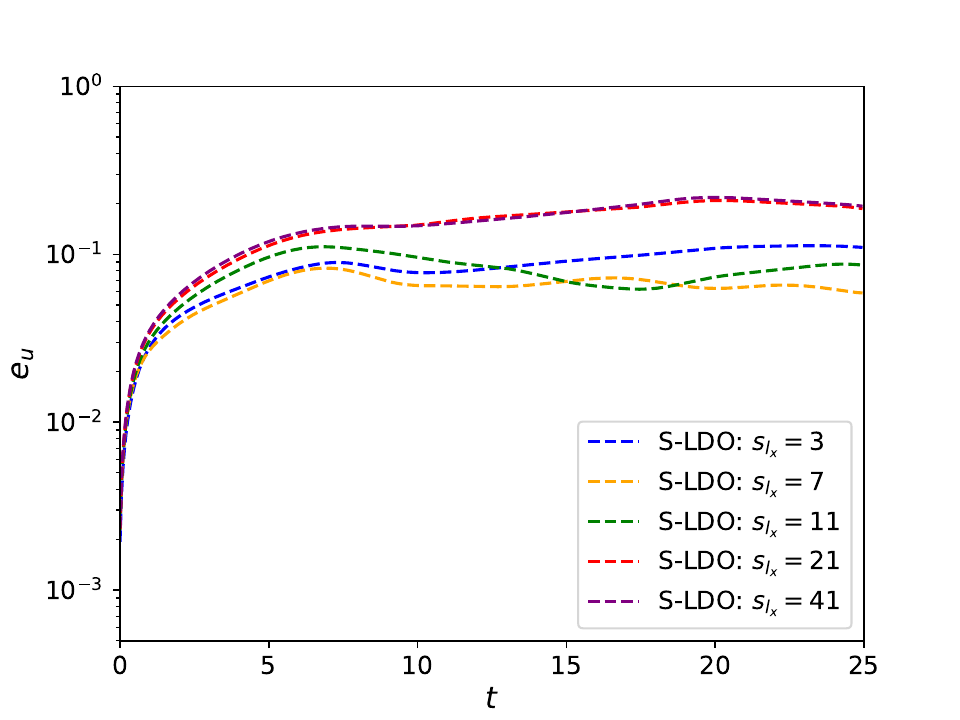}}
    \vspace{-3mm}
    \caption{2-D advection problem: Temporal variation of errors in the predicted solution obtained using (a) LDOs (with $\beta_1 = 10^{-2}$) and (b) S-LDOs for several stencil sizes $s_{l_x}$ while keeping $s_{l_y} = 7$.}
    \label{fig:error2DAdvec_comp_slx}
\end{figure}

The solution predicted using an LDO and an S-LDO is compared to the reference data in \figref{uCont_t250}. We observe that prediction by LDO is very different than the reference data as the unstable behavior of these operators leads to large deviations in the results. The predicted solution by the S-LDO exhibits a similar behavior and magnitude as the reference data. This behavior echoes other test cases, where S-LDOs consistently gave much more accurate results than LDOs. The temporal variation of errors in the solutions predicted by LDOs and S-LDOs for different stencil sizes $s_{l_x}$ is shown in \figref{error2DAdvec_comp_slx}. We observe that LDOs exhibit a large error and grow quickly in time for different stencil sizes. These results confirm that LDOs are ill-equipped to infer operators for advection-dominated problems, as already assessed from the 1-D advection case. Conversely, S-LDOs exhibit a low error for all the stencil sizes considered in this study. After an initial rise, the error reaches a near-constant value, which does not change drastically as time increases. This behavior also numerically verifies the stability properties of S-LDOs and demonstrates its applicability for stable prediction of solution field. Although these results are shown for different values of stencil sizes $s_{l_x}$ while keeping $s_{l_y} = 7$, this behavior also holds other values of $s_{l_y}$.  

\begin{figure}
    \centering
    \subfigure[\label{fig:Materror2DAdvec_LDO}]{\includegraphics[width=0.42\textwidth,trim={1.0cm 0.0cm 4.0cm 0cm},clip]{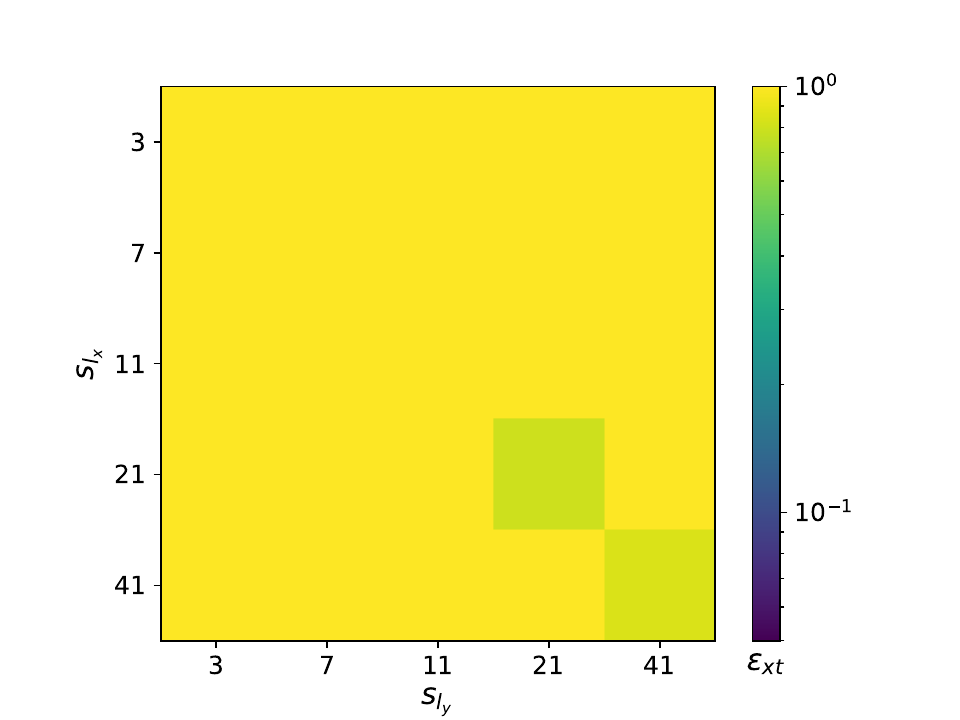}}        \subfigure[\label{fig:Materror2DAdvec_SLDO}]{\includegraphics[width=0.42\textwidth,trim={1.0cm 0.0cm 4.0cm 0cm},clip]{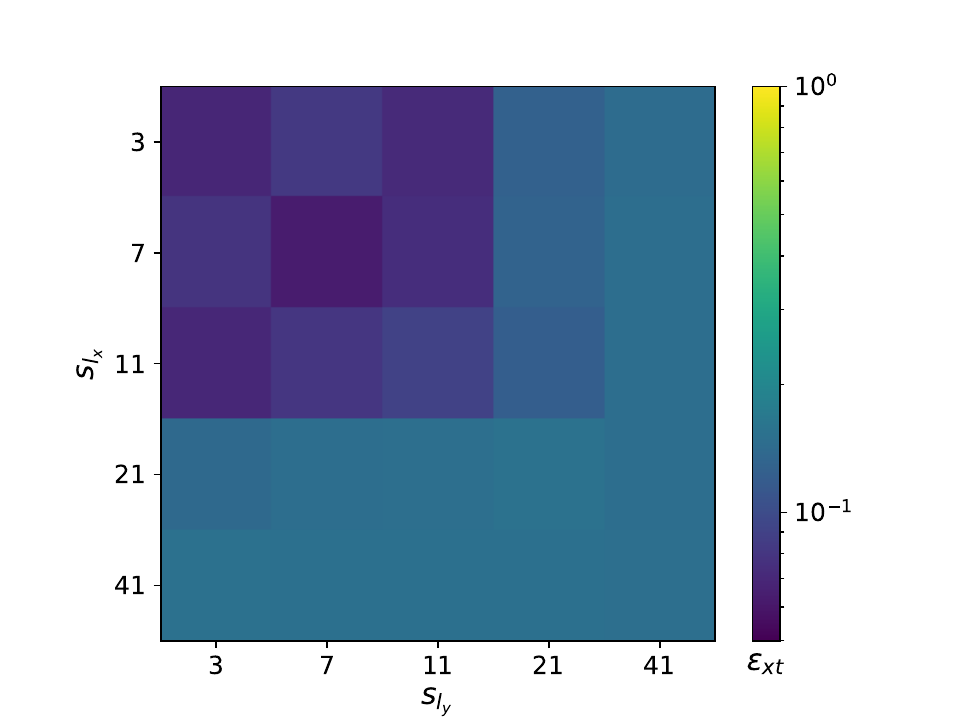}}    \includegraphics[width=0.11\textwidth, trim={12.5cm 0.5cm 1cm 0cm},clip]{MatErrorUComp_SLDO_Ieq_2.pdf}  
    \vspace{-3mm}
    \caption{2-D advection problem: Total error in the solution obtained using (a) LDOs (with $\beta_1 = 10^{-2}$) and (b) S-LDOs for different stencil sizes $s_{l_x}$ and $s_{l_y}$.}
    \label{fig:Materror_2DAdvec_comp}
\end{figure}

The total error in the predicted solution for LDOs and S-LDOs of different stencil sizes is shown in \figref{Materror_2DAdvec_comp}. We observe that the total error is very high for LDOs for all combinations of stencil sizes. This behavior renders LDOs impractical for their use in dynamics predictions. On the contrary, S-LDOs exhibit a comparatively low error for all the stencil sizes. The error appears to be lower for the smaller stencil sizes and becomes higher for larger stencils. This behavior is different from those observed for the 1-D advection case, where the error appeared to decrease with a larger stencil. Nevertheless, the low prediction error and theoretical linear stability make S-LDOs suitable for dynamics forecasting. 

\section{Conclusions}

Several applications, such as system identification and nonintrusive reduced order modeling, motivate the need to identify discrete approximations of PDEs from data. Most studies have focused on modern machine learning techniques to obtain accurate approximations of a PDE. However, these approaches often yield noninterpretable and dense representations, which may not be scalable for large-scale applications. Even fewer studies have addressed stability concerns for learned semi-discrete differential equations.

In this article, we propose a novel methodology for determining sparse semi-discrete approximations of PDEs from data while ensuring the linear stability of learned approximations. This approach is inspired by common spatial discretizations that have a sparse structure for computational efficiency and are stable for accurate dynamics prediction. We demonstrate that the standard regression approach, which is common in the literature, does not yield theoretically stable differential operators even for simple 1-D PDEs. We overcome this drawback by identifying stability conditions on local differential operators using dynamical system theory. These conditions are then added to the regression problem as constraints on the learned differential operators. Solving these constrained regression problems yields theoretically stable differential operators for linear PDEs. We also extend this approach to nonlinear PDEs by formulating constraints using linearized differential equations. The applicability of the proposed approach is demonstrated using three examples: 1-D scalar advection-diffusion equation, 1-D Burgers equation and 2-D advection equation. The numerical experiments indicated that differential operators learned using the standard regression approach yielded unstable differential operators for different combinations of stencil sizes and regularization parameters. Consequently, the learned differential operator often produces highly oscillatory results in the initial time window and blows up in finite time, even for linear problems. In contrast, the learned differential operators obtained using the proposed constrained regression approach yielded stable and highly accurate results for linear and nonlinear PDEs considered in this study.

The proposed approach to incorporate constraints for ensuring stability while learning differential operators can be extended to several practical applications. These learned differential operators are ideally suited for developing nonintrusive reduced order models. We plan to demonstrate the proposed approach for this application while comparing it with other nonintrusive approaches \cite{Audouze2013, Hesthaven2018, Peherstorfer2016}. A key benefit of the proposed approach is the possibility of enabling nonintrusive reduced order modeling for standard projection-based techniques such as those that rely on stabilization \cite{Carlberg2011} and closure modeling \cite{Ahmed2021, Prakash2024}. We also plan to extend this approach to determine stable coarse-grained discretizations from data, enabling high-fidelity simulations at a lower computational cost. 

\section{Acknowledgements}

The authors would like to acknowledge the support from the National Science Foundation (NSF) grant, CMMI-1953323, for the funds used towards this project. The research in this paper was also sponsored by the Army Research Laboratory and was accomplished under Cooperative Agreement Number W911NF-20-2-0175. The views and conclusions contained in this document are those of the authors and should not be interpreted as representing the official policies, either expressed or implied, of the Army Research Laboratory or the U.S. Government. The U.S. Government is authorized to reproduce and distribute reprints for Government purposes notwithstanding any copyright notation herein.

\bibliographystyle{unsrt}
\bibliography{main_bbl.bib}

\end{document}